\documentclass[a4paper,11pt,reqno]{amsart}
\usepackage{amssymb}
\usepackage{amscd}
\usepackage{amsthm}
\usepackage{amsmath}
\usepackage[all]{xy}
\usepackage{extarrows}
\usepackage{rotating}
\usepackage{tikz-cd}
\usepackage{leftidx}
\usepackage{bm}
\usepackage{dynkin-diagrams}
\usepackage{mathrsfs}
\usepackage{rotating}
\usepackage{upgreek}
\usepackage{comment}
\numberwithin{equation}{section}
\usepackage[pagebackref=true]{hyperref}
\renewcommand*{\backref}[1]{}
\renewcommand*{\backrefalt}[4]{[{	%
		\ifcase #1 Not cited.%
		\or Cited on page~#2.%
		\else Cited on pages #2.%
		\fi%
	}]}

\usepackage{aliascnt}
\usepackage{cleveref}
\newcommand{\fd}{\mathsf{fd}}

\renewenvironment{proof}{{ \textbf{Proof}.}}{\qed}
\newtheorem{Thm}{Theorem}[section]
\newaliascnt{Lem}{Thm}
\newtheorem{Lem}[Lem]{Lemma}
\aliascntresetthe{Lem}

\newaliascnt{Def}{Thm}
\newtheorem{Def}[Def]{Definition}
\aliascntresetthe{Def}

\newaliascnt{Cor}{Thm}
\newtheorem{Cor}[Cor]{Corollary}
\aliascntresetthe{Cor}

\newaliascnt{Prop}{Thm}
\newtheorem{Prop}[Prop]{Proposition}
\aliascntresetthe{Prop}

\newaliascnt{Ex1}{Thm}
\newtheorem{Ex1}[Ex1]{Example}
\aliascntresetthe{Ex1}

\newaliascnt{Rem1}{Thm}
\newtheorem{Rem1}[Rem1]{Remark}
\aliascntresetthe{Rem1}

\newaliascnt{Conj}{Thm}

\aliascntresetthe{Conj}

\newaliascnt{assumption}{Thm}
\newtheorem{assumption}[assumption]{Assumption}
\aliascntresetthe{assumption}

\crefname{Thm}{Theorem}{Theorems}
\Crefname{Thm}{Theorem}{Theorems}

\crefname{Lem}{Lemma}{Lemmas}
\Crefname{Lem}{Lemma}{Lemmas}

\crefname{Def}{Definition}{Definitions}
\Crefname{Def}{Definition}{Definitions}

\crefname{Cor}{Corollary}{Corollaries}
\Crefname{Cor}{Corollary}{Corollaries}

\crefname{Prop}{Proposition}{Propositions}
\Crefname{Prop}{Proposition}{Propositions}

\crefname{Ex1}{Example}{Examples}
\Crefname{Ex1}{Example}{Examples}

\crefname{Rem1}{Remark}{Remarks}
\Crefname{Rem1}{Remark}{Remarks}

\crefname{Conj}{Conjecture}{Conjectures}
\Crefname{Conj}{Conjecture}{Conjectures}

\crefname{assumption}{Assumption}{Assumptions}
\Crefname{assumption}{Assumption}{Assumptions}

\Crefname{section}{Section}{Sections}
\Crefname{subsection}{section}{sections}
\Crefname{appendix}{Appendix}{Appendices}
\newcommand{\CrefLemmaSpace}{\nobreakspace}
\Crefformat{Lem}{Lemma\CrefLemmaSpace#2#1#3}
\makeatletter
\Crefformat{equation}{\textup{\tagform@{#2#1#3}}}
\makeatother
\Crefformat{enumi}{#2#1#3}
\Crefmultiformat{section}{Sections~#2#1#3}{--#2#1#3}{, #2#1#3}{ and~#2#1#3}

\newcommand{\cone}{\mathrm{Cone}}
\newcommand{\Hom}{\mathrm{Hom}}
\newcommand{\thick}{\mathrm{thick}\,}
\newcommand{\add}{\mathrm{add}}
\newcommand{\per}{\mathrm{per}}
\newcommand{\Ext}{\mathrm{Ext}}
\newcommand{\End}{\mathrm{End}}
\newcommand{\Mod}{\mathrm{Mod}}

\newcommand{\bijar}[1][]{%
	\ar[#1]
	\ar@<0.7ex>@{}[#1]|-*[@]{\sim}}

\newcommand{\pvd}{\mathrm{pvd}}

\newenvironment{Rem}{\begin{Rem1}\rm}{\end{Rem1}}
\newenvironment{Ex}{\begin{Ex1}\rm}{\end{Ex1}}

\usepackage{geometry}
\newcommand{\ten}{\otimes}
\newcommand{\lten}{\overset{\boldmath{L}}{\ten}}

\newcommand{\RHom}{\mathrm{\mathbf{R}Hom}}
\newcommand{\REnd}{\mathrm{\mathbf{R}End}}
\newcommand{\ra}{\rightarrow}

\newcommand{\iso}{\xrightarrow{_\sim}}

\newcommand{\id}{\mathbf{1}}

\newcommand{\ca}{{\mathcal A}}

\newcommand{\cc}{{\mathcal C}}
\newcommand{\cd}{{\mathcal D}}

\newcommand{\cf}{{\mathcal F}}

\newcommand{\ch}{{\mathcal H}}

\newcommand{\cm}{{\mathcal M}}

\newcommand{\cp}{{\mathcal P}}
\newcommand{\cR}{{\mathcal R}}
\newcommand{\cq}{{\mathcal Q}}
\newcommand{\cs}{{\mathcal S}}
\newcommand{\ct}{{\mathcal T}}
\newcommand{\cu}{{\mathcal U}}

\newcommand{\cx}{{\mathcal X}}
\newcommand{\cy}{{\mathcal Y}}
\newcommand{\cz}{{\mathcal Z}}

\renewcommand{\k}{\mathbf{k}}
\newcommand{\uA}{\overline{A}}

\newcommand{\Kar}{\operatorname{Kar}}



\begin{document}
	\title{Graded Higgs categories}
	
	\author{Li Fan, Yilin Wu}
	\address{Morningside Center of Mathematics\\ Chinese Academy of Sciences\\ 100190\\ Beijing\\ China
	}
	\email{fan-l17@tsinghua.org.cn}
	\address{
	University of Luxembourg\\                           
	6, avenue de la Fonte, L-4364\\
	Esch-sur-Alzette\\            
	Luxembourg
	}
	\email{yilin.wu@uni.lu}
	
	
	\dedicatory{}
	
	
\begin{abstract}
	Let
$	A=\bigoplus_{p\leqslant0,q\leqslant0}A_p^q$
	be a differential bigraded algebra and let $e\in A_0^0$ be a closed
	idempotent.  We develop a graded version of relative cluster categories
	and Higgs categories in which the cohomological and Adams gradings are
	coupled by a non-decreasing function
$	f:\mathbb Z\longrightarrow\mathbb Z.$ The associated family
	$\mathcal M_f = \operatorname{add}\{A(-p)[f(p)]\mid p\in\mathbb Z\}$ defines a silting subcategory and an adjacent $t$-structure. 
	
	Under suitable relative Serre-duality and approximation hypotheses, we
	construct a graded Higgs category
	$\mathcal H_{f}^{\mathbb Z}(A,e)$
	and prove that it is Frobenius extriangulated, with
	projective-injective subcategory
	$\mathcal P_f = \operatorname{add}\{eA(-p)[f(p)]\mid p\in\mathbb Z\}.$ Its stable category is triangle equivalent to the graded cluster category
	of the graded derived quotient
	$B=A/^{\mathbf L}AeA.$ It has relative Serre functor
	$(-)_\alpha(a)[d].$ If $f(p)-f(p-a)=r$, $\forall p\in\mathbb{Z}$, then
	$\mathcal M_f$ is a $(d+r)$-cluster-tilting subcategory. Examples arise from Keller--Scherotzke's work on graded singular Nakajima categories and from Iyama--Takahashi's work on graded maximal Cohen-Macaulay modules over isolated singularities.
	
	Finally, when $B$ has positive Gorenstein
	parameter and $B_{0}$ is proper, we construct a finite Beilinson-type silting generator of its
	graded cluster category and obtain a derived Morita description.
\end{abstract}
	
	\maketitle
	\tableofcontents

	\section{Introduction}
\subsection{Motivations}
Cluster categories were introduced by Buan--Marsh--Reineke--Reiten--Todorov
\cite{BMRRT} to categorify the combinatorics of acyclic cluster algebras.
Their orbit category construction was later generalized by Amiot
\cite{AmiotCluster} to algebras of global dimension at most two and to quivers
with potential, and extended to higher dimensions by Guo
\cite{GuoHigherCluster}. In the dg setting, a generalized version of cluster categories is defined as the cosingularity category
\[\cc_{\Gamma}=\per(\Gamma)/\pvd{\Gamma},\]
where $\Gamma$ is a Calabi--Yau completion or a Ginzburg dg algebra \cite{KellerCYCompletion}. Here $\per(\Gamma)$ denotes the perfect derived category and $\pvd{\Gamma}$ the perfectly valued derived category, whose objects are the dg $\Gamma$-modules whose underlying complexes of vector spaces are perfect. 

Cluster algebras arising in geometry typically carry coefficients. In the Frobenius categorifications of Geiss--Leclerc--Schr\"oer and in the general framework of Fu--Keller \cite{FuKellerCoefficients}, frozen variables correspond to projective--injective objects. Relative cluster categories and Higgs categories extend this approach to relative Calabi--Yau structures \cite{Wu2023,KellerWang2023}, and the construction also applies to certain ice quivers with potential beyond the relative Jacobi-finite case \cite{KellerWuHiggs}.

In this paper, we study Adams-graded analogues of these relative constructions. Additional gradings arise naturally in cluster theory, for instance in Amiot--Oppermann's work on cluster equivalence \cite{AmiotOppermann}, Keller--Scherotzke's graded Nakajima categories \cite{KellerScherotzke2016}, and Hanihara's study of cluster and singularity categories \cite{Hanihara2022}. In the relative setting, however, there is no single canonical choice of frozen subcategory, since the Adams shift and the cohomological shift can be combined according to an arbitrary non-decreasing function. We show that this choice determines both the family of frozen projective--injectives and the cluster-tilting dimension, whereas all the resulting Frobenius categories share the same stable category.

\subsection{Main results}
Let $k$ be a field and $A=\bigoplus_{p,q\leq0}A_p^q$ a differential bigraded $k$-algebra, with Adams degree $p$, cohomological degree $q$ and differential of bidegree $(0,1)$. Let $e\in A_0^0$ be a sum of pairwise orthogonal closed primitive idempotents and put $B=\overline A=A/^{\mathbf L}AeA$. For a non-decreasing function $f\colon\mathbb Z\to\mathbb Z$, put
\[\cm_f=\add\{A(-p)[f(p)]\mid p\in\mathbb Z\},\qquad \cp_f=\add\{eA(-p)[f(p)]\mid p\in\mathbb Z\}.\]
The subcategory $\cm_f$ is silting in $\per^{\mathbb Z}(A)$. The choice of $f$ may change the subcategory $\cp_f$ of frozen objects, but not the thick subcategory it generates.

The $t$-structure on $\cd^{\mathbb{Z}}(A)$ is cut out by the shifted degree $q+f(p)$, with
aisle and coaisle
\[ \begin{aligned} \cd_{\cm_f}^{\leq 0} &=\{X\mid H^q(X)_p=0 \text{ whenever } q+f(p)>0\},\\
\cd_{\cm_f}^{\geq 0} &=\{X\mid H^q(X)_p=0 \text{ whenever } q+f(p)<0\}. \end{aligned} \]
In particular, $f=0$ recovers the usual cohomological $t$-structure, while
$f(p)=p$ gives the diagonal one.

Unlike in the ungraded case, this $t$-structure need not restrict to $\per^{\mathbb Z}(A)$. For example, take $A=k[x]$ with $x$ in bidegree $(-1,0)$ and zero differential.  For $f=0$, each $H^q(A)_p$ is finite dimensional, but the truncation $\tau^{\geq 0}_0 A=A$ has infinite-dimensional total cohomology. For $f(p)=p$, on the other hand, the truncation $\tau^{\geq 0}_f A$ leaves only a single copy of $k$ in bidegree $(0,0)$. Under the assumptions $\dim_k H^0(A)_0<\infty$ and $\pvd^{\mathbb Z}(A)\subseteq\per^{\mathbb Z}(A)$, we prove in \Cref{Thm:restriction-to-perfect} that
\[ \begin{gathered} \tau_f^{\geq n}P\in\pvd^{\mathbb Z}(A) \quad\text{for all }P\in\per^{\mathbb Z}(A) \text{ and }n\in\mathbb Z\\
\Longleftrightarrow\\[-2pt]
\lim_{p\to-\infty}f(p)=-\infty \quad\text{or}\quad \dim_k H^0(A)<\infty . \end{gathered} \]
Whenever either condition holds, the $t$-structure restricts to $\per^{\mathbb Z}(A)$.

The relative perfectly valued derived category $\pvd_e^{\mathbb Z}(A)$ is the full subcategory of $\pvd^{\mathbb Z}(A)$ consisting of objects $X$ such that $Xe\simeq0$. Suppose that $A$ is graded homologically smooth. The \emph{graded relative cluster category} is defined as the idempotent completion of the relative graded cosingularity category:
\[\cc^{\mathbb Z}(A,e)=\Kar\left(\per^{\mathbb Z}(A)/\pvd_e^{\mathbb Z}(A)\right).\]
For $e=0$, we obtain the graded cluster category $\cc^{\mathbb Z}(A):=\cc^{\mathbb Z}(A,0)$.
We still write $\cm_f$ and $\cp_f$ for their images in $\cc^{\mathbb Z}(A,e)$. The \emph{graded Higgs category} associated with $f$ is the full subcategory
\[\ch_f^{\mathbb Z}(A,e)=\{X\in\cc^{\mathbb Z}(A,e)\mid\Hom(X,\cp_f[>0])=0=\Hom(\cp_f,X[>0])\}.\]
To establish the Frobenius structure of $\ch_f^{\mathbb Z}(A,e)$, we make the following assumptions.

\begin{assumption}\label{Ass:intro}
\begin{enumerate}
\item $\dim_k H^0(A)_0<\infty$.
\item The pair $(\per^{\mathbb Z}(A),\pvd_e^{\mathbb Z}(A))$ admits a relative Serre functor $\mathbb S=(-)_\alpha(a)[d+1]$, where $\alpha$ is a bigraded dg automorphism permuting the chosen frozen idempotents. More precisely, $\mathbb S$ preserves $\pvd_e^{\mathbb Z}(A)$ and there are bifunctorial isomorphisms
\[D\Hom(X,Y)\cong\Hom(Y,\mathbb S X)\qquad(X\in\pvd_e^{\mathbb Z}(A),\ Y\in\per^{\mathbb Z}(A)).\]
\item $\cp_f$ is functorially finite in $\cm_f$.
\end{enumerate}
\end{assumption}
In fact, condition \textup{(3)} holds whenever $H^0(A)$ is graded left and right Noetherian or $f$ is strictly increasing by \Cref{lem:graded-noetherian-functorial-finiteness}. We obtain the following results.

\begin{Thm}[{\normalfont\Cref{Thm:CTC-in-gradedHiggs}}]
Under \Cref{Ass:intro}, suppose in addition that either $\dim_k H^0(B)<\infty$ or $\lim_{p\to+\infty}f(p)=+\infty$.
\begin{enumerate}
\item The graded Higgs category $\ch_f^{\mathbb Z}(A,e)$ is Frobenius extriangulated with projective--injective subcategory $\cp_f$. The derived quotient induces a triangle equivalence
\[\ch_f^{\mathbb Z}(A,e)/[\cp_f]\xrightarrow{\sim}\cc^{\mathbb Z}(B),\]
under which the Serre functor is given by $(-)_{\bar\alpha}(a)[d]$, where $\bar\alpha$ is induced by $\alpha$ on $B$.

\item If there exists $r\in\mathbb Z$ such that
\[f(p)-f(p-a)=r\quad\text{for all }p\in\mathbb Z,\qquad d+r\geq1,\]
then $\cm_f$ is a $(d+r)$-cluster-tilting subcategory of $\ch_f^{\mathbb Z}(A,e)$.
\end{enumerate}
\end{Thm}
For fixed $(A,e)$, the theorem gives Frobenius categories of the same triangulated category $\cc^{\mathbb Z}(B)$, with different cluster-tilting dimensions depending on the choice of $f$. The proof uses the theory of relative Serre five-tuples developed in \Cref{app:relative-Serre-five-tuple}, which builds on Iyama and Yang's results on silting reduction and triangulated quotients \cite{IY2018,IY2020}. 

Assume further that $a>0$ and that $B_0$ is proper. Then the category $\cc^{\mathbb Z}(B)$ admits a finite silting generator. Put
\[T=\bigoplus_{\ell=0}^{a-1}B(-\ell),\qquad\Lambda=\REnd_B(T),\qquad U=\RHom_B(T,T(-1)).\]
Here $\REnd_B$ denotes morphisms of Adams degree zero. Applying Hanihara's Beilinson-type theorem \cite[Theorem~A.7]{Hanihara24}, we obtain a triangle equivalence
\[\per(\Lambda)\xrightarrow{\sim}\cc^{\mathbb Z}(B)\]
induced by $-\lten_\Lambda T$, where the Adams shift $(-1)$ corresponds to $-\lten_\Lambda U$ (see \Cref{Thm:finite-Adams-window}).


As an illustration, let $Q$ be a finite acyclic quiver and consider the relative $2$-Calabi--Yau completion of its framed quiver, where $d=a=1$ and $B\simeq\Pi_2(kQ)$ with $B_0=kQ$. In this case, we have
\[\cc^{\mathbb Z}(B)\simeq\cd^b(\operatorname{mod}kQ).\]
For $f(p)=p$, the graded Higgs category $\ch_f^{\mathbb Z}(A,e)$ admits a $2$-cluster-tilting subcategory $\cm_f$, where in its stable category, the Adams shift $(1)$ corresponds to the Auslander--Reiten translation (see \Cref{Ex:acyclic-diagonal-grading}). For $f=0$ and $Q$ of Dynkin type, $\cm_0$ is a $1$-cluster-tilting subcategory of $\ch_0^{\mathbb Z}(A,e)$. In this case, the graded Higgs category $\ch_0^{\mathbb Z}(A,e)$ is additively equivalent to the category of finitely presented Gorenstein-projective modules over the singular Nakajima category of Keller--Scherotzke \cite{KellerScherotzke2016} (see \Cref{Thm:Dynkin-constant-grading}). Another class of examples arises from skew group algebras, where $f=0$ recovers the categories of graded maximal Cohen--Macaulay modules studied by Iyama--Takahashi \cite{IyamaTakahashiQuotient}; see \Cref{cor:higgs-mcm-mckay}.

Relative Koszul duality gives a singularity-category description of these categories. Set
\[ S=H^0(A)_0/\operatorname{rad} H^0(A)_0,\qquad E\simeq\REnd_A^{\mathbb Z}(S), \]
and let $u$ be the idempotent in $E$ corresponding to the simple summands of $S$ annihilated by $e$. Put $E_{\mathrm u}=uEu$ and define the \emph{relative graded singularity category} by
\[\operatorname{sg}^{\mathbb Z}(E;u)=\Kar\left(\pvd^{\mathbb Z}(E)/\thick^{\mathbb Z}(uE)\right).\]
For $u=1$, it is the graded singularity category $\operatorname{sg}^{\mathbb Z}(E)$.  By \Cref{thm:relative-koszul,prop:frozen-quotient}, we obtain a commutative diagram of triangulated categories via Koszul duality
\[\begin{tikzcd}[column sep=large,row sep=large]
\cc^{\mathbb Z}(A,e)\arrow[r,"\sim"]\arrow[d,"\pi"']&\operatorname{sg}^{\mathbb Z}(E;u)\arrow[d,"(-)u"]\\
\cc^{\mathbb Z}(B)\arrow[r,"\sim"]&\operatorname{sg}^{\mathbb Z}(E_{\mathrm u}).
\end{tikzcd}\]
In particular, under the hypotheses of \Cref{Thm:CTC-in-gradedHiggs}, the lower equivalence realizes the stable Higgs category as a graded singularity category.

Next, we relate the graded version to the ungraded one via \emph{weighted totalizations} and orbit categories. For integers $N\geq0$ and $c$, put \[f_N(p)=Np+c,\qquad F_N=(-1)[N].\]
The functor $F_N$ permutes the generators of $\cp_{f_N}$. Totalization in degree $q+Np$ gives dg algebras $A_N^{\mathrm{tot}}=\operatorname{Tot}_N(A)$ and $B_N^{\mathrm{tot}}=\operatorname{Tot}_N(B)$, with a natural isomorphism $\operatorname{Tot}_N F_N\simeq\operatorname{Tot}_N$. Under \Cref{Ass} for $f=f_N$, we obtain the following commutative diagram of triangulated categories (see \Cref{thm:weighted-orbit}):
\[\begin{tikzcd}[column sep=large,row sep=large]
(\cc^{\mathbb Z}(A,e)/F_N)_\triangle^\sharp\arrow[r,"\operatorname{Tot}_N","\sim"']\arrow[d,"\pi"']&\cc(A_N^{\mathrm{tot}},e)\arrow[d,"\pi_N"]\\
(\cc^{\mathbb Z}(B)/F_N)_\triangle^\sharp\arrow[r,"\operatorname{Tot}_N","\sim"']&\cc(B_N^{\mathrm{tot}}).
\end{tikzcd}\]
To get an analogous result for Higgs categories, one needs an orbit construction that retains the extriangulated structure. For this purpose, using Chen's theory of exact dg categories \cite{Chen20261,Chen20262}, we construct the \emph{exact dg orbit categories} in \Cref{appB:exact-dg-orbit}. We prove that the construction preserves the Frobenius property and that taking stable categories commutes with taking orbit up to idempotent completion. Applying this to graded Higgs categories, we have the following result.

\begin{Thm}[{\normalfont\Cref{cor:weighted-totalization-higgs-orbit}}]
Under \Cref{Ass} for $f=f_N$ and conditions \eqref{eq:weighted-bimodule-hypothesis} and \eqref{eq:totalized-Krull-Schmidt}, the functor $J_N=[-c]\circ\operatorname{Tot}_N$ is fully faithful on the ordinary orbit and induces an exact equivalence
\[(\ch_{f_N}^{\mathbb Z}(A,e)/F_N)^{\mathrm{ex},\sharp}\xrightarrow{\sim}\ch(A_N^{\mathrm{tot}},e).\]
The stable category is triangle equivalent to $\cc(B_N^{\mathrm{tot}})$. Here $(-)^{\mathrm{ex},\sharp}$ denotes the homotopy category of the exact dg orbit hull constructed in \Cref{appB:exact-dg-orbit}.
\end{Thm}

Finally, the construction of graded Higgs categories extends to algebras graded by an abelian group $L$ with support in a submonoid $\Sigma$. We refer to \Cref{sec:other-gradings} for details.

\subsection{Contents}
The paper is organized as follows. In \Cref{sec:preliminaries}, we recall some basic concepts about bigraded dg algebras, dg quotients, and graded relative Calabi--Yau completions. In \Cref{sec:graded-t-structure}, we give the $t$-structures on perfect derived categories. In Section~\ref{sec:graded-higgs}, we study the graded Higgs categories and show the Frobenius and cluster-tilting properties, together with the construction of a silting generator for the corresponding graded cluster category. Moreover, we work out examples from framed quivers and graded maximal Cohen--Macaulay modules over isolated quotient singularities, and describe graded relative cluster categories via Koszul duality. In Section~\ref{sec:weighted-totalization-singularity}, we relate graded and ungraded Higgs categories via weighted totalizations and orbit constructions. In Section~\ref{sec:other-gradings}, we extend our constructions to abelian grading groups with a given support monoid. Finally, in Appendix~\ref{app:relative-Serre-five-tuple}, we treat the relative Serre five-tuples and in Appendix~\ref{appB:exact-dg-orbit}, we construct exact dg orbit hulls and study their relation to triangulated orbit categories via stable categories.
	
\subsection*{Acknowledgments}
Part of this work was completed during the second author's visit to Paris.
He thanks Bernhard Keller for the invitation and for helpful discussions.
	
	\subsection*{Use of generative AI}
	The authors used GPT for English-language and \LaTeX{} editing,
	manuscript organization, discussion of the presentation of proofs,
	and preliminary checks of proofs and calculations. The authors take
	responsibility for the content of the manuscript.
	
	\subsection*{Notation and conventions}
	For a homogeneous element $v\in V_p^q$, write
	$\operatorname{bideg}(v)=(p,q)$, $\deg(v)=p$, and $|v|=q$.
	Our shift conventions are
	\[ V(m)_p^q=V_{p+m}^q,\qquad V[m]_p^q=V_p^{q+m}. \]
	We write $D$ for the graded $k$-dual, with $(DV)^q_p=\Hom_k(V^{-q}_{-p},k)$ and the usual differential.
	We write $\Kar(\mathcal V)$ for the idempotent completion of an additive category $\mathcal V$.

	For a triangulated category $\mathcal V=H^0(\mathcal E)$ with a
	pretriangulated dg enhancement and an autoequivalence $F$ with dg
	lift $\widetilde F$, we denote the triangulated hull of the orbit category $\mathcal V/F$ in the sense of \cite[Definition~3.19]{FKQ2024} by
	\[ (\mathcal V/F)_\triangle =H^0\!\left(\operatorname{pretr} (\mathcal E/\widetilde F^{\mathbb Z})\right). \]
	A superscript $\sharp$ on the triangulated hull denotes its
	idempotent completion.

\section{Preliminaries}\label{sec:preliminaries}
	In this section, we recall some background on differential bigraded algebras and their derived
	categories.
	
	\begin{Def}\rm
		A \emph{differential bigraded $k$-algebra} is a bigraded $k$-algebra $A=\bigoplus_{(i,j)\in\mathbb{Z}^2}A^j_i$ equipped with
		a degree $(0,1)$ $k$-linear map $d_A$ that squares to $0$ and satisfies the Leibniz rule:
		\[ d_A(xy) = d_A(x)y+(-1)^{|x|}xd_A(y). \]
		
		A \emph{morphism} of differential bigraded
		$k$-algebras is a degree $(0,0)$ $k$-algebra homomorphism that commutes with differentials. Such a
		morphism is a \emph{quasi-isomorphism} if it induces an isomorphism on cohomology.
	\end{Def}
	
	We will refer to differential bigraded $k$-algebras as simply “dbg-algebras”. We denote
	by $A^{\natural}$ the underlying bigraded $k$-algebra of a dg-algebra $A$.

	\begin{Def}\rm
		A dbg-algebra $A$ is \emph{connected} if $A_{i}^{j}=0$ when $i>0$ or $j>0$.
	\end{Def}	

\subsection{Differential bigraded modules and derived categories}

Let $A$ be a dbg algebra. A \emph{differential bigraded right $A$-module} is a bigraded right $A$-module $M=\bigoplus_{(i,j)\in\mathbb{Z}^2}M_i^j$ endowed with a differential $d_M$ of bidegree $(0,1)$ such that $d_M^2=0$ and
\[ d_M(ma)=d_M(m)a+(-1)^{|m|}m d_A(a) \]
for all homogeneous elements $m\in M$ and $a\in A$. We use the notation
\[ M_i=\bigoplus_{j\in\mathbb{Z}}M_i^j\qquad\text{and}\qquad M^j=\bigoplus_{i\in\mathbb{Z}}M_i^j. \]
We denote by $\mathcal{C}^{\mathbb{Z}}(A)$ the category of differential bigraded right $A$-modules. Its morphisms are right $A$-linear maps of bidegree $(0,0)$ commuting with the differentials. A morphism $s:M\to L$ in $\mathcal{C}^{\mathbb{Z}}(A)$ is called a \emph{quasi-isomorphism} if $s_p:M_p\to L_p$ is a quasi-isomorphism of complexes of $k$-vector spaces for every $p\in\mathbb{Z}$. The \emph{bigraded derived category} of $A$ is the localization
\[ \mathcal{D}^{\mathbb{Z}}(A)=\mathcal{C}^{\mathbb{Z}}(A)[\mathrm{Qis}^{-1}]. \]
The category $\mathcal{D}^{\mathbb{Z}}(A)$ is compactly generated by $\{A(i)\mid i\in\mathbb{Z}\}$. We define the \emph{perfect derived category} of $A$ as
\[ \operatorname{per}^{\mathbb{Z}}(A)\coloneqq \operatorname{thick}_{\mathcal{D}^{\mathbb{Z}}(A)}\{A(i)\mid i\in\mathbb{Z}\}. \]

The \emph{graded enveloping algebra} of $A$ is defined as
$A^e:=A\otimes_k A^{op},$
equipped with the graded tensor product multiplication. We regard $A$ as
a graded $A^e$-module via its natural graded $A$-$A$-bimodule structure.

\begin{Def}\rm
	We say that $A$ is \emph{graded homologically smooth} if $A$, regarded as
	a graded $A^e$-module, is perfect, that is,
	$A\in \per^{\mathbb Z}(A^e).$
\end{Def}

If $k$ is regarded as a dbg algebra with trivial grading shift and we denote it $\k$, then we define the \emph{perfectly valued derived category} as $\per^{\mathbb{Z}}(\k)$, which is the thick subcategory generated by $\k(p)$ for $p\in\mathbb{Z}$, i.e.
\[ \pvd^{\mathbb{Z}}(\k)=\per^{\mathbb{Z}}(\k)=\{M\in\cd^{\mathbb{Z}}(\k)\mid\sum_{i,j\in\mathbb{Z}}\dim H^j(M)_i<\infty\}. \]
Then we define the \emph{perfectly valued derived category $\pvd^{\mathbb{Z}}$(A)} as
\[ \operatorname{pvd}^{\mathbb Z}(A)= \{M\in D^{\mathbb Z}(A)\mid M|_{\k}\in\pvd^{\mathbb{Z}}(\k)\}. \]

\bigskip

For $L,M\in\mathcal{D}^{\mathbb{Z}}(A)$ and $p\in\mathbb{Z}$, we have
\[ \operatorname{Hom}_{\mathcal{D}^{\mathbb{Z}}(A)}(L,M[p])=H^p\operatorname{RHom}_A(L,M). \]
We also use the graded Hom space
\[ \operatorname{Hom}^{\mathbb{Z}}_{\mathcal{D}^{\mathbb{Z}}(A)}(L,M[p])=\bigoplus_{q\in\mathbb{Z}}\operatorname{Hom}_{\mathcal{D}^{\mathbb{Z}}(A)}(L,M[p](q)), \]
and set
\[ \operatorname{RHom}^{\mathbb{Z}}_A(L,M)=\bigoplus_{q\in\mathbb{Z}}\operatorname{RHom}_A(L,M(q)). \]
Thus
\[ \operatorname{Hom}^{\mathbb{Z}}_{\mathcal{D}^{\mathbb{Z}}(A)}(L,M[p])=H^p\operatorname{RHom}^{\mathbb{Z}}_A(L,M). \]

There are two shift functors in $\cd^{\mathbb{Z}}(A)$.
One is the shift of differential degree:
\[ M\mapsto M[1]~\mbox{with}~(M[1])^p=M^{p+1}~\mbox{and}~d_{M[1]}=-d_M, \]
and the other is the shift of Adams degree:
\[ M\mapsto M(1)~\mbox{with}~(M(1))_p=M_{1+p}~\mbox{and}~d_{M(1)}=d_M. \]
\begin{Rem}
If $A$ is a dg algebra considered as a dbg algebra concentrated in Adams degree 0, we have
\begin{equation} \per^{\mathbb{Z}}(A)=\thick(A(p), p\in\mathbb{Z})\subseteq\cd^{\mathbb{Z}}(A). \end{equation}
\end{Rem}	
Let $e\in A_0^0$ be an idempotent of $A$. We define
\[ \pvd_e^{\mathbb{Z}}(A)=\left\{M\in\pvd^{\mathbb{Z}}(A)\middle|Me=0\text{ in }\mathcal{D}^{\mathbb{Z}}(eAe)\right\}, \]
which is a full triangulated subcategory of $\pvd^{\mathbb{Z}}(A)$. Equivalently, $M\in\pvd_e^{\mathbb{Z}}(A)$ if and only if $M\overset{\mathbf{L}}{\otimes}_A Ae=0$. We denote by
\[ \thick^{\mathbb{Z}}(eA)=\operatorname{thick}_{\mathcal{D}^{\mathbb{Z}}(A)}\left\{eA(i)\mid i\in\mathbb{Z}\right\}\subseteq\per^{\mathbb{Z}}(A) \]
the thick subcategory generated by all internal shifts of $eA$.

	\subsection{Dg quotient in bigraded setting}\label{Subsection:dg quotient}
	Let $A$ be a bigraded dg algebra and $e$ a closed idempotent of $A$, i.e. $e\in A_0^0$ and $d_A(e)=0$. We use the following Adams-graded version of Drinfeld's construction \cite[Section~7.2]{Kalck-Yang2018}.
	
	Define a dg $ k $-category $ \ca $ with two objects $ \epsilon $ and $ \gamma $ such that the morphism spaces are given by:
	\begin{align*} \text{Hom}_\ca(\epsilon, \epsilon) &= eAe, \\
\text{Hom}_\ca(\epsilon, \gamma) &= (1 - e)Ae, \\
\text{Hom}_\ca(\gamma, \epsilon) &= eA(1 - e), \\
\text{Hom}_\ca(\gamma, \gamma) &= (1 - e)A(1 - e), \end{align*}
	
	and the composition of morphisms is induced from the multiplication of $A$. Let $ \ca' $ be the dg subcategory of $ \ca $ consisting of the object $ \epsilon $.  We form the dg quotient $ \ca / \ca' $ by formally adjoining a morphism $ x: \epsilon \to \epsilon $ of degree $(0,-1)$ such that $ d(x) = \text{id}_{\epsilon} = e $.
In the finite additive hull of $\ca/\ca'$, we take the endomorphism dg
algebra of $\epsilon\oplus\gamma$ and denote it by $\overline A$.
The quotient functor induces a canonical unital morphism
$p:A\to\overline A$. Explicitly,
\[ \overline A\cong A\langle x\rangle/(ex=xe=x),\qquad \operatorname{bideg}(x)=(0,-1),\qquad d(x)=e. \]
Since $\epsilon$ is contractible in the quotient, this choice gives the same
derived category as $\operatorname{End}_{\ca/\ca'}(\gamma)$.

	As a differential bigraded module, the underlying graded module of $\overline{A}$ has the form
	\begin{align*} \cdots & \longrightarrow Ae \otimes (eAe)^{\otimes p} \otimes eA \overset{d^{-p-1}}{\longrightarrow} \cdots \longrightarrow Ae \otimes eAe \otimes eA \overset{d^{-2}}{\longrightarrow} Ae \otimes eA \overset{d^{-1}}{\longrightarrow} A, \end{align*}
	The term $Ae\otimes(eAe)^{\otimes p}\otimes eA$ corresponds to words
	$a_0xa_1x\cdots xa_{p+1}$ containing $p+1$ copies of $x$. For homogeneous
	$a_0,\ldots,a_{p+1}$, such a word has bidegree
	\[ \left(\sum_{j=0}^{p+1}\deg(a_j),\, \sum_{j=0}^{p+1}|a_j|-(p+1)\right), \]
	since $\operatorname{bideg}(x)=(0,-1)$. Thus the displayed position $-p-1$
	records only the contribution of the copies of $x$; the factors $a_j$ retain
	their own cohomological degrees and differentials.
	
	For homogeneous $a_0,\ldots,a_{p+1}$, the differential of $\overline{A}$ is
	\begin{align*} d_{\overline{A}}(a_0\otimes\cdots\otimes a_{p+1}) ={}&\sum_{j=0}^{p+1}(-1)^{\sigma_j} a_0\otimes\cdots\otimes d_A(a_j)\otimes\cdots\otimes a_{p+1}\\
&+\sum_{i=0}^{p}(-1)^{\tau_i} a_0\otimes\cdots\otimes a_i a_{i+1}\otimes\cdots\otimes a_{p+1}, \end{align*}
	where
	$\sigma_j=\sum_{t=0}^{j-1}|a_t|-j, \tau_i=\sum_{t=0}^{i}|a_t|-i.$
	The two sums are induced by $d_A$ and $d(x)=e$, respectively.
	The multiplication of $ \overline{A} $ is induced from that of $ A $
	\[ (a_0 \otimes \cdots \otimes a_{p+1})(b_0 \otimes \cdots \otimes b_{q+1}) = a_0 \otimes \cdots \otimes a_p \otimes a_{p+1} b_0 \otimes b_1 \otimes \cdots \otimes b_{q+1}. \]

\begin{Rem}
 Suppose that $A=(kQ,d)$ is a differential bigraded quiver algebra, where $Q$ has finitely many vertices, and that $e=\sum_{i\in I}e_i$ is a sum of vertex idempotents. Since $d(e)=0$, the ideal $AeA$ is a differential bigraded ideal, and hence $A/AeA$ inherits a differential bigraded algebra structure. The proof of \cite[Theorem 7.1]{Kalck-Yang2018} is compatible with the Adams grading. Hence the dg quotient $\overline{A}$ is quasi-equivalent, as a differential bigraded algebra, to
	$A/AeA$.

Equivalently, $A/AeA$ is obtained from $A$ by deleting the vertices belonging to $I$, together with all paths passing through these vertices, and equipping the resulting quotient with the induced differential and bigrading.
\end{Rem}

	\begin{Prop}\label{prop:rec}
		We have the following recollement of bigraded derived categories	
\[\begin{tikzpicture}[>=latex] \node (L) at (0,0) {$\mathcal D^\mathbb{Z}(\uA)$}; \node (M) at (4,0) {$\mathcal D^\mathbb{Z}(A)$}; \node (R) at (8,0) {$\mathcal D^\mathbb{Z}(eAe)$}; \draw[->] (L) -- node[above] {$p_* = p_!$} (M); \draw[->] (M) -- node[above] {$i_* = i_!$} (R); \draw[<-] (L) to[out=30, in=150] node[above] {$p^*$} (M); \draw[->] (M) to[out=-150, in=-30] node[below] {$p^!$} (L); \draw[<-] (M) to[out=30, in=150] node[above] {$i^*$} (R); \draw[->] (R) to[out=-150, in=-30] node[below] {$i^!$} (M); \end{tikzpicture}\]
		where
		\[ \begin{aligned} p^*&=-\overset{\mathbf L}{\otimes}_A {}_A\uA_{\uA},\\
p_*=p_!&=\operatorname{Res}_p \simeq -\overset{\mathbf L}{\otimes}_{\uA}{}_{\uA}\uA_A,\\
p^!&=\mathbf R\!\operatorname{Hom}^{\mathbb{Z}}_A({}_{\uA}\uA_A,-), \end{aligned} \qquad \begin{aligned} i^*&=-\overset{\mathbf L}{\otimes}_{eAe}{}_{eAe}eA_A,\\
i_*=i_!&=-\overset{\mathbf L}{\otimes}_A{}_A Ae_{eAe} \simeq \mathbf R\!\operatorname{Hom}^{\mathbb{Z}}_A({}_{eAe}eA_A,-),\\
i^!&=\mathbf R\!\operatorname{Hom}^{\mathbb{Z}}_{eAe}({}_A Ae_{eAe},-). \end{aligned} \]
	\end{Prop}

\begin{proof}
The proof of \cite[Theorem~5.1.3]{ChenChen2021} also applies
in the Adams-graded setting over the field $k$, since the
resolutions and homotopies used there can be chosen to preserve
Adams degree.
Applying it to
$\mathcal A'^{\mathrm{op}}\subseteq\mathcal A^{\mathrm{op}}$
from the preceding quotient construction, and identifying
the adjoints by graded tensor--Hom adjunction, gives the
stated recollement, with
\[\operatorname{Im}(p_*)=\ker(i_*)=\{M\in\mathcal D^{\mathbb Z}(A)\mid Me\simeq0\}.\]

Moreover, full faithfulness of $p_*$ makes the counit at
$\overline A$ an isomorphism; this is the multiplication map
$\overline A\otimes_A^{\mathbf L}\overline A \xrightarrow{\sim}\overline A$.
Since $p^*(A(r))\cong\overline A(r)$, the functor $p^*$ preserves
perfect objects. The compact-localization assertion of the same
theorem, followed by idempotent completion, gives
\[\operatorname{Kar}\!\left(\operatorname{per}^{\mathbb Z}(A)/\operatorname{thick}^{\mathbb Z}(eA)\right)\xrightarrow{\sim}\operatorname{per}^{\mathbb Z}(\overline A).\]
All these constructions commute with Adams shifts.
\end{proof}
\par

	\subsection{Relative Calabi--Yau completions in the graded setting}
\label{subsec:graded-relative-Ginzburg}

Let $n$ be a positive integer. Let $(Q,F)$ be a finite ice quiver,
where $Q$ is a finite quiver and $F$ is a subquiver of $Q$.
Fix an Adams grading and a cohomological grading on the arrows of $Q$,
\[ \deg_{\mathrm A}\colon Q_1\longrightarrow\mathbb Z, \qquad |-|\colon Q_1\longrightarrow\mathbb Z, \]
and extend them additively to paths and to $kQ$.
For homogeneous paths $u,v$, write
\[ [u,v]_{\mathrm{gr}} := uv-(-1)^{|u||v|}vu, \]
and set
$(kQ)^{\mathrm{gr}}_{\mathrm{cyc}}:=kQ/[kQ,kQ]_{\mathrm{gr}}$.
Let $W\in (kQ)^{\mathrm{gr}}_{\mathrm{cyc}}$ be a potential. We assume only
that $W$ is bihomogeneous of some fixed
bidegree
\[ \operatorname{bideg}(W)=(\ell,3-n), \qquad \ell\in\mathbb Z. \]
Only the cohomological grading enters the Koszul signs; the Adams grading
is an auxiliary grading preserved by the differential.

Let $\widetilde Q$ be the bigraded quiver obtained from $Q$ by adjoining
\begin{itemize}
\item an opposite arrow
$a^\vee\colon t(a)\longrightarrow s(a)$
for every unfrozen arrow $a\in Q_1\setminus F_1$;
\item a loop
$t_i\colon i\longrightarrow i$
for every mutable vertex $i\in Q_0\setminus F_0$.
\end{itemize}
The first component below is the Adams degree and the second is the
cohomological degree. We set

\begin{equation} \label{eq:relative-Ginzburg-arbitrary-weight-bidegrees} \begin{aligned} \operatorname{bideg}(a) &= \bigl(\deg_{\mathrm A}(a),|a|\bigr), &&a\in Q_1,\\
\operatorname{bideg}(a^\vee) &= \bigl(\ell-\deg_{\mathrm A}(a),2-n-|a|\bigr), &&a\in Q_1\setminus F_1,\\
\operatorname{bideg}(t_i) &= (\ell,1-n), &&i\in Q_0\setminus F_0. \end{aligned} \end{equation}

For an arrow $a\in Q_1$, define the graded cyclic derivative by
\[ \partial_a[w]=\sum_{w=uav}(-1)^{|u|(|a|+|v|)}vu \]
for a homogeneous cyclic path $w$, and extend it linearly to potentials.
\begin{equation} \label{eq:degree-of-cyclic-derivative} \operatorname{bideg}(\partial_aW) = \bigl(\ell-\deg_{\mathrm A}(a),3-n-|a|\bigr). \end{equation}

\begin{Def}\rm
\label{def:graded-complete-relative-Ginzburg}
The \emph{graded relative $n$-dimensional Ginzburg dg algebra} associated
with
	the bigraded ice quiver with potential $(Q,F,W)$
is
\[ \Gamma^{\mathbb Z}_{\mathrm{rel},\ell,n}(Q,F,W) :=\bigl(k\widetilde Q,d\bigr), \]
where $d$ is the unique derivation of bidegree $(0,1)$ such that
\begin{align} d(a) &= 0, &&a\in Q_1, \label{eq:graded-relative-Ginzburg-d-a}\\
d(a^\vee) &=(-1)^{|a|} \partial_aW, &&a\in Q_1\setminus F_1, \label{eq:graded-relative-Ginzburg-d-dual}\\
d(t_i) &= e_i \left( \sum_{a\in Q_1\setminus F_1} [a,a^\vee]_{\mathrm{gr}} \right) e_i, &&i\in Q_0\setminus F_0. \label{eq:graded-relative-Ginzburg-d-loop} \end{align}
Here $[-,-]_{\mathrm{gr}}$ is the supercommutator for the cohomological
grading.
\end{Def}

The degree choices are forced by the requirement that $d$ have bidegree
$(0,1)$. Indeed,
\[ \operatorname{bideg}(\partial_aW) = \operatorname{bideg}(a^\vee)+(0,1), \qquad \operatorname{bideg}\bigl([a,a^\vee]_{\mathrm{gr}}\bigr) = \operatorname{bideg}(t_i)+(0,1). \]
If all original arrows have cohomological degree zero, a nonzero
homogeneous potential forces $n=3$. For $n\neq3$, one may instead take
$W=0$, which is homogeneous of any degree, or equip $Q$ with a nontrivial
cohomological grading.

Let
$\Pi^{\mathbb Z}_{n-1,\ell}(F)$
be the bigraded derived $(n-1)$-dimensional preprojective dg
algebra of $F$. 
Its underlying algebra is the ordinary path algebra of the
bigraded quiver described below. The quiver has,
in addition to every $b\in F_1$, an opposite arrow $\widetilde b$ and a loop $r_i$ at
every $i\in F_0$, with

\[ \begin{aligned} \operatorname{bideg}(b) &= \bigl(\deg_{\mathrm A}(b),|b|\bigr),\\
\operatorname{bideg}(\widetilde b) &= \bigl(\ell-\deg_{\mathrm A}(b),3-n-|b|\bigr),\\
\operatorname{bideg}(r_i) &= (\ell,2-n), \end{aligned} \]
and differential
\[ d(b)=d(\widetilde b)=0, \qquad d(r_i) = e_i \left( \sum_{b\in F_1} [b,\widetilde b]_{\mathrm{gr}} \right) e_i. \]

\begin{Prop}
\label{prop:graded-Ginzburg-morphism-arbitrary-weight}
There is a morphism of bigraded dg algebras
\begin{equation} \label{eq:graded-Ginzburg-morphism-arbitrary-weight} G_{\mathrm{rel},\ell}\colon \Pi^{\mathbb Z}_{n-1,\ell}(F) \longrightarrow \Gamma^{\mathbb Z}_{\mathrm{rel},\ell,n}(Q,F,W) . \end{equation}
It is given on generators by
\begin{align*} G_{\mathrm{rel},\ell}(e_i) &= e_i, & G_{\mathrm{rel},\ell}(b) &= b,\\
G_{\mathrm{rel},\ell}(\widetilde b) &= -\partial_bW, & G_{\mathrm{rel},\ell}(r_i) &= e_i \left( \sum_{a\in Q_1\setminus F_1} [a,a^\vee]_{\mathrm{gr}} \right) e_i. \end{align*}
The morphism
\Cref{eq:graded-Ginzburg-morphism-arbitrary-weight}
is the $W$-deformed reduced relative $n$-Calabi--Yau completion
of the inclusion $kF\hookrightarrow kQ$.
It carries the canonical homogeneous left relative
$n$-Calabi--Yau structure by
\cite[Section~3.6 and Proposition~3.18]{Wu2023}.

\end{Prop}
\begin{Rem}
\label{rem:potential-weight-versus-Gorenstein-parameter}
With our shift convention
$M(r)^q_p=M^q_{p+r},$
an Adams-homogeneous Calabi--Yau structure of weight $\ell$ corresponds to
the internal shift $(-\ell)$. Thus, in the absolute case
$F=\varnothing$, one has
\[ \mathbf R\!\operatorname{Hom}_{\Gamma^e} (\Gamma,\Gamma^e)(-\ell)[n] \simeq \Gamma, \qquad \Gamma = \Gamma^{\mathbb Z}_{\mathrm{rel},\ell,n} (Q,\varnothing,W) . \]
Consequently, in the terminology of the preceding subsection, the
Gorenstein parameter is
\[ a=-\ell. \]

In particular, if $\ell=-a<0$ and the arrow grading is chosen so that
the resulting relative Ginzburg dg algebra is non-positively Adams
graded, then its Gorenstein parameter is $a>0$. For example, this holds
when all original arrows have Adams degree zero. If $\ell=0$ and all
original arrows have Adams degree zero, then the whole relative
Ginzburg dg algebra is concentrated in Adams degree zero, but not, in
general, in cohomological degree zero.

\end{Rem}

\section{$t$-structure on $\per^{\mathbb{Z}}A$}\label{sec:graded-t-structure}
Let $\ct$ be a triangulated category with arbitrary (set-indexed) coproducts. An object $X\in\ct$ is \emph{compact} if the natural map $\bigoplus_{i\in I}\Hom_{\ct}(X,X_{i})\ra\Hom_{\ct}(X,\bigoplus_{i\in I}X_{i})$ is an isomorphism for every set $\{X_{i}\,|\,i\in I\}$ of objects. A subcategory is said to \emph{generate} $\ct$ if the smallest triangulated subcategory containing it and closed under arbitrary coproducts is the whole $\ct$ . We call a subcategory of $\ct$ a \emph{compact set of generators} if it is skeletally small, consists of compact objects and generates $\ct$.

We call a subcategory $\mathcal{M}$ of $\mathcal{T}$ \emph{silting} if it is a compact set of generators and
\[ \operatorname{Hom}_{\mathcal{T}}(\mathcal{M},\mathcal{M}[k])=0 \]
for every $k>0$.

Let $A=\displaystyle\bigoplus_{i\leqslant 0,j\leqslant 0} A^{j}_{i}$ be a dbg-algebra over $k$. Let $f\colon \mathbb{R}\to\mathbb{R}$ be a non-decreasing function such that $f(\mathbb{Z})\subseteq \mathbb{Z}.$ Then its restriction $f|_{\mathbb{Z}}\colon \mathbb{Z}\to\mathbb{Z}$ is also non-decreasing, i.e. for all $i\leqslant j$, we have $f(i)\leqslant f(j)$. 
	
\begin{Prop}\label{prop:graded-silting-subcategory}
	The subcategory $\cm_{f}=\add\{A(-i)[f(i)]\,|\,i\in\mathbb{Z}\}$ is a silting subcategory of $\cd^{\mathbb{Z}}(A)$.
	\end{Prop}
	
	\begin{proof}
	It is clear that $\cm_{f}$ is a compact set of generators for $\cd^{\mathbb{Z}}(A)$ and  $\thick(\cm_{f})=\per^{\mathbb{Z}}(A)$. It remains to show that $$\Hom_{\per^{\mathbb{Z}}(A)}(A,A(i-j)[f(j)-f(i)+k])$$ vanishes for each $i,j\in\mathbb{Z}$ and $k>0$. Since $A_{i}^{j}=0$ for $i>0$ or $j>0$ and $f$ is non-decreasing, we see that $\Hom_{\per^{\mathbb{Z}}(A)}(A,A(i-j)[f(j)-f(i)+k])$ vanishes for each $i,j\in\mathbb{Z}$ and $k>0$.	
	\end{proof}

\bigskip
	
We have the standard $t$-structure $(\cd_{\cm_{f}}^{\leqslant0},\cd_{\cm_{f}}^{\geqslant0})$ in $\cd^{\mathbb{Z}}(A)$ associated to $\cm_{f}$, which is given by
$$\cd_{\cm_{f}}^{\leqslant0}=\cm_{f}[<\!0]^{\perp_{\cd^{\mathbb{Z}}({A})}}=\{X\in\cd^{\mathbb{Z}}(A)\,|\,\Hom_{\cd^{\mathbb{Z}}(A)}(M,X[k])=0,\,\text{for all $M\in\cm_{f}$ and $k>0$}\},$$	
$$\cd_{\cm_{f}}^{\geqslant0}=\cm_{f}[>\!0]^{\perp_{\cd^{\mathbb{Z}}({A})}}=\{X\in\cd^{\mathbb{Z}}(A)\,|\,\Hom_{\cd^{\mathbb{Z}}(A)}(M,X[k])=0,\,\text{for all $M\in\cm_{f}$ and $k<0$}\}.$$	

We write $\tau_f^{\leq n}$ and $\tau_f^{\geq n}$ for its truncation functors. They satisfy
\begin{equation}\label{eq:f-truncation-cohomology} H^m\bigl(\tau_f^{\geq n}X\bigr)_i \cong \begin{cases} H^m(X)_i, & m+f(i)\geq n,\\
0, & m+f(i)<n. \end{cases} \end{equation}

\begin{equation} H^m\bigl(\tau_f^{\leq n}X\bigr)_i \cong \begin{cases} H^m(X)_i, & m+f(i)\leq n,\\
0, & m+f(i)>n. \end{cases} \end{equation}

\begin{Lem}\label{lem:t-stru-onD(A)}
	We have
	$$\cd_{\cm_{f}}^{\leqslant0}=\{X\in\cd^{\mathbb{Z}}(A)\,|\,H^{k-f(i)}(X)_{i}=0,\,\text{for all $k>0$ and $i\in\mathbb{Z}$}\},$$
	$$\cd_{\cm_{f}}^{\geqslant0}=\{X\in\cd^{\mathbb{Z}}(A)\,|\,H^{k-f(i)}(X)_{i}=0,\,\text{for all $k<0$ and $i\in\mathbb{Z}$}\}.$$
Moreover, suppose that
$\lim_{i\to-\infty} f(i)=c$ and
$\lim_{i\to\infty} f(i)=c',$
where $c\leq c'$ in $\mathbb{R}\cup\{\pm\infty\}$. For each
$m\in\mathbb{Z}$ with $-c'<m\leq -c$, define
\[ b_m=\max\{i\in\mathbb{Z}\mid f(i)\leq -m\}, \]
and for each $m\in\mathbb{Z}$ with $-c'\leq m<-c$, define
\[ a_m=\min\{i\in\mathbb{Z}\mid f(i)\geq -m\}. \]
Then
\[ \mathcal{D}_{\mathcal{M}_f}^{\leqslant 0} = \left\{ X\in\mathcal{D}^{\mathbb{Z}}(A) \ \middle|\ \begin{array}{l} H^m(X)\in\Mod(H^0(A))_{\leq b_m} \text{ for all }-c'<m\leq -c;\\
H^m(X)=0\text{ for all }m>-c \end{array} \right\} \]
and
\[ \mathcal{D}_{\mathcal{M}_f}^{\geqslant 0} = \left\{ X\in\mathcal{D}^{\mathbb{Z}}(A) \ \middle|\ \begin{array}{l} H^m(X)\in\Mod(H^0(A))_{\geq a_m} \text{ for all }-c'\leq m<-c;\\
H^m(X)=0\text{ for all }m<-c' \end{array} \right\}. \]
Here $\Mod(H^0(A))_{\leq b}$ (respectively,
$\Mod(H^0(A))_{\geq a}$) denotes the full subcategory of graded
$H^0(A)$-modules supported in Adams degrees at most $b$
(respectively, at least $a$).

\end{Lem}	
\begin{proof}
We have
\begin{align*} \cd_{\cm_{f}}^{\leqslant0}&=\{X\in\cd^{\mathbb{Z}}(A)\,|\,\Hom_{\cd^{\mathbb{Z}}(A)}(A(-i)[f(i)-k],X)=0,\,\text{for all $i\in\mathbb{Z}$ and $k>0$}\}\\
&=\{X\in\cd^{\mathbb{Z}}(A)\,|\,H^{k-f(i)}(X)_{i}=0,\,\text{for all $i\in\mathbb{Z}$ and $k>0$}\}\\
&=\{X\in\cd^{\mathbb{Z}}(A)\,|\,H^{m}(X)_{i}=0,\,\text{for all $i\in\mathbb{Z}$ and $f(i)>-m$}\} \end{align*}

\begin{align*} \cd_{\cm_{f}}^{\geqslant0}&=\{X\in\cd^{\mathbb{Z}}(A)\,|\,\Hom_{\cd^{\mathbb{Z}}(A)}(A(-i)[f(i)-k],X)=0,\,\text{for all $i\in\mathbb{Z}$ and $k<0$}\}\\
&=\{X\in\cd^{\mathbb{Z}}(A)\,|\,H^{k-f(i)}(X)_{i}=0,\,\text{for all $i\in\mathbb{Z}$ and $k<0$}\}\\
&=\{X\in\cd^{\mathbb{Z}}(A)\,|\,H^{m}(X)_{i}=0,\,\text{for all $i\in\mathbb{Z}$ and $f(i)<-m$}\}. \end{align*}
{
For $-c'<m\le -c$, we have $c\le -m<c'$. Hence $f(i)\le -m$ for $i\ll 0$ and
$f(i)>-m$ for $i\gg 0$, so $b_m$ is well-defined. If $-m=c$ and $c$ is finite,
this follows from the fact that the integer-valued monotone sequence
$f|_{\mathbb Z}$ is eventually constant with value $c$. The same argument gives
$a_m$ for $-c'\le m<-c$. If $c=-\infty$ or $c'=+\infty$, the corresponding
assertions follow directly from the limits.

By monotonicity, $f(i)>-m$ if and only if $i>b_m$, and $f(i)<-m$ if and only if
$i<a_m$. Hence the support descriptions follow. If $m>-c$, then $f(i)>-m$ for
all $i\in\mathbb Z$; if $m<-c'$, then $f(i)<-m$ for all $i\in\mathbb Z$. This
proves the vanishing statements. When $c=-\infty$ or $c'=+\infty$, the
corresponding vanishing condition is void.
}
\end{proof}
\begin{Lem}\label{lem:3.5}
Let $X'\ra X\ra X''\ra X'[1]$ be a triangle in $\cd^{\mathbb{Z}}(A)$
with $X'\in\cd_{\cm_{f}}^{\leqslant0}$ and $X''\in\cd_{\cm_{f}}^{\geqslant1}$. Then for each $i\in\mathbb{Z}$, we have the following short exact sequence
\[0\ra H^iX'\ra H^iX\ra H^iX''\ra 0\]
in $\Mod^{\mathbb{Z}}(H^0(A))$.
\end{Lem}
\begin{proof}

Since $X''\in\mathcal{D}_{\mathcal{M}_f}^{\geqslant 1} =\mathcal{D}_{\mathcal{M}_f}^{\geqslant 0}[-1]$, we have
$X''[1]\in\mathcal{D}_{\mathcal{M}_f}^{\geqslant 0}$. It is enough to show that the connecting homomorphism
$H^i(X'')\longrightarrow H^{i+1}(X')$ is zero for each $i\in\mathbb Z$.

By \Cref{lem:t-stru-onD(A)}, if $i-1<-c'$, then
$H^i(X'')=H^{i-1}(X''[1])=0$, while if $i+1>-c$, then
$H^{i+1}(X')=0$. Otherwise $-c'\leq i-1<-c$ and
$-c'<i+1\leq -c$, and hence
$H^i(X'')\in \Mod(H^0(A))_{\geq a_{i-1}}$ and
$H^{i+1}(X')\in \Mod(H^0(A))_{\leq b_{i+1}}$.
By the definitions of $a_{i-1}$ and $b_{i+1}$, we have
\[ f(a_{i-1})\geq -(i-1)>-(i+1)\geq f(b_{i+1}). \]
Since $f$ is non-decreasing, $a_{i-1}>b_{i+1}$. Thus the internal supports of
$H^i(X'')$ and $H^{i+1}(X')$ are disjoint. The connecting homomorphism is a
degree-zero morphism in $\Mod^{\mathbb Z}(H^0(A))$, and is therefore zero.
This gives the required short exact sequence.

\end{proof}

\bigskip

We make the following assumption$\colon$
\begin{assumption}\label{Ass1}
	\begin{itemize}
		\item[(A1)]
		The space $H^{0}(A)_{0}$ is finite dimensional.
		
		\item[(A2)] $\pvd^{\mathbb{Z}}(A)$ is contained in $\per^{\mathbb{Z}}(A)$.
	\end{itemize}
\end{assumption}	

The condition $(A2)$ is satisfied if $A$ is graded homologically smooth.

\begin{Lem}
\label{lem:automatic-bigraded-local-finiteness}
Assume that $A$ satisfies \textup{(A1)} and is graded homologically smooth.
Then $H^q(A)_p$ is finite-dimensional over $k$ for all $p,q\in\mathbb Z$.
\end{Lem}
\begin{proof}
Put $B=\operatorname{Tot}_2(A)$, with $B^n=\bigoplus_{q+2p=n}A_p^q$.
Since $d_A$ preserves Adams degree and $q+2p\equiv q\pmod 2$,
the differential and multiplication of $A$ give $B$ the structure
of a dg algebra. As $A_p^q=0$ unless $p,q\leq0$, $B$ is non-positive.
Taking cohomology in each Adams degree gives
\begin{equation}\label{eq:tot2-cohomology} H^n(B)\cong\bigoplus_{q+2p=n}H^q(A)_p. \end{equation}
In particular, $H^0(B)\cong H^0(A)_0$ is finite-dimensional by \textup{(A1)}.

The functor $\operatorname{Tot}_2$ extends naturally to differential bigraded $A$-bimodules. Since $2p$ is even, totalization leaves the Koszul signs in opposite algebras and tensor products unchanged, giving a natural isomorphism of dg algebras
\[\operatorname{Tot}_{2}(A^e)\cong\operatorname{Tot}_{2}(A)^e=B^e.\]
Under this identification, the exact functor $\operatorname{Tot}_2$ sends each generator $A^e(r)[s]$ to $B^e[2r+s]$, and hence sends $\per^{\mathbb Z}(A^e)$ into $\per(B^e)$. Since $A\in\per^{\mathbb Z}(A^e)$, we obtain $B=\operatorname{Tot}_2(A)\in\per(B^e)$.

Thus $B$ is homologically smooth, so $\pvd(B)\subseteq\per(B)$.
By \cite[Proposition~2.5]{Kalck-Yang2016}, each $H^n(B)$ is finite-dimensional.
Since $H^q(A)_p$ is a direct summand of $H^{q+2p}(B)$
by \Cref{eq:tot2-cohomology}, the result follows.
\end{proof}

\begin{Thm}\label{Thm:restriction-to-perfect}
Let $A=\bigoplus_{i,j\le0}A_i^j$ and $f$ be as above, and assume
\textup{(A1)} and \textup{(A2)}. The following are equivalent.
\begin{enumerate}
\item For every $P\in\per^{\mathbb Z}(A)$, one has $\tau_f^{\ge0}P\in\pvd^{\mathbb Z}(A)$.
\item\label{item:perfect-truncation-all-degrees} For every $P\in\per^{\mathbb Z}(A)$ and $n\in\mathbb Z$, one has $\tau_f^{\ge n}P\in\pvd^{\mathbb Z}(A)$.
\item\label{item:perfect-truncation-finite-sum} For every $r,n\in\mathbb Z$, one has
\[\sum_{\substack{p,q\le0\\q+f(p-r)\ge n}}\dim_k H^q(A)_p<\infty.\]
\item 

	One of the following two conditions holds:
		\begin{enumerate}
			\item[\rm(a)]
			$\lim_{i\to-\infty} f(i)=-\infty;$
			\item[\rm(b)]
			$\dim_k H^0(A) = \sum_{p\le0}\dim_k H^0(A)_p <\infty.$
		\end{enumerate}\end{enumerate}
		Under these conditions, the $t$-structure $(\mathcal D_{\cm_f}^{\le0},\mathcal D_{\cm_f}^{\ge0})$ restricts to $\per^{\mathbb Z}(A)$, and for every $n\in\mathbb Z$ we have
\[\per^{\mathbb Z}(A)\cap\mathcal D_{\cm_f}^{\ge n}\subseteq\pvd^{\mathbb Z}(A).\]
\end{Thm}
\begin{proof}
For $X\in\mathcal D^{\mathbb Z}(A)$ and $n\in\mathbb Z$, put
\[\lambda_{f,n}(X)=\sum_{q+f(p)\ge n}\dim_k H^q(X)_p.\]
By \Cref{eq:f-truncation-cohomology}, this is the total dimension of the cohomology of $\tau_f^{\ge n}X$. Hence
\[\tau_f^{\ge n}X\in\pvd^{\mathbb Z}(A)\quad\Longleftrightarrow\quad\lambda_{f,n}(X)<\infty.\]
Since $\tau_f^{\ge n}P\simeq\tau_f^{\ge0}(P[n])[-n]$, conditions \textup{(1)} and \textup{(2)} are equivalent.

The sum in \textup{(3)} is $\lambda_{f,n}(A(r))$, since $H^q(A(r))_i=H^q(A)_{i+r}$. Evaluating \textup{(2)} on the generators $A(r)$ gives \textup{(3)}. For the converse, the objects $X$ for which $\lambda_{f,n}(X)$ is finite for every $n$ form a thick subcategory: the condition is preserved by direct summands, $\lambda_{f,n}(X[s])=\lambda_{f,n+s}(X)$, and the long exact cohomology sequence of a triangle $X\to Y\to Z\to X[1]$ gives
\[\lambda_{f,n}(Z)\le\lambda_{f,n}(Y)+\lambda_{f,n+1}(X).\]
Under \textup{(3)}, this subcategory contains every $A(r)$ and hence all of $\per^{\mathbb Z}(A)$, proving \textup{(2)}.

Assume \textup{(3)} and suppose that \textup{(4)(a)} fails. Since $f$ is non-decreasing, it has an integer lower bound $c$. Taking $r=0$ and $n=c$ in \textup{(3)}, we obtain
\[\dim_k H^0(A)=\sum_{p\le0}\dim_k H^0(A)_p\le\sum_{\substack{p,q\le0\\q+f(p)\ge c}}\dim_k H^q(A)_p<\infty,\]
since $f(p)\ge c$ for every $p\le0$. Thus \textup{(4)(b)} holds.

For the converse, put $C=\operatorname{Tot}_2(A)$. As in the proof of \Cref{lem:automatic-bigraded-local-finiteness}, $C$ is non-positive and $H^0(C)\cong H^0(A)_0$ is finite-dimensional by \textup{(A1)}. Let $S$ be a simple right $H^0(C)$-module. Via $A\to H^0(A)_0$, regard $S$ as a bigraded $A$-module concentrated in bidegree $(0,0)$. Since $S$ is finite-dimensional, \textup{(A2)} gives
\[S\in\pvd^{\mathbb Z}(A)\subseteq\per^{\mathbb Z}(A).\]
Totalization preserves perfect objects, since it sends each generator $A(r)[s]$ to $C[2r+s]$. Hence $\operatorname{Tot}_2(S)=S$ belongs to $\per(C)$.

By \cite[Proposition~2.1(b)]{Kalck-Yang2016}, the finite-dimensional $H^0(C)$-modules generate $\pvd(C)$. Taking their composition series gives
\[\pvd(C)=\operatorname{thick}\{S\mid S\text{ is a simple right }H^0(C)\text{-module}\}\subseteq\per(C).\]
Thus \cite[Proposition~2.5]{Kalck-Yang2016} gives $\dim_k H^m(C)<\infty$ for every $m$. By \Cref{eq:tot2-cohomology}, we obtain
\begin{equation}\label{eq:local-bigraded-finiteness-from-A1-A2} \dim_k H^q(A)_p\le \dim_k H^{q+2p}(C)<\infty \qquad(p,q\in\mathbb Z). \end{equation}

Under \textup{(4)(a)}, fix $r,n\in\mathbb Z$. In the sum in
\textup{(3)}, the inequality $q\le0$ gives $f(p-r)\ge n$.
Since $f(i)\to-\infty$ as $i\to-\infty$, only finitely many $p\le0$
are possible. For each such $p$, one has $n-f(p-r)\le q\le0$.
Thus \textup{(3)} follows from \Cref{eq:local-bigraded-finiteness-from-A1-A2}.

Under \textup{(4)(b)}, let $A^{\mathrm u}$ denote $A$ with its Adams grading forgotten. The ideal $\bigoplus_{p<0}H^0(A)_p$ of the finite-dimensional algebra $H^0(A^{\mathrm u})$ is nilpotent. Every simple $H^0(A^{\mathrm u})$-module therefore factors through $H^0(A)_0$ and, placed in bidegree $(0,0)$, belongs to $\pvd^{\mathbb Z}(A)\subseteq\per^{\mathbb Z}(A)$ by \textup{(A2)}. Forgetting the Adams grading sends $A(r)[s]$ to $A^{\mathrm u}[s]$, so these simples are perfect over $A^{\mathrm u}$. By the same truncation and composition-series argument as for $C$, the simple $H^0(A^{\mathrm u})$-modules generate $\pvd(A^{\mathrm u})$, which implies $\pvd(A^{\mathrm u})\subseteq\per(A^{\mathrm u})$. Hence \cite[Proposition~2.5]{Kalck-Yang2016} gives
\[\sum_{p\le0}\dim_k H^q(A)_p=\dim_k H^q(A^{\mathrm u})<\infty \qquad(q\in\mathbb Z).\]
Fix $r,n\in\mathbb Z$. Since $f(p-r)\le f(-r)$ for $p\le0$, only the finitely many rows $n-f(-r)\le q\le0$ contribute to the sum in \textup{(3)}. Thus
\[\sum_{\substack{p,q\le0\\q+f(p-r)\ge n}}\dim_k H^q(A)_p\le\sum_{\substack{p\le0\\n-f(-r)\le q\le0}}\dim_k H^q(A)_p=\sum_{n-f(-r)\le q\le0}\dim_k H^q(A^{\mathrm u})<\infty.\]
This proves \textup{(3)}.

Finally, let $P\in\per^{\mathbb Z}(A)$ and $n\in\mathbb Z$. By \textup{(2)} and \textup{(A2)}, we have $\tau_f^{\ge n}P\in\pvd^{\mathbb Z}(A)\subseteq\per^{\mathbb Z}(A)$. The truncation triangle
\[\tau_f^{\le n-1}P\longrightarrow P\longrightarrow\tau_f^{\ge n}P\longrightarrow\tau_f^{\le n-1}P[1]\]
then shows that $\tau_f^{\le n-1}P$ is also perfect. Hence the $t$-structure restricts to $\per^{\mathbb Z}(A)$. If moreover $P\in\mathcal D_{\cm_f}^{\ge n}$, then $P\simeq\tau_f^{\ge n}P\in\pvd^{\mathbb Z}(A)$.
\end{proof}

\section{Cluster-tilting in graded Higgs category}\label{sec:graded-higgs}
\subsection{Relative $t$-structures and graded Higgs categories}\label{subsec:relative-t-structures}
\emergencystretch=3em 

Let $A=\bigoplus_{i,j\leq0}A_i^j$ be a dbg algebra over $k$, and let
$e=e_1+\cdots+e_n\in A_0^0$ be a sum of pairwise orthogonal primitive
idempotents. Let $f\colon\mathbb R\to\mathbb R$ be non-decreasing with $f(\mathbb Z)\subseteq\mathbb Z$. Its restriction to $\mathbb Z$ is therefore non-decreasing and integer-valued.
Put
\[\cm_f=\add\{A(-i)[f(i)]\mid i\in\mathbb Z\},\qquad \cp_f=\add\{eA(-i)[f(i)]\mid i\in\mathbb Z\}\subseteq\cm_f.\]
The functor $-\lten_{eAe}eA$ identifies $\per^{\mathbb Z}(eAe)$ with
$\thick(\cp_f)$.

Let $\overline A$ be the derived idempotent
quotient $A/^{\mathbf L}AeA$ from \Cref{Subsection:dg quotient}.
Let $\pi\colon\per^{\mathbb Z}(A)\to\per^{\mathbb Z}(\overline A)$ denote the
restriction of $p^*=-\lten_A \overline A$ to perfect objects. We retain the
standing assumption $\dim_k H^0(A)_0<\infty$ from
\Cref{Ass1}\textup{(A1)}.

\begin{assumption}\label{Ass}
\begin{itemize}
\item[\textup{(A0)}]
$A$ is graded homologically smooth, that is, $A\in\per^{\mathbb Z}(A^e)$.

\item[(A1)] The pair $(\per^{\mathbb{Z}}A,\pvd_{e}^{\mathbb{Z}}A)$ is relative twisted $(d + 1)$-Calabi--Yau of Gorenstein parameter $a$ in the sense that there exists a triangle equivalence $$\mathbb{S}=()_{\alpha}(a)[d+1]\colon\per^{\mathbb{Z}}(A)\ra\per^{\mathbb{Z}}(A)$$ such that it restricts to an
		autoequivalence on $\pvd_{e}^{\mathbb{Z}}(A)$ and there is a bifunctorial isomorphism for any $X\in\pvd_{e}^{\mathbb{Z}}A$ and $Y\in\per^{\mathbb{Z}}A$
		$$D\Hom_{\per^{\mathbb{Z}}A}(X,Y)\cong\Hom_{\per^{\mathbb{Z}}A}(Y,X_{\alpha}(a)[d+1]),$$ where $\alpha\colon A\ra A$ is an Adams graded dg automorphism which preserves the set $\{e_{1},\cdots,e_{n}\}$ and $X_\alpha$ has the same underlying complex as $X$, with right
action $x\cdot_\alpha b=x\cdot\alpha(b)$.

\item[\textup{(A2)}]
$H^0(A)$ is graded left and right Noetherian, and at least one of the
following conditions holds:
\begin{enumerate}
\item[\textup{(a)}] $\displaystyle\lim_{i\to-\infty}f(i)=-\infty$.
\item[\textup{(b)}] $\dim_k H^0(\overline A)<\infty$.
\end{enumerate}
\end{itemize}
\end{assumption}

\begin{Lem}\label{lem:graded-noetherian-functorial-finiteness}
Put $R=H^0(A)=\bigoplus_{d\le0}R_d$, where $R_d=H^0(A)_d$, and let $R_e=eRe$.
\begin{enumerate}
\item\label{item:frozen-finiteness-noetherian} If $Re$ (resp.\ $eR$) is finitely generated as a graded right (resp.\ left) $R_e$-module, then $\cp_f$ is contravariantly (resp.\ covariantly) finite in $\cm_f$. In particular, this holds if $R$ is graded right (resp.\ left) Noetherian.
\item If $f$ is strictly increasing on $\mathbb Z$, then $\cp_f$ is functorially finite in $\cm_f$, without any Noetherian hypothesis. This applies in particular to $f(p)=Np+c$ with $N>0$.
\end{enumerate}
\end{Lem}

\begin{proof}
Write $M_i=A(-i)[f(i)]$ and $P_i=eA(-i)[f(i)]$. All Hom spaces are taken in $\mathcal D^{\mathbb Z}(A)$. We have $\Hom(P_j,M_i)\cong H^{f(i)-f(j)}(A)_{j-i}e$. Since $A$ is concentrated in non-positive Adams and cohomological degrees, this space vanishes unless $j\le i$ and $f(i)\le f(j)$. By monotonicity of $f$, we obtain
\begin{equation}\label{eq:hom-Pj-Mi-graded-noetherian} \Hom(P_j,M_i)\cong\begin{cases}R_{j-i}e,&j\le i\text{ and }f(j)=f(i),\\
0,&\text{otherwise}.\end{cases} \end{equation}
Similarly,
\begin{equation}\label{eq:hom-Pl-Pj-graded-noetherian} \Hom(P_j,P_i)\cong\begin{cases}eR_{j-i}e,&j\le i\text{ and }f(j)=f(i),\\
0,&\text{otherwise}.\end{cases} \end{equation}
Under these identifications, composition is given by multiplication in $R$.

For \textup{(1)}, choose homogeneous generators $x_1,\ldots,x_m$ of $Re$ as a right $R_e$-module, with $x_s\in R_{d_s}e$. Fix $i$. By \Cref{eq:hom-Pj-Mi-graded-noetherian}, the generators satisfying $f(i+d_s)=f(i)$ define a morphism
\[u_i=(x_s)_s\colon\bigoplus_{\substack{1\le s\le m\\f(i+d_s)=f(i)}}P_{i+d_s}\longrightarrow M_i.\]
We claim that $u_i$ is a right $\cp_f$-approximation. Any nonzero morphism $P_j\to M_i$ is represented by an element $x\in R_{j-i}e$, with $j\le i$ and $f(j)=f(i)$. Write
\begin{equation}\label{eq:right-factorization-graded-noetherian} x=\sum_{s=1}^m x_s a_s,\qquad a_s\in eR_{j-i-d_s}e. \end{equation}
If $a_s\ne0$, then $j-i-d_s\le0$. Since $d_s\le0$, we have $j\le i+d_s\le i$, so $f(i+d_s)=f(i)$ by monotonicity. Thus $x_s$ is a component of $u_i$, and $a_s$ represents a morphism $P_j\to P_{i+d_s}$ by \Cref{eq:hom-Pl-Pj-graded-noetherian}. Hence \Cref{eq:right-factorization-graded-noetherian} gives the required factorization through $u_i$. Passing to finite direct sums and direct summands proves contravariant finiteness.

The left approximation follows dually, using homogeneous generators of $eR$ as a left $R_e$-module.
Finally, if $R$ is graded right Noetherian, finitely many homogeneous elements $x_1,\ldots,x_m\in Re$ generate $ReR$ as a right ideal, so $Re=(ReR)e=\sum_{s=1}^m x_sR_e$. The left Noetherian case is dual.

For \textup{(2)}, strict monotonicity gives $\Hom(P_j,M_i)=0=\Hom(M_i,P_j)$ whenever $j\ne i$. For $j=i$, these spaces are $R_0e$ and $eR_0$, respectively. Both are finite-dimensional by the standing assumption $\dim_k R_0<\infty$. A basis of the first space gives a right $\cp_f$-approximation of $M_i$, and a basis of the second gives a left one. Hence $\cp_f$ is functorially finite in $\cm_f=\add\{M_i\mid i\in\mathbb Z\}$.
\end{proof}

\begin{Ex}
Let $A=kQ$, where $Q$ is the quiver
\[ \begin{tikzcd} 1\arrow[r,"a"]&2&3\arrow[l,"b",swap] \end{tikzcd} \]
with zero differential and all generators in bidegree $(0,0)$.
Write $e_v$ for the vertex idempotents and put $e=e_1+e_3$.
Define a dg automorphism $\alpha$ by
\[ \alpha(e_1)=e_3,\quad \alpha(e_3)=e_1,\quad \alpha(e_2)=e_2,\qquad \alpha(a)=b,\quad \alpha(b)=a. \]
Then $\alpha^2=\id$ and $\alpha(e)=e$, although $\alpha$
interchanges the two frozen idempotents. In particular,
$(e_1A)_\alpha\cong e_3A$ and $(e_3A)_\alpha\cong e_1A$,
so $(-)_\alpha$ preserves both $\cm_f$ and $\cp_f$.

Put $X_{v,i}=e_vA(-i)[f(i)]$. Since $A$ is concentrated
in bidegree $(0,0)$, we have
\[ \Hom(X_{v,j},X_{w,i})\cong \begin{cases} e_wAe_v,&j=i,\\
0,&j\ne i. \end{cases} \]
In particular, $\Hom(X_{1,i},X_{2,i})=ka$ and
$\Hom(X_{3,i},X_{2,i})=kb$. Hence
\[ X_{1,i}\oplus X_{3,i} \xrightarrow{\ (a,b)\ }X_{2,i} \]
is a right $\cp_f$-approximation. Since
$e_1Ae_2=e_3Ae_2=0$, the map $X_{2,i}\to0$
is a left $\cp_f$-approximation. For the frozen summands $X_{1,i}$ and $X_{3,i}$, the identity maps give both approximations. These extend to $\cm_f$ by finite direct
sums and direct summands, as in \Cref{lem:graded-noetherian-functorial-finiteness}.

This example illustrates the permutation of frozen idempotents and the frozen approximations; it does not verify the relative Serre-duality condition in \Cref{Ass}(A1).
\end{Ex}

Let $\overline{\cm}_{f}=\add\{\overline{A}(-i)[f(i)]\,|\,i\in\mathbb{Z}\}$. Then $\overline{\cm}_{f}$ is a silting subcategory of $\cd^{\mathbb{Z}}(\overline{A})$. 
The algebra $\overline A$ is graded homologically smooth: indeed, extension of scalars on both sides sends the perfect diagonal of $A$ to that of $\overline A$, since $\overline A\lten_A \overline A\simeq \overline A$. Moreover, $H^0(\overline A)_0$ is a quotient of $H^0(A)_0$ and hence finite-dimensional. 
Thus $\overline A$ satisfies \Cref{Ass1}.
Moreover, the alternatives in \Cref{Ass}(A2)
verify condition~\textup{(4)} of
\Cref{Thm:restriction-to-perfect} for $\overline A$.

Let
$$\overline{t}^{\leqslant0}=\overline{\cm}_{f}[<\!0]^{\perp_{\cd^{\mathbb{Z}}(\overline{A})}}=\{X\in\cd^{\mathbb{Z}}(\overline{A})\,|\,H^{k-f(i)}(X)_{i}=0,\,\text{for all $k>0$ and $i\in\mathbb{Z}$}\},$$
$$\overline{t}^{\geqslant0}=\overline{\cm}_{f}[>\!0]^{\perp_{\cd^{\mathbb{Z}}(\overline{A})}}=\{X\in\cd^{\mathbb{Z}}(\overline{A})\,|\,H^{k-f(i)}(X)_{i}=0,\,\text{for all $k<0$ and $i\in\mathbb{Z}$}\}.$$ It follows from {\let\CrefLemmaSpace\space\Cref{lem:t-stru-onD(A)}} and \Cref{Thm:restriction-to-perfect} that the pair $(\overline{t}^{\leqslant0}, \overline{t}^{\geqslant0})$ is a $t$-structure on $\cd^{\mathbb{Z}}(\overline{A})$, which restricts to a $t$-structure on $\per^{\mathbb{Z}}(\overline{A})$. 

On $\cd^{\mathbb{Z}}(eAe)$ we take the
trivial $t$-structure whose left aisle is $\cd^{\mathbb{Z}}(eAe)$. We deduce the following Corollary from \Cref{prop:rec}.

\begin{Cor}\cite[Corollary 4.8]{Wu2023}
	There is a $t$-structure $(t_{rel}^{\leqslant0}, t_{rel}^{\geqslant0})$ on $\cd^{\mathbb{Z}}(A)$ obtained by gluing the $t$-structure $(\overline{t}^{\leqslant0}, \overline{t}^{\geqslant0})$ on $\cd^{\mathbb{Z}}(\overline{A})$ with
	the trivial t-structure on $\cd^{\mathbb{Z}}(eAe)$ through the recollement in \Cref{prop:rec}, where
	\[ t_{\mathrm{rel}}^{\leq0}=(p^*)^{-1}(\overline t^{\leq0}) \quad\text{and}\quad t_{\mathrm{rel}}^{\geq0}=p_*(\overline t^{\geq0}). \]
\end{Cor}
For each $k\in\mathbb{Z}$, the relative truncation functor is given by \cite[Definition 4.9]{Wu2023}
$$\tau^{>k}=\tau^{\geqslant k+1}\colon\cd^{\mathbb{Z}}(A)\xrightarrow{p^{*}}\cd^{\mathbb{Z}}(\overline{A})\xrightarrow{\overline{\tau}_{f}^{>k}}\cd^{\mathbb{Z}}(\overline{A})\xrightarrow{p_{*}}\cd^{\mathbb{Z}}(A).$$
For any $X\in \cd^{\mathbb{Z}}(A)$, we have a canonical triangle in $\cd^{\mathbb{Z}}(A)$
\[ \tau^{\leq k}X\longrightarrow X\longrightarrow\tau^{>k}X\longrightarrow\tau^{\leq k}X[1], \]
with $\tau^{\leq k}X\in t_{rel}^{\leq k}$ and $\tau^{>k}X\in t_{rel}^{\geq k+1}$.

\begin{Prop}\label{Prop:relative-t-on per}
The $t$-structure $(t_{\mathrm{rel}}^{\leq0}, t_{\mathrm{rel}}^{\geq0})$ on $\cd^{\mathbb Z}(A)$ restricts to $\per^{\mathbb Z}(A)$, and
\[ \tau^{\geq0}(\per^{\mathbb Z}(A)) \subseteq\pvd_e^{\mathbb Z}(A). \]
Moreover, the coaisle of the restricted $t$-structure is
\[ t_{\mathrm{rel}}^{\geq0}\cap\per^{\mathbb Z}(A)={}\cm_f[>\!0]^{\perp_{\per^{\mathbb Z}(A)}}\cap\thick(\cp_f)^{\perp_{\per^{\mathbb Z}(A)}}\cap{}^{\perp_{\per^{\mathbb Z}(A)}}\thick(\cp_f)\subseteq\pvd_e^{\mathbb Z}(A). \]
\end{Prop}

\begin{proof}
Let $X\in\per^{\mathbb Z}(A)$. Consider the relative truncation triangle
\[ \tau^{\leqslant-1}X\ra X\ra\tau^{\geqslant0}X \ra\tau^{\leqslant-1}X[1], \]
where $\tau^{\geqslant0}X=p_*\overline{\tau}_f^{>-1}(p^*X)$.
Since $p^*X\in\per^{\mathbb Z}(\overline A)$, \Cref{Thm:restriction-to-perfect} gives $\overline{\tau}_f^{>-1}(p^*X)\in\pvd^{\mathbb Z}(\overline A)$. Hence
\[ \tau^{\geqslant0}X\in p_*(\pvd^{\mathbb Z}(\overline A))=\pvd_e^{\mathbb Z}(A)\subseteq\per^{\mathbb Z}(A). \]
The triangle then gives $\tau^{\leqslant-1}X\in\per^{\mathbb Z}(A)$. Thus the $t$-structure restricts to $\per^{\mathbb Z}(A)$, with $\tau^{\geqslant0}(\per^{\mathbb Z}(A)) \subseteq\pvd_e^{\mathbb Z}(A)$.

We now prove the coaisle equality. Suppose that
\[ X\in\cm_f[>\!0]^{\perp_{\per^{\mathbb Z}(A)}} \cap\thick(\cp_f)^{\perp_{\per^{\mathbb Z}(A)}} \cap{}^{\perp_{\per^{\mathbb Z}(A)}}\thick(\cp_f). \]
Since $X\in\thick(\cp_f)^\perp$, \Cref{prop:rec} gives $X\simeq p_*p^*X$. For $M\in\cm_f$ and $r>0$, adjunction yields
\[ \Hom_{\per^{\mathbb Z}(\overline A)} (p^*M[r],p^*X) \simeq \Hom_{\per^{\mathbb Z}(A)}(M[r],X)=0. \]
As $p^*\cm_f=\overline\cm_f$, we have $p^*X\in\overline t^{\geq0}$, and therefore $X\in t_{\mathrm{rel}}^{\geq0} \cap\per^{\mathbb Z}(A)$.

Conversely, let $X\in t_{\mathrm{rel}}^{\geq0} \cap\per^{\mathbb Z}(A)$.
By adjunction, $X\in\cm_f[>\!0]^\perp$. The first part of the proof gives $X\in\pvd_e^{\mathbb Z}(A)$, so $\Hom(P,X)=0$ for every $P\in\thick(\cp_f)$. Moreover, \Cref{Ass}(A1) gives
\[ D\Hom_{\per^{\mathbb Z}(A)}(X,P)\cong \Hom_{\per^{\mathbb Z}(A)}(P,\mathbb S X)=0, \]
because $\mathbb S X\in\pvd_e^{\mathbb Z}(A) \subseteq\thick(\cp_f)^\perp$. Thus $X$ also belongs to ${}^\perp\thick(\cp_f)$, proving the claimed equality.
\end{proof}

\begin{Lem}[Dualizing criterion]\label{lem:dualizing-criterion}
Put $\overline M_p=\overline A(-p)[f(p)]$. Then $\overline{\cm}_f=\add\{\overline M_p\mid p\in\mathbb Z\}$ is a dualizing $k$-variety if and only if
\[ \dim_k H^0(\overline A)<\infty\quad\text{or}\quad \lim_{p\to-\infty}f(p)=-\infty\ \text{and}\ \lim_{p\to+\infty}f(p)=+\infty. \]
Under \Cref{Ass}, this is equivalent to
\begin{equation}\label{eq:reduced-dualizing-condition} \tag{D} \dim_k H^0(\overline A)<\infty\quad\text{or}\quad\lim_{p\to+\infty}f(p)=+\infty. \end{equation}
\end{Lem}

\begin{proof}
The third-quadrant condition and monotonicity of $f$ give
\begin{equation}\label{eq:reduced-silting-Hom} \Hom(\overline M_p,\overline M_q)= \begin{cases} H^0(\overline A)_{p-q},&p\leq q\text{ and }f(p)=f(q),\\
0,&\text{otherwise}. \end{cases} \end{equation}
Each graded component of $H^0(\overline A)$ is finite dimensional. If $\dim_k H^0(\overline A)<\infty$ or every fiber of $f$ is finite, \Cref{eq:reduced-silting-Hom} gives, for every $p\in\mathbb Z$,
\[ \sum_{q\in\mathbb Z}\dim_k\Hom(\overline M_p,\overline M_q)<\infty,\qquad \sum_{q\in\mathbb Z}\dim_k\Hom(\overline M_q,\overline M_p)<\infty. \]
Thus every finitely presented covariant or contravariant functor has finite total dimension. Conversely, let $F$ be a finite-dimensional contravariant functor. Since its support is finite, Yoneda gives an epimorphism
\[ \bigoplus_{\substack{p\in\mathbb Z\\
F(\overline M_p)\ne0}} \Hom(-,\overline M_p)^{\oplus\dim_k F(\overline M_p)} \twoheadrightarrow F. \]
The kernel is again finite dimensional, so applying the same construction to it yields a finite presentation of $F$. The covariant case is dual. Hence
\[ \operatorname{mod}\overline{\cm}_f=\operatorname{fd}\overline{\cm}_f,\qquad \operatorname{mod}\overline{\cm}_f^{\mathrm{op}}=\operatorname{fd}\overline{\cm}_f^{\mathrm{op}}. \]
Pointwise $k$-duality therefore induces a duality between the categories of finitely presented functors. Since $f$ is non-decreasing and integer-valued, all its fibers are finite if and only if
\[ \lim_{p\to-\infty}f(p)=-\infty\quad\text{and}\quad\lim_{p\to+\infty}f(p)=+\infty. \]
This proves sufficiency.

For necessity, suppose that $H^0(\overline A)$ is infinite dimensional. Since its graded components are finite dimensional and concentrated in non-positive degrees, $H^0(\overline A)_d\ne0$ for arbitrarily negative $d$. If $f$ is bounded above, it is constant on a right tail. Choose $p$ on this tail. Then
\[ D\Hom(\overline M_p,\overline M_q) \cong D H^0(\overline A)_{p-q}\ne0 \qquad\text{for arbitrarily large }q. \]
Thus the contravariant functor $D\Hom(\overline M_p,-)$ has unbounded-above support. However, every finitely generated contravariant functor is a quotient of a finite sum of representables $\Hom(-,\overline M_j)$, each supported only at indices $q\leq j$ by \Cref{eq:reduced-silting-Hom}. Hence $D\Hom(\overline M_p,-)$ is not finitely generated, contradicting dualizingness. The case where $f$ is bounded below is dual: on its constant left tail, $D\Hom(-,\overline M_p)$ has unbounded-below support, whereas every finitely generated covariant functor has support bounded below. Thus dualizingness forces both limits in the statement.

Finally, the equivalence with \Cref{eq:reduced-dualizing-condition} follows immediately from \Cref{Ass}\textup{(A2)}.
\end{proof}

\bigskip

The dualizing condition is not automatic. For example, take $A=\overline A=k[x,y]$ with zero differential, $\operatorname{bideg}(x)=\operatorname{bideg}(y)=(-1,0)$, $e=0$, and $f(p)=\min(p,0)$. The bimodule Koszul resolution gives
\[ \mathbf R\operatorname{Hom}^{\mathbb Z}_{A^e}(A,A^e)(2)[2]\simeq A, \]
where the dual carries the inner bimodule action. Thus \Cref{Ass}\textup{(A1)} holds with $a=2$, $d=1$, and $\alpha=\id$. Moreover, $A$ is graded homologically smooth, $H^0(A)_0=k$, $H^0(A)$ is graded Noetherian, and $\lim_{p\to-\infty}f(p)=-\infty$. Hence \Cref{Ass} is satisfied. However, $H^0(\overline A)$ is infinite dimensional and $f$ is bounded above, so \Cref{eq:reduced-dualizing-condition} fails. By \Cref{lem:dualizing-criterion}, $\overline{\cm}_f$ is not a dualizing $k$-variety.


\begin{Prop}\label{Prop:relative Serre five-tuple}
Under \Cref{Ass}, the five-tuple
$(\per^{\mathbb Z}(A),\pvd_e^{\mathbb Z}(A),\mathbb S,\cm_f,\cp_f)$
is a relative Serre five-tuple in the sense of \Cref{Def:Serre five-tuple} if and only if \Cref{eq:reduced-dualizing-condition} holds. The sufficiency also holds if the Noetherianity condition in \textup{(A2)} is replaced by the functorial finiteness of $\cp_f$ in $\cm_f$.
\end{Prop}

\begin{proof}
We verify the conditions of \Cref{Def:Serre five-tuple}. By \Cref{lem:automatic-bigraded-local-finiteness}, we have $\dim_k H^q(A)_p<\infty$ and $\dim_k H^q(\overline A)_p<\infty$ for all $p,q$.
Hence $\per^{\mathbb Z}(A)$ and $\per^{\mathbb Z}(\overline A)$ are Hom-finite, and the former is Krull--Schmidt. 

Put $\cu=\per^{\mathbb Z}(A)/\thick(\cp_f)$. Compact localization gives $\Kar(\cu)\simeq\per^{\mathbb Z}(\overline A),$
so $\cu$ is also Hom-finite. The functorial finiteness of $\cp_f$ in $\cm_f$ follows from \Cref{lem:graded-noetherian-functorial-finiteness}, or is assumed directly in the alternative formulation. This proves \textup{(RS1)}.

For \textup{(RS2)}, $X\in\pvd_e^{\mathbb Z}(A)$ implies $Xe\simeq0$, hence $\Hom(\thick(\cp_f),X)=0$. Condition \textup{(RS3)} is precisely \Cref{Ass}\textup{(A1)}.

For \textup{(RS4)}, let $q_f:\per^{\mathbb Z}(A)\to\cu$ be the quotient functor. By compact localization, $\cu$ embeds fully faithfully into $\per^{\mathbb Z}(\overline A)$. By the Hom calculation in \Cref{lem:graded-noetherian-functorial-finiteness}, the map on Hom spaces between generators of $\cm_f$ is
\[H^0(A)_{i-j}\longrightarrow H^0(\overline A)_{i-j}=\bigl(H^0(A)/H^0(A)eH^0(A)\bigr)_{i-j}\]
when $i\le j$ and $f(i)=f(j)$; otherwise both Hom spaces vanish. Every element of the kernel factors through a finite sum of objects $eA(-k)[f(k)]$ with $i\le k\le j$, since the Adams grading is non-positive and $f$ is non-decreasing. Hence
\[\cm_f/[\cp_f]\xrightarrow{\sim}q_f(\cm_f).\]
The subcategory $q_f(\cm_f)$ is silting in $\cu$, since its image $\overline{\cm}_f$ is presilting and it thickly generates $\cu$. Moreover, $\cm_f/[\cp_f]$ is Krull--Schmidt, hence $q_f(\cm_f)$ is idempotent complete. By \cite[Theorem~2.9]{IY2018}, $\cu$ is idempotent complete. Thus
\[\cu\xrightarrow{\sim}\per^{\mathbb Z}(\overline A),\qquad q_f(\cm_f)\simeq\overline{\cm}_f.\]

By \Cref{Thm:restriction-to-perfect}, the equivalence $\cu\simeq\per^{\mathbb Z}(\overline A)$ gives the required $t$-structure on $\cu$, with $\cu^{\geq0}\subseteq\pvd^{\mathbb Z}(\overline A)$. Moreover, \Cref{Prop:relative-t-on per} and the relative truncation formula give
\[\cu^{\geq0}\subseteq q_f\bigl(\pvd_e^{\mathbb Z}(A)\bigr).\]
Finally, \Cref{lem:dualizing-criterion} and \Cref{eq:reduced-dualizing-condition} show that $\overline{\cm}_f$ is a dualizing $k$-variety. This proves \textup{(RS4)}.

Condition \textup{(RS5)} follows from \Cref{Prop:relative-t-on per}. Conversely, the dualizing requirement in \Cref{Def:Serre five-tuple}, together with \Cref{lem:dualizing-criterion}, forces \Cref{eq:reduced-dualizing-condition}.

Finally, Noetherianity is used only to obtain the functorial finiteness of $\cp_f$ in $\cm_f$ via \Cref{lem:graded-noetherian-functorial-finiteness}, so it can be replaced by that condition.
\end{proof}

\begin{Def}\rm\label{Def:grcc-ghiggs}
The \emph{graded relative cluster category} of $(A,e)$ is defined as the idempotent completion of the Verdier quotient
\[\cc^{\mathbb Z}(A,e)=\Kar\bigl(\per^{\mathbb Z}(A)/\pvd_e^{\mathbb Z}(A)\bigr).\]
The \emph{graded Higgs category} is the full subcategory
\[\ch_{f}^{\mathbb{Z}}(A,e)=\{X\in\cc^{\mathbb{Z}}(A,e)\,|\,\Hom_{\cc^{\mathbb{Z}}(A,e)}(X,\cp_{f}[>0])=0=\Hom_{\cc^{\mathbb{Z}}(A,e)}(\cp_{f},X[>0])\},\]
where $\cp_f$ denotes its image in $\cc^{\mathbb Z}(A,e)$.
The \emph{graded relative fundamental domain} $\cf_f^{\mathbb Z}(A,e)$ is the full subcategory of $\per^{\mathbb Z}(A)$ defined by
\[\cf_f^{\mathbb Z}(A,e)=\cx^\perp\cap{}^\perp\cy\cap\cz,\]
where
\[\cx=t_{\mathrm{rel}}^{\leq1}\cap\pvd_e^{\mathbb Z}(A),\qquad \cy=t_{\mathrm{rel}}^{\geq1},\qquad \cz={}^\perp\cp_f[>0]\cap\cp_f[<0]^\perp.\]
Here the relative $t$-structure is restricted to $\per^{\mathbb Z}(A)$, and all orthogonals are taken in this category.
The \emph{graded cluster category} of $\overline A$ is defined as
\[\cc^{\mathbb Z}(\overline A)=\Kar\bigl(\per^{\mathbb Z}(\overline A)/\pvd^{\mathbb Z}(\overline A)\bigr).\]
\end{Def}

\begin{Rem}
Let $f(i)=0$ for each $i\in\mathbb{Z}$. Then the graded Higgs category is given by
\begin{equation} \begin{split} \ch_{0}^{\mathbb{Z}}(A,e)=\{X\in\cc^{\mathbb{Z}}(A,e)\,&|\,\Hom_{\cc^{\mathbb{Z}}(A,e)}(X,eA(i)[>0])=0\\
&\quad\Hom_{\cc^{\mathbb{Z}}(A,e)}(eA(i),X[>0])=0,\,\forall i\in\mathbb{Z}\}. \end{split} \end{equation}
Suppose in addition that $A$ is quasi-isomorphic to the graded algebra $R=H^0(A)$ and that $R$ is graded left and right Noetherian. Smoothness then gives finite graded global dimension, and hence
$\per^{\mathbb Z}(A)\simeq\cd^b(\mathrm{mod}^{\mathbb Z}R)$.
Put ${R_e}=eRe$. If, moreover, ${R_e}$ is graded Iwanaga--Gorenstein and the restriction functor $(-)e$ induces a triangle equivalence
\[ \cc^{\mathbb Z}(A,e)\xrightarrow{\sim}\cd^b(\mathrm{mod}^{\mathbb Z}{R_e}), \]
then $eA(i)$ corresponds to $R_e(i)$. Thus the graded version of \cite[Theorem 3.10]{IY2018} identifies the graded Higgs category with the category of graded Cohen--Macaulay ${R_e}$-modules:
\[ \ch_{0}^{\mathbb Z}(A,e)\simeq \mathrm{CM}^{\mathbb Z}({R_e}). \]
\end{Rem}

\begin{Prop}\label{Prop:rf^Z=add}
Assume \Cref{eq:reduced-dualizing-condition} and that $f(p)-f(p-a)=r$ for all $p\in\mathbb Z$, and put $n=d+r\geq1$. Then
\[\cf_f^{\mathbb Z}(A,e)=(\cm_f*\cm_f[1]*\cdots*\cm_f[n-1])\cap\cz.\]
In particular, $\cf_f^{\mathbb Z}(A,e)=\cm_f$ when $n=1$.
\end{Prop}

\begin{proof}
Put $\mathcal W_n=\cm_f*\cm_f[1]*\cdots*\cm_f[n-1]$ and $\mathbb S_{n+1}=\mathbb S[-n-1]$. The compatibility condition gives
\[\mathbb S_{n+1}\bigl(A(-p)[f(p)]\bigr)\simeq A(-(p-a))[f(p-a)].\]
Together with $\alpha(e)=e$, this shows that $\mathbb S_{n+1}$ preserves $\cm_f$ and $\cp_f$, and hence the reduced $t$-structure. Since it also preserves $\pvd_e^{\mathbb Z}(A)$, it preserves $\cx$.

Let $T\in\mathcal W_n\cap\cz$. Since $\cy\subseteq\cm_f[\geq0]^\perp$, we have $T\in{}^\perp\cy$. The definition of $\cx$ and the Hom comparison under $\pi$ for targets in $\pvd_e^{\mathbb Z}(A)$ give $\cx\subseteq\cm_f[<-1]^\perp$. For $X\in\cx$, relative Serre duality therefore gives
\[D\Hom_{\per^{\mathbb Z}(A)}(X,T)\cong\Hom_{\per^{\mathbb Z}(A)}(T[-n-1],\mathbb S_{n+1}X)=0,\]
because $\mathbb S_{n+1}X\in\cx$ and $T[-n-1]\in\cm_f[-n-1]*\cdots*\cm_f[-2]$. Thus $T\in\cx^\perp\cap{}^\perp\cy\cap\cz=\cf_f^{\mathbb Z}(A,e)$.

Conversely, let $U\in\cf_f^{\mathbb Z}(A,e)$. By \Cref{Lem:F/P=F_U}, its image $\pi U$ belongs to $\cf_{\cu}$. Since $\mathbb S_{n+1}(\overline{\cm}_f)=\overline{\cm}_f$, the argument in \cite[proof of Theorem~2.12, case~(i)]{IY2020}, shifted by $[-1]$, together with \cite[Proposition~2.8(b)]{IY2018}, gives
\[\cf_{\cu}=\overline{\cm}_f*\overline{\cm}_f[1]*\cdots*\overline{\cm}_f[n-1].\]
The lifting construction in the proof of \Cref{Lem:F/P=F_U} preserves this range of shifts, since its added frozen terms lie in $\cp_f\subseteq\cm_f$. It therefore gives an object
\[U'\in\mathcal W_n\cap\cz,\qquad \pi U'\simeq\pi U.\]

By the first inclusion, $U'\in\cf_f^{\mathbb Z}(A,e)$. The equivalence $\cf_f^{\mathbb Z}(A,e)/[\cp_f]\simeq\cf_{\cu}$ of \Cref{Lem:F/P=F_U} then shows that $U$ and $U'$ are isomorphic modulo $[\cp_f]$. Consequently, $U$ is a direct summand of $U'\oplus P$ for some $P\in\cp_f$. Since $\mathcal W_n$ contains $\cp_f$ and is closed under direct summands by \cite[Proposition~2.7(a)]{IY2018}, we obtain $U\in\mathcal W_n\cap\cz$. This proves the equality.

Finally, $\cm_f\subseteq\cz$ by the presilting property, so the assertion for $n=1$ follows.
\end{proof}

We also write $\pi$ for the induced functor on cluster categories. We have the following commutative diagram
\[\begin{tikzcd} &\thick^{\mathbb{Z}}(eA)\arrow[r,equal]\arrow[d,hook]&\thick^{\mathbb{Z}}(eA)\arrow[d,hook]\\
\pvd_{e}^{\mathbb{Z}}(A)\arrow[d,equal]\arrow[r,hook]&\per^{\mathbb{Z}}(A)\arrow[r,two heads,""]\arrow[d,two heads,"\pi"] &\cc^{\mathbb{Z}}(A,e)\arrow[d,two heads,"\pi"]\\
\pvd^{\mathbb{Z}}(\uA)\arrow[r,hook]&\per^{\mathbb{Z}}(\uA)\arrow[r,two heads,""]&\cc^{\mathbb{Z}}(\uA). \end{tikzcd}\]
\begin{Thm}\label{Thm:CTC-in-gradedHiggs}
Assume \Cref{eq:reduced-dualizing-condition}, in addition to \Cref{Ass}.
\begin{itemize}
\item[(1)] The graded Higgs category $\ch_f^{\mathbb Z}(A,e)$ is a Frobenius extriangulated category with projective--injective subcategory $\cp_f$. Its stable category admits a triangle equivalence
\[\ch_f^{\mathbb Z}(A,e)/[\cp_f]\xrightarrow{\sim}\cc^{\mathbb Z}(\overline A)\]
and a Serre functor given by $\mathbb S[-1]=(-)_\alpha(a)[d]$.
\item[(2)] If there exists $r\in\mathbb Z$ such that $f(i)-f(i-a)=r$ for all $i\in\mathbb Z$, and $d+r\geq1$, then $\cm_f$ is a $(d+r)$-cluster-tilting subcategory of $\ch_f^{\mathbb Z}(A,e)$.
\end{itemize}
Both assertions remain valid if graded left and right Noetherianity in \Cref{Ass}\textup{(A2)} is replaced by the functorial finiteness of $\cp_f$ in $\cm_f$, with all other assumptions unchanged.
\end{Thm}

\begin{proof}
By \Cref{Prop:relative Serre five-tuple}, the five-tuple $(\per^{\mathbb Z}(A),\pvd_e^{\mathbb Z}(A),\mathbb S,\cm_f,\cp_f)$ is a relative Serre five-tuple. Since $\alpha$ preserves the set of frozen primitive idempotents, $\mathbb S^{\pm1}(\cp_f)\subseteq\thick(\cp_f)$. Thus condition~\textup{(a)} of \Cref{Thm:H is Frobenius and H/P=U/Ufd} holds, giving \textup{(1)}.

For \textup{(2)}, we have $\mathbb S_{d+r+1}=\mathbb S[-d-r-1]=(-)_\alpha(a)[-r]$. Since $A_\alpha\simeq A$, the compatibility condition gives
\[\mathbb S_{d+r+1}\bigl(A(-i)[f(i)]\bigr)\simeq A(-(i-a))[f(i-a)].\]
Hence $\mathbb S_{d+r+1}(\cm_f)=\cm_f$, and \textup{(2)} follows from the same theorem.

Noetherianity is used only through \Cref{Prop:relative Serre five-tuple}, which also covers the alternative hypothesis.
\end{proof}

\begin{Rem}
In \Cref{Thm:CTC-in-gradedHiggs}\textup{(2)}, condition~\Cref{eq:reduced-dualizing-condition} is automatic when $a\ne0$. Indeed, if $r\ne0$, iteration of the compatibility condition together with monotonicity gives
\[\lim_{p\to-\infty}f(p)=-\infty,\qquad \lim_{p\to+\infty}f(p)=+\infty.\]
If $r=0$, then $f$ is periodic and non-decreasing, hence constant, so \Cref{Ass}\textup{(A2)} forces $\dim_k H^0(\overline A)<\infty$.

Under the hypotheses of \Cref{Thm:CTC-in-gradedHiggs}, we obtain the following special cases.
\begin{itemize}
\item[(1)] If $e=0$, then $\cp_f=0$, $\overline A=A$, and
\[\ch_f^{\mathbb Z}(A,0)=\cc^{\mathbb Z}(A).\]
Thus $\cc^{\mathbb Z}(A)$ has Serre functor $(-)_\alpha(a)[d]$, and $\cm_f$ is $(d+r)$-cluster-tilting whenever $f(i)-f(i-a)=r$ for all $i\in\mathbb Z$ and $d+r\geq1$.

\item[(2)] For $f(i)=i$, one has $r=a$. Thus, if $d+a\geq1$, then $\cm_f=\add\{A(-i)[i]\mid i\in\mathbb Z\}$ is $(d+a)$-cluster-tilting. For $e=0$ and $A$ a graded algebra satisfying the hypotheses of \cite[Theorem~4.6]{Hanihara2022}, this recovers part~\textup{(4)} of that theorem.

\item[(3)] For $f=0$, one has $r=0$, and \Cref{Ass}\textup{(A2)} gives $\dim_k H^0(\overline A)<\infty$. Hence \Cref{eq:reduced-dualizing-condition} holds, and $\cm_0=\add\{A(-i)\mid i\in\mathbb Z\}$ is $d$-cluster-tilting when $d\geq1$.
\end{itemize}
\end{Rem}

\subsection{Silting objects in the stable graded Higgs categories}\label{subsec:stable-graded-higgs}\
Assume $a>0$ and write $B=\overline A=\bigoplus_{p\le0}B_p$. Let $\bar\alpha$ be the dg automorphism of $B$ induced by $\alpha$. The relative Serre functor on $(\per^{\mathbb Z}(B),\pvd^{\mathbb Z}(B))$ induced by \Cref{Ass}\textup{(A1)} is $(-)_{\bar\alpha}(a)[d+1]$.
\begin{Lem}\label{lem:Adams-zero-part}
Under \Cref{Ass}\textup{(A0)}, suppose that $B_0$ is proper. Then $B_0$ is homologically smooth and perfect as both a right and a left $B$-module via $B\to B_0$. Moreover,
\[\pvd^{\mathbb Z}(B)=\thick\{B_0(i)\mid i\in\mathbb Z\}.\]
\end{Lem}

\begin{proof}
The algebra $B$ is graded homologically smooth, as established in \Cref*{subsec:relative-t-structures}. Consider the exact functor
\[F_0(M):=\bigl(B_0\lten_B M\lten_B B_0\bigr)_0\colon\cd^{\mathbb Z}(B^e)\longrightarrow\cd(B_0^e).\]
It sends $B^e(i)$ to $0$ for $i\ne0$ and to $B_0^e$ for $i=0$. Since $B$ is non-positively Adams graded, the Adams-degree-zero part of the two-sided bar model for $B_0\lten_B B_0$ is the usual bar model over $B_0$. Thus $F_0(B)\simeq B_0$, and $B\in\per^{\mathbb Z}(B^e)$ gives $B_0\in\per(B_0^e)$.

Properness of $B_0$ and graded homological smoothness of $B$ give the asserted right and left perfectness. Since $B_0$ is smooth and proper, $\pvd^{\mathbb Z}(B_0)=\per^{\mathbb Z}(B_0)$. Restriction along $B_0\hookrightarrow B$ preserves the underlying bigraded complex, so $\pvd^{\mathbb Z}(B)$ is precisely the relative perfectly valued
subcategory $\pvd^{\mathbb Z}(B,B_0)$ of
\cite[Definition~4.1]{FKQ2024}. Hence \cite[Lemma~4.2]{FKQ2024} gives the claimed equality.
\end{proof}

\bigskip

Under the hypotheses of \Cref{lem:Adams-zero-part}, put
\[\cc^{\mathbb Z}_{B_0}(B):=\per^{\mathbb Z}(B)/\thick\{B_0(i)\mid i\in\mathbb Z\},\qquad \cc_{B_0}(B):=\per(B)/\thick(B_0).\]
Here $\per(B)$ is formed after forgetting the Adams grading, and neither quotient includes idempotent completion. The lemma gives
\[\cc^{\mathbb Z}(B)\simeq\Kar\!\left(\cc^{\mathbb Z}_{B_0}(B)\right).\]
Forgetting the Adams grading induces a functor $\cc^{\mathbb Z}_{B_0}(B)\to\cc_{B_0}(B)$.

\begin{Thm}\label{Thm:finite-Adams-window}\cite[Theorem~A.7]{Hanihara24}
Under conditions \textup{(A0)} and \textup{(A1)} of \Cref{Ass}, suppose that $B_0$ is proper. Put $T=\bigoplus_{l=0}^{a-1}B(-l)$, and define the dg algebra $\Lambda$ and the $\Lambda$-$\Lambda$-bimodule $U$ by
\[\Lambda=\REnd_B(T)=\begin{pmatrix}B_0&0&\cdots&0\\B_{-1}&B_0&\cdots&0\\\vdots&\vdots&\ddots&\vdots\\B_{-(a-1)}&B_{-(a-2)}&\cdots&B_0\end{pmatrix},\]
\[U=\RHom_B(T,T(-1))=\begin{pmatrix}B_{-1}&B_0&\cdots&0\\B_{-2}&B_{-1}&\ddots&\vdots\\\vdots&\vdots&\ddots&B_0\\B_{-a}&B_{-(a-1)}&\cdots&B_{-1}\end{pmatrix}.\]
Then the following hold.
\begin{enumerate}
\item The functor $-\lten_\Lambda T$ induces a triangle equivalence $\per\Lambda\xrightarrow{\sim}\cc^{\mathbb Z}(B)$.
\item There is a semi-orthogonal decomposition
\[\cc^{\mathbb Z}(B)=\thick(B(-(a-1)))\perp\cdots\perp\thick(B(-1))\perp\thick(B).\]
\item We have the following commutative diagram, where $\mathrm{For}$ is the forgetful functor and $(-)_{\triangle}$ the triangulated hull:
\[\begin{tikzcd}\cc^{\mathbb Z}(B)\arrow[r,"(-1)"]\arrow[d,"\simeq"{rotate=90,anchor=south}]&\cc^{\mathbb Z}(B)\arrow[r,"\mathrm{For}"]\arrow[d,"\simeq"{rotate=90,anchor=south}]&\cc_{B_0}(B)\arrow[d,"\simeq"{rotate=90,anchor=south}]\\\per\Lambda\arrow[r,"-\lten_\Lambda U"]&\per\Lambda\arrow[r]&\bigl(\per\Lambda/(-\lten_\Lambda U)\bigr)_\triangle.\end{tikzcd}\]
\end{enumerate}
\end{Thm}

\begin{proof}
By \Cref{lem:Adams-zero-part}, $B_0$ is perfect as both a right and a left $B$-module. Since $B_0$ restricts to an object of $\pvd_e^{\mathbb Z}(A)$, localization preserves the relevant Hom spaces by two-sided orthogonality to $\thick(\cp_f)$. For $t,s\in\mathbb Z$, the induced relative Serre duality gives
\[D\Hom_{\cd^{\mathbb Z}(B)}(B_0,B(t)[s])\cong H^{d+1-s}\!\left((B_0)_{\bar\alpha}\right)_{a-t}.\]
The right-hand side vanishes for $t\ne a$, so $\RHom_B^{\mathbb Z}(B_0,B)$ has cohomology concentrated in Adams degree $a$. After reversing the Adams grading, \cite[Lemma~A.4(2)]{Hanihara24} therefore gives Gorenstein parameter $a$, and the hypotheses of \cite[Setting~A.5]{Hanihara24} are satisfied.

After reversing the Adams grading, the object in Hanihara's theorem becomes $T$, and his Adams shift $(1)$ becomes $(-1)$. Thus \cite[Theorem~A.7(1)--(2)]{Hanihara24} gives
\[\cc^{\mathbb Z}_{B_0}(B)\simeq\per\Lambda\]
and the semi-orthogonal decomposition in \textup{(2)}. Since $\per\Lambda$ is idempotent complete, the equivalence above and \Cref{lem:Adams-zero-part} give
\[\cc^{\mathbb Z}(B)\simeq\Kar\!\left(\cc^{\mathbb Z}_{B_0}(B)\right)\simeq\cc^{\mathbb Z}_{B_0}(B).\] Finally, \cite[Theorem~A.7(3)]{Hanihara24} identifies $(-1)$ with $-\lten_\Lambda U$ and gives the diagram in \textup{(3)}.
\end{proof}
\begin{Rem}
Under the hypotheses of \Cref{Thm:finite-Adams-window}, if $\dim_k H^0(B)<\infty$, then
\[\pvd(B)=\thick(B_0),\qquad \cc_{B_0}(B)=\per(B)/\pvd(B).\]
Indeed, the ideal $\bigoplus_{p<0}H^0(B)_p$ of $H^0(B)$ is nilpotent, so every simple $H^0(B)$-module factors through $H^0(B_0)$. Such a module is finite-dimensional and hence perfect over the smooth dg algebra $B_0$. Viewed as a $B$-module via $B\to B_0$, it therefore belongs to $\thick(B_0)$. Standard truncations and composition series, as in the proof of \Cref{Thm:restriction-to-perfect}, give $\pvd(B)\subseteq\thick(B_0)$; properness of $B_0$ gives the reverse inclusion.

With the idempotent-completion convention of \Cref{Def:grcc-ghiggs}, the corresponding ungraded cluster category is $\Kar(\cc_{B_0}(B))$.
\end{Rem}
Let $Q=(Q_0,Q_1)$ be a finite acyclic quiver. Its framed quiver $Q^{\mathrm{fr}}$ is obtained by adding a vertex $i'$ and an arrow $i'\to i$ for each $i\in Q_0$. The associated ice quiver is $(Q^{\mathrm{fr}},F)$, where $F$ is the discrete subquiver on the new vertices.

\begin{Ex}
For the quiver $Q\colon 1\to2\to3$, the associated ice quiver is
\[\begin{tikzcd}1\arrow[r]&2\arrow[r]&3\\\color{blue}\boxed{1'}\arrow[u]&\color{blue}\boxed{2'}\arrow[u]&\color{blue}\boxed{3'}\arrow[u]\end{tikzcd}\]
\end{Ex}
Write $\eta_i\colon i'\to i$ for the framing arrow and put $e=\sum_{i\in Q_0}e_{i'}$. The relative $2$-Calabi--Yau completion of $kF\hookrightarrow kQ^{\mathrm{fr}}$ is the dg functor
\[\phi\colon\Pi_1(kF)=\prod_{i\in Q_0}k[x_{i'}]\longrightarrow\Gamma_Q^{\mathrm{fr}}:=\Pi_2(kQ^{\mathrm{fr}},kF).\]
Its target is obtained from $kQ^{\mathrm{fr}}$ by adjoining a reverse arrow $\beta^*$ for every arrow $\beta$ and a loop $t_i$ at every unfrozen vertex $i\in Q_0$; see \cite[Sections~6--7]{KellerWang2023}. We use the bidegrees
\[\operatorname{bideg}(\beta)=(0,0),\qquad \operatorname{bideg}(\beta^*)=\operatorname{bideg}(x_{i'})=(-1,0),\qquad \operatorname{bideg}(t_i)=(-1,-1).\]
The differential vanishes on $\beta$ and $\beta^*$ and satisfies
\begin{equation}\label{eq:relative-preprojective-differential}d(t_i)=e_i\left(\sum_{\beta\in(Q^{\mathrm{fr}})_1}[\beta,\beta^*]\right)e_i.\end{equation}
The functor $\phi$ fixes the frozen vertices and sends
\[\phi(x_{i'})=e_{i'}\left(\sum_{\beta\in(Q^{\mathrm{fr}})_1}[\beta,\beta^*]\right)e_{i'}.\]

Put $A=\Gamma_Q^{\mathrm{fr}}$. The graded relative Calabi--Yau structure give
\[H^0(A)_0=kQ^{\mathrm{fr}},\qquad a=d=1,\qquad \mathbb S=(1)[2],\]
with the Gorenstein-parameter convention of \Cref{rem:potential-weight-versus-Gorenstein-parameter}. 
In particular, the dg quotient by the frozen vertices is
\begin{equation}\label{eq:framed-quotient-derived-preprojective} B:=\overline{\Gamma_Q^{\mathrm{fr}}} \simeq\Pi_2(kQ). \end{equation}
It has the same presentation with the frozen vertices and the arrows $\eta_i,\eta_i^*$ deleted, and moreover,
\begin{equation}\label{eq:degree-zero-derived-preprojective} B_0=kQ. \end{equation}

\begin{Thm}\label{Ex:acyclic-diagonal-grading}
Let $Q$ be a finite acyclic quiver, with $A=\Gamma_Q^{\mathrm{fr}}$, $B=\Pi_2(kQ)$, and $e$ as above. Take $f(p)=p$, put $\ch_{\mathrm{diag}}^{\mathbb Z}(Q):=\ch_f^{\mathbb Z}(\Gamma_Q^{\mathrm{fr}},e)$, and set
\[\cp_{\mathrm{diag}}=\add\{e\Gamma_Q^{\mathrm{fr}}(-p)[p]\mid p\in\mathbb Z\},\qquad \cm_{\mathrm{diag}}=\add\{\Gamma_Q^{\mathrm{fr}}(-p)[p]\mid p\in\mathbb Z\}.\]
Let $q\colon\per^{\mathbb Z}(B)\to\cc^{\mathbb Z}(B)$ be the canonical quotient functor. Then the following hold.
\begin{enumerate}
\item The $t$-structure associated with $f$ restricts to $\per^{\mathbb Z}(B)$, with
\[\overline{\tau}_f^{\ge0}\bigl(\per^{\mathbb Z}(B)\bigr)\subseteq\pvd^{\mathbb Z}(B).\]

\item The object $qB$ is a silting generator of $\cc^{\mathbb Z}(B)$. In particular,
\[\cc^{\mathbb Z}(B)=\thick(qB).\]

\item In the canonical dg enhancement of $\cc^{\mathbb Z}(B)$, there is a quasi-isomorphism of dg algebras
\[\REnd_{\cc^{\mathbb Z}(B)}(qB)\simeq kQ.\]
There is a dg Morita equivalence inducing the triangle equivalence
\begin{equation}\label{eq:acyclic-derived-Morita}q\circ\bigl(-\lten_{kQ}B\bigr)\colon\per(kQ)\xrightarrow{\sim}\cc^{\mathbb Z}(B).\end{equation}
Consequently, $\cc^{\mathbb Z}(B)\simeq \cd^b(\operatorname{mod}kQ)$.

\item The category $\ch_{\mathrm{diag}}^{\mathbb Z}(Q)$ is a Frobenius extriangulated category with projective--injective subcategory $\cp_{\mathrm{diag}}$ and $2$-cluster-tilting subcategory $\cm_{\mathrm{diag}}$. There are triangle equivalences
\begin{equation}\label{eq:stable-diagonal-Higgs-derived-Q}\ch_{\mathrm{diag}}^{\mathbb Z}(Q)/[\cp_{\mathrm{diag}}]\simeq\cc^{\mathbb Z}(B)\simeq \cd^b(\operatorname{mod}kQ).\end{equation}
Its Serre functor is $(1)[1]$, and the Adams shift $(1)$ on $\cc^{\mathbb Z}(B)$ corresponds to the Auslander--Reiten translation $\tau$ of $\cd^b(\operatorname{mod}kQ)$.
\end{enumerate}
\end{Thm} 

\begin{proof}
The relative completion construction and the quotient argument in \Cref*{subsec:relative-t-structures} give the graded homological smoothness of $A$ and $B$. Since $H^0(B)_0=kQ$ is finite-dimensional and $f(p)=p$ tends to $-\infty$ as $p\to-\infty$, \Cref{Thm:restriction-to-perfect} gives \textup{(1)}.

For \textup{(2)} and \textup{(3)}, we have $a=1$ and $B_0=kQ$ is proper. \Cref{Thm:finite-Adams-window} therefore applies with $T=B$ and $\Lambda=(\REnd_B^{\mathbb Z}(B))_0=kQ$, giving \Cref{eq:acyclic-derived-Morita}. The induced dg functor sends $kQ$ to $qB$ and identifies their endomorphism dg algebras. It also takes the silting generator $kQ$ to $qB$, proving \textup{(2)}. Finally, $\per(kQ)=\cd^b(\operatorname{mod}kQ)$ because $kQ$ is finite-dimensional and hereditary.

For \textup{(4)}, \Cref{lem:graded-noetherian-functorial-finiteness}\textup{(2)} gives the functorial finiteness of $\cp_{\mathrm{diag}}$ in $\cm_{\mathrm{diag}}$. The function $f(p)=p$ satisfies both \Cref{Ass}\textup{(A2)(a)} and \Cref{eq:reduced-dualizing-condition}. Together with the properties of $A$ established above, this allows us to apply the functorial-finiteness form of \Cref{Thm:CTC-in-gradedHiggs}. Since $f(p)-f(p-1)=1$, we have $r=1$ and $d+r=2$. That theorem gives the asserted Frobenius structure, projective--injectives and cluster-tilting subcategory, together with the first equivalence in \Cref{eq:stable-diagonal-Higgs-derived-Q}; the second follows from \textup{(3)}. The stable Serre functor corresponds to $(1)[1]$ on $\cc^{\mathbb Z}(B)$. Comparing this with $\tau[1]$ on $\cd^b(\operatorname{mod}kQ)$, it implies that the Adams shift $(1)$ on $\cc^{\mathbb Z}(B)$ corresponds to the Auslander--Reiten translation $\tau$ of $\cd^b(\operatorname{mod}kQ)$.
\end{proof}

\begin{Thm}\label{Thm:Dynkin-constant-grading}
Let $Q$ be a Dynkin quiver, with $A=\Gamma_Q^{\mathrm{fr}}$, $B=\Pi_2(kQ)$, and $e$ as above. Take $f(p)=0$, put $\ch_0^{\mathbb Z}(Q):=\ch_f^{\mathbb Z}(\Gamma_Q^{\mathrm{fr}},e)$, and set
\[\cp_0=\add\{e\Gamma_Q^{\mathrm{fr}}(-p)\mid p\in\mathbb Z\},\qquad \cm_0=\add\{\Gamma_Q^{\mathrm{fr}}(-p)\mid p\in\mathbb Z\}.\]
Then the following hold.
\begin{enumerate}
\item The $t$-structure associated with $f$ restricts to $\per^{\mathbb Z}(B)$, with
\[\overline{\tau}_f^{\ge0}\bigl(\per^{\mathbb Z}(B)\bigr)\subseteq\pvd^{\mathbb Z}(B).\]

\item The category $\ch_0^{\mathbb Z}(Q)$ is a Frobenius extriangulated category with projective--injective subcategory $\cp_0$.

\item There are triangle equivalences
\begin{equation}\label{eq:stable-constant-Higgs-Dynkin}\ch_0^{\mathbb Z}(Q)/[\cp_0]\simeq\cc^{\mathbb Z}(B)\simeq \cd^b(\operatorname{mod}kQ).\end{equation}
The Serre functor of the stable category is $(1)[1]$.

\item The subcategory $\cm_0$ is a $1$-cluster-tilting subcategory of $\ch_0^{\mathbb Z}(Q)$. In particular,
\[\ch_0^{\mathbb Z}(Q)=\cm_0.\]
\end{enumerate}
\end{Thm}

\begin{proof}
Let $\overline Q$ be the doubled quiver of $Q$. The presentation above gives
\[H^0(B)\cong\Pi(Q):=k\overline Q/\langle\rho\rangle,\qquad \rho=\sum_{\beta\in Q_1}[\beta,\beta^*];\]
see also \cite[Section~5]{KellerWang2023}. Since $Q$ is Dynkin, $\Pi(Q)$ is finite dimensional. As $B$ is graded homologically smooth, the finite-dimensional alternative of \Cref{Thm:restriction-to-perfect} applies to $f=0$ and proves \textup{(1)}.

We next verify the functorial finiteness of $\cp_0$ in $\cm_0$. Put $R=H^0(A)$ and $R_e=eRe$. Since $\Pi(Q)$ is finite dimensional and graded by path length, choose $N\ge1$ such that every path in $\overline Q$ of length at least $N$ lies in $\langle\rho\rangle$.

Cutting at the last frozen vertex shows that $Re$ is generated over $R_e$ by $e$ and the paths $v\eta_i$, where $v$ is a path in $\overline Q$ starting at $i$. If $v$ has length at least $N$, it is a linear combination of terms $u(e_j\rho e_j)w$ with the same endpoints and length as $v$. The relations in $R$ give
\[u(e_j\rho e_j)w\eta_i=-(u\eta_j)(\eta_j^*w\eta_i),\qquad \eta_j^*w\eta_i\in R_e.\]
Since $u$ is strictly shorter than $v$, induction on path length yields
\[Re=R_e+\sum_{\substack{i\in Q_0,\ v\text{ a path in }\overline Q\\v\text{ starts at }i,\ \operatorname{length}(v)<N}}(v\eta_i)R_e.\]

These generators are finite in number and homogeneous. Applying the same argument to the opposite algebra shows that $eR$ is finitely generated as a graded left $R_e$-module. \Cref{lem:graded-noetherian-functorial-finiteness}\textup{(1)} therefore gives the functorial finiteness of $\cp_0$ in $\cm_0$.

The finite-dimensionality of $H^0(B)$ verifies \Cref{Ass}\textup{(A2)(b)} and \Cref{eq:reduced-dualizing-condition}. Together with the properties of $A$ established above, the functorial-finiteness form of \Cref{Thm:CTC-in-gradedHiggs} therefore applies with $r=0$ and $d+r=1$. It gives \textup{(2)}, the first equivalence in \Cref{eq:stable-constant-Higgs-Dynkin}, the stable Serre functor $(1)[1]$, and the $1$-cluster-tilting property in \textup{(4)}. Since the defining orthogonality range for $1$-cluster-tilting is empty, $\ch_0^{\mathbb Z}(Q)=\cm_0$.

The second equivalence in \Cref{eq:stable-constant-Higgs-Dynkin} follows from \Cref{eq:acyclic-derived-Morita}, which is independent of $f$, and $\per(kQ)=\cd^b(\operatorname{mod}kQ)$.
\end{proof}

The Dynkin assumption is necessary for the properness assertion with $f=0$: if $Q$ is non-Dynkin, then $\Pi(Q)$ is infinite dimensional and
\[\overline{\tau}_f^{\ge0}B\simeq H^0(B)\cong\Pi(Q)\notin\pvd^{\mathbb Z}(B).\]
For $f(p)=p$, this properness assertion holds for every finite acyclic $Q$ by \Cref{Ex:acyclic-diagonal-grading}\textup{(1)}.

\begin{Rem}\label{Rem:Keller-Scherotzke-comparison}
Let $\cR$ and $\cs$ be the regular and singular Nakajima categories of \cite{KellerScherotzke2016}. Thus $\cR$ is the mesh category of the framed repetition quiver $\mathbb Z\widetilde Q$, with relations imposed only at the non-frozen vertices, and $\cs$ is its full subcategory on the frozen vertices.

Consider the $k$-category
\begin{equation}\label{eq:unfolded-relative-completion}\cR_{\Gamma}:=H^0\!\left(\add\{e_v\Gamma_Q^{\mathrm{fr}}(-p)\mid v\in Q_0^{\mathrm{fr}},\ p\in\mathbb Z\}\right).\end{equation}
Here the additive closure is taken in the dg category of bigraded right $\Gamma_Q^{\mathrm{fr}}$-modules, with morphism complexes of Adams degree zero, and $H^0$ takes their zeroth cohomology. Unfolding the Adams grading gives the repetition quiver, while \Cref{eq:relative-preprojective-differential} gives its mesh relations. After reindexing the frozen vertices and rescaling arrows, the generators $e_i\Gamma_Q^{\mathrm{fr}}(-p)$ and $e_{i'}\Gamma_Q^{\mathrm{fr}}(-p)$ correspond to the vertices $(i,p)$ and $(i',p-1)$, respectively. Taking finite direct sums and direct summands therefore gives
\begin{equation}\label{eq:unfolding-is-regular-Nakajima}\cm_0\simeq\cR_{\Gamma}\simeq\operatorname{proj}\cR,\end{equation}
where $\operatorname{proj}\cR$ denotes the finitely generated projective right $\cR$-modules, and $e_v\Gamma_Q^{\mathrm{fr}}(-p)$ corresponds to the representable projective of $\cR$ at the vertex $(v,p)$.

Now assume that $Q$ is Dynkin. \Cref{Thm:Dynkin-constant-grading} and restriction along $\cs\hookrightarrow\cR$ give equivalences of additive categories
\begin{equation}\label{eq:Higgs-Nakajima-Gorenstein}\ch_0^{\mathbb Z}(Q)=\cm_0\simeq\operatorname{proj}\cR\xrightarrow{\sim}\operatorname{gpr}\cs,\end{equation}
where the last equivalence is \cite[Theorem~2.17]{KellerScherotzke2016}. Here $\operatorname{gpr}\cs$ consists of finitely presented Gorenstein-projective right $\cs$-modules; its projective--injective objects are $\operatorname{proj}\cs$, corresponding to $\cp_0$.

Write $\operatorname{mod}\cs$ for the category of finitely presented right $\cs$-modules and put
\[\cd_{\mathrm{sg}}(\cs):=\cd^b(\operatorname{mod}\cs)/K^b(\operatorname{proj}\cs).\]
The categories $\cR$ and $\cs$ are coherent, and finitely presented $\cR$-modules have projective resolutions of length at most two; this follows from \cite[Section~3.1, Lemma~3.5 and Proposition~5.10]{KellerScherotzke2016}. Lifting a finite presentation from $\cs$ to $\cR$ and restricting such a resolution gives a finite Gorenstein-projective resolution over $\cs$, by \cite[Theorem~2.17]{KellerScherotzke2016}. Hence the Buchweitz equivalence identifies $\underline{\operatorname{gpr}}\,\cs$ with $D_{\mathrm{sg}}(\cs)$.

Let $\Xi\colon\underline{\operatorname{gpr}}\,\cs\xrightarrow{\sim}\cd^b(\operatorname{mod}kQ)$ be the triangle equivalence of \cite[Theorem~2.16]{KellerScherotzke2016}. Composing $\Xi^{-1}$ with \Cref{eq:stable-constant-Higgs-Dynkin}, and using the preceding identification, gives triangle equivalences
\begin{equation}\label{eq:stable-Higgs-KS-singularity}\frac{\ch_0^{\mathbb Z}(Q)}{[\cp_0]}\simeq\underline{\operatorname{gpr}}\,\cs\simeq D_{\mathrm{sg}}(\cs)\simeq \cd^b(\operatorname{mod}kQ).\end{equation}

Let $H\colon k(\mathbb ZQ)\to \cd^b(\operatorname{mod}kQ)$ be Happel's embedding, normalized by $H(i,0)=e_i kQ$, and put $\sigma(i,p)=(i',p-1)$. Keller--Scherotzke's stratifying functor $\Phi$ takes finite-dimensional right $\cs$-modules to $\cd^b(\operatorname{mod}kQ)$ and satisfies
\begin{equation}\label{eq:KS-stratifying-functor-in-Higgs-model}\Phi(M)\simeq\Xi(\Omega_{\cs}M),\qquad \Phi(S_{\sigma(x)})\simeq H(x),\end{equation}
by \cite[Theorems~2.7 and~2.16]{KellerScherotzke2016}. Here $S_{\sigma(x)}$ is the simple module at $\sigma(x)$, and $\Omega_{\cs}M$ is a first syzygy, viewed in $\underline{\operatorname{gpr}}\,\cs$. In Dynkin type, the arrow rescalings can be chosen independently of $p$, so the Adams shift $(1)$ corresponds under unfolding to $\tau(i,p)=(i,p-1)$, consistently with \Cref{Ex:acyclic-diagonal-grading}.

Over $k=\mathbb C$, the finite-dimensional $\cs$-modules of a fixed dimension vector form an affine graded quiver variety. By \cite[Theorem~2.7]{KellerScherotzke2016}, two such modules lie in the same stratum if and only if their images under $\Phi$ are isomorphic. In Dynkin type, these images can therefore be viewed as objects of the stable Higgs category through \Cref{eq:stable-Higgs-KS-singularity}.

For connected non-Dynkin $Q$, the category $\cs$ is hereditary, so $D_{\mathrm{sg}}(\cs)=0$, while $\Phi$ takes values in the proper additive subcategory $\add(\operatorname{Im}H)\subset \cd^b(\operatorname{mod}kQ)$; see \cite[Theorem~2.7 and Section~5.1]{KellerScherotzke2016}. Thus $\operatorname{gpr}\cs$ has zero stable category. By contrast, the diagonal Higgs category provides a Frobenius extriangulated category for the whole derived category:
\[\ch_{\mathrm{diag}}^{\mathbb Z}(Q)/[\cp_{\mathrm{diag}}]\simeq \cd^b(\operatorname{mod}kQ),\]
as established in \Cref{eq:stable-diagonal-Higgs-derived-Q}.
\end{Rem}

\subsection{Graded maximal Cohen--Macaulay modules over isolated singularities}
\label{subsec:higgs-graded-mcm}

When $f=0$, the frozen objects are precisely the Adams shifts of $eA$.
Below we give a criterion under which restriction along $e$ identifies
$\ch^{\mathbb Z}_0(A,e)$ with the category of graded maximal
Cohen--Macaulay modules over $eAe$, and apply it to isolated quotient
singularities to determine the cluster-tilting subcategory and the silting
object from \Cref{Thm:CTC-in-gradedHiggs}.

Let $A$ be an ordinary graded $k$-algebra concentrated in nonpositive
Adams degrees, viewed as a dg algebra in cohomological degree zero; let
$e\in A_0$ be an idempotent and set $R=eAe$. Denote by
$\operatorname{gr}R$ the category of finitely generated graded right
$R$-modules with degree-zero morphisms. When $R$ is commutative,
connected, and Gorenstein, we denote by
$\operatorname{CM}^{\mathbb Z}(R)\subseteq\operatorname{gr}R$ the exact
subcategory of graded maximal Cohen--Macaulay modules.

\begin{Prop}\label{prop:higgs-mcm-corner}
	Assume that $A$ is graded right Noetherian of finite graded right
	global dimension, that $Ae$ is finitely generated as a right
	$R$-module, and that $A/AeA$ is finite-dimensional over $k$.
	Assume further that $R$ is a commutative Noetherian connected graded
	Gorenstein $k$-algebra of finite Krull dimension. Then restriction
	$F(X)=Xe$ induces a triangle equivalence
	\[ \cc^{\mathbb Z}(A,e)\xrightarrow{\ \sim\ }\cd^b(\operatorname{gr}R), \]
	which restricts to an equivalence of Frobenius exact categories
	\[ \ch^{\mathbb Z}_0(A,e)\xrightarrow{\ \sim\ } \operatorname{CM}^{\mathbb Z}(R),\qquad eA(t)\longmapsto R(t). \]
	In particular, we obtain triangle equivalences
	\[ \underline{\ch^{\mathbb Z}_0(A,e)} \simeq \underline{\operatorname{CM}}^{\mathbb Z}(R) \simeq D_{\mathrm{sg}}^{\mathbb Z}(R) :=\cd^b(\operatorname{gr}R)/\per^{\mathbb Z}(R). \]
\end{Prop}

\begin{proof}
	Since $A$ is graded right Noetherian of finite graded right global
	dimension, $\per^{\mathbb Z}(A)$ identifies with
	$\cd^b(\operatorname{gr}A)$. Because $Ae$ is finitely generated over
	$R$, the exact functor $F$ preserves finite generation and induces
	a triangle functor $\cd^b(\operatorname{gr}A)\to \cd^b(\operatorname{gr}R)$.
	A complex $X\in \cd^b(\operatorname{gr}A)$ lies in $\ker F$ if and
	only if $H^i(X)e=0$ for all $i$; each $H^i(X)$ is then a
	finitely generated module over $A/AeA$, hence finite-dimensional over
	$k$. Thus $\ker F=\pvd_e^{\mathbb Z}(A)$.
	
	To construct a quasi-inverse on the Verdier quotient, apply
	$L(N)=N\otimes_R eA$ termwise to bounded complexes of finitely
	generated graded $R$-modules. Since $FL\simeq\mathrm{id}$ on
	complexes and $F$ is exact, $L$ maps the cone of any
	quasi-isomorphism into $\ker F$, and therefore descends to a triangle
	functor $\cd^b(\operatorname{gr}R)\to \cd^b(\operatorname{gr}A)/\pvd_e^{\mathbb Z}(A)$.
	For any complex $X$, the multiplication map $(Xe)\otimes_R eA\to X$
	becomes an isomorphism after applying $F$, so its cone also lies in
	$\ker F$; this yields the natural isomorphism $LF\simeq\mathrm{id}$
	in the quotient. As $\cd^b(\operatorname{gr}R)$ is idempotent complete,
	the resulting equivalence extends to $\cc^{\mathbb Z}(A,e)$.
	
	Under this equivalence, $\cp_0$ corresponds to
	$\add\{R(t)\mid t\in\mathbb Z\}$. For $Y\in \cd^b(\operatorname{gr}R)$,
	put $Y^\vee=\mathbf R\underline{\operatorname{Hom}}_R(Y,R)$, where
	$\underline{\operatorname{Hom}}$ denotes the graded Hom. For all
	$j,t\in\mathbb Z$, we have
	\[ \Hom(R(-t),Y[j])=H^j(Y)_t,\qquad \Hom(Y,R(t)[j])=H^j(Y^\vee)_t. \]
	The two orthogonality conditions in \Cref{Def:grcc-ghiggs}
	thus say precisely that both $Y$ and $Y^\vee$ lie in $\cd^{\leq 0}$.
	Since the dual of a complex in $D^{\leq 0}$ lies in $\cd^{\geq 0}$,
	graded Gorenstein biduality gives $Y\simeq (Y^\vee)^\vee\in D^{\geq 0}$,
	forcing $Y$ to be concentrated in degree zero as a module
	$M\in\operatorname{gr}R$. For such a module, the condition
	$M^\vee\in D^{\leq 0}$ amounts to
	$\underline{\operatorname{Ext}}_R^j(M,R)=0$ for all $j>0$, which
	characterizes $\operatorname{CM}^{\mathbb Z}(R)$ inside
	$\operatorname{gr}R$.
	
	Triangles in $\cd^b(\operatorname{gr}R)$ with all three terms in
	$\operatorname{CM}^{\mathbb Z}(R)$ correspond to graded short exact
	sequences, and the projective--injective objects are precisely the
	finitely generated graded projectives $\add\{R(t)\mid t\in\mathbb Z\}$.
	Passing to stable categories and applying the graded Buchweitz
	equivalence \cite{BuchweitzMCM} yields the final assertion.
\end{proof}

\bigskip

If $R$ has an isolated singularity with graded canonical module
$\omega_R\simeq R(\gamma)$, then graded Auslander--Reiten duality
\cite[Corollary~2.5]{IyamaTakahashiQuotient} identifies the Serre functor
on $\underline{\ch^{\mathbb Z}_0(A,e)}$ with $(\gamma)[\dim R-1]$.

In particular, suppose $R$ is regular away from its homogeneous maximal
ideal $\mathfrak m$, and let $T=R\oplus M$ with
$M\in\operatorname{CM}^{\mathbb Z}(R)$. Set
$A=\underline{\operatorname{End}}_R(T) =\bigoplus_{t\in\mathbb Z}\Hom_{\operatorname{gr}R}(T,T(t))$,
and let $e\in A_0$ be the projection onto $R$, so that $eAe=R$ and
$Ae\simeq T$. Since $T$ is finitely generated over $R$, the algebra
$A$ is module-finite over $R$, hence graded right Noetherian. The
two-sided ideal $AeA$ consists of sums of homogeneous endomorphisms
factoring through shifts of $R$. Because $T$ is locally free outside
$\mathfrak m$, the finite $R$-module $A/AeA$ is supported at
$\mathfrak m$ and therefore finite-dimensional over $k$. Thus,
whenever $A$ is concentrated in nonpositive Adams degrees and has finite
graded right global dimension, \Cref{prop:higgs-mcm-corner}
applies with $F(A)=T$.

\begin{Cor}\label{cor:higgs-mcm-mckay}
	Let $k$ be an algebraically closed field of characteristic zero, and
	let $G\subseteq\mathrm{SL}_n(k)$ ($n\geq 2$) be a finite subgroup
	acting freely on $k^n\setminus\{0\}$. Set
	\[ S=k[x_1,\ldots,x_n],\qquad R=S^G,\qquad A=S\#G,\qquad e=|G|^{-1}\sum_{g\in G}g,\qquad B=A/^{\mathbf L}AeA, \]
	where each $x_i$ has Adams degree $-1$ and each $g\in G$ has
	degree $0$. Then $F$ identifies $\ch^{\mathbb Z}_0(A,e)$ with
	$\operatorname{CM}^{\mathbb Z}(R)$ and, via
	\Cref{Thm:CTC-in-gradedHiggs}, induces a triangle equivalence
	\[ \cc^{\mathbb Z}(B)\simeq\underline{\operatorname{CM}}^{\mathbb Z}(R) \]
	sending the image of $B$ to the class $\underline{S}$ of $S$ in
	the stable category. Under this identification, the canonical
	$(n-1)$-cluster-tilting subcategory of $\ch^{\mathbb Z}_0(A,e)$
	corresponds to $\add\{S(t)\mid t\in\mathbb Z\}$, while the Serre
	functor on $\underline{\operatorname{CM}}^{\mathbb Z}(R)$ and the
	silting object of \Cref{Thm:finite-Adams-window} are given
	respectively by
	\[ (n)[n-1]\qquad\text{and}\qquad \bigoplus_{\ell=0}^{n-1}\underline{S(-\ell)}. \]
\end{Cor}

\begin{proof}
	Since $G\subseteq\mathrm{SL}_n(k)$ acts freely on
	$k^n\setminus\{0\}$, the invariant ring $R=S^G$ is Gorenstein and
	regular away from its homogeneous maximal ideal. By Auslander's graded
	isomorphism \cite[Theorem~3.2]{IyamaTakahashiQuotient}, we have
	$A=S\#G\simeq\underline{\operatorname{End}}_R(S)$, under which $e$
	corresponds to the Reynolds projection onto $R$. Since
	$S\in\operatorname{CM}^{\mathbb Z}(R)$ decomposes as
	$R\oplus(1-e)S$, the preceding discussion shows that $A$ is graded
	right Noetherian and $\dim_k(A/AeA)<\infty$. Because $kG$ is
	separable over $k$, the $G$-equivariant Koszul complex gives a
	finite resolution of $A$ by finitely generated graded projective
	$A^e$-modules of length $n$; in particular, $A$ has graded global
	dimension $n$, so \Cref{prop:higgs-mcm-corner} applies.
	
	The top exterior power of $\operatorname{span}_k\{x_1,\ldots,x_n\}$
	is the trivial $G$-representation in Adams degree $-n$, since
	$G\subseteq\mathrm{SL}_n(k)$. Dualizing the Koszul bimodule
	resolution therefore yields
	$\mathbf R\underline{\operatorname{Hom}}_{A^e}(A,A^e) \simeq A(-n)[-n]$, so the relative Serre functor is $(n)[n]$.
	With $A_0=kG$ and $\dim_k H^0(B)=\dim_k(A/AeA)<\infty$,
	\Cref{Thm:CTC-in-gradedHiggs} applies for $d=n-1$, $a=n$,
	and $r=0$. Its stable equivalence identifies the images of $A$ and
	$B$; combined with $F(A)=Ae\simeq S$, this yields the stated
	cluster-tilting subcategory and Serre functor.
	
	In the dg quotient $B=A/^{\mathbf L}AeA$, the contracting morphism
	for $e$ has Adams degree zero. Since $A$ is concentrated in
	nonpositive Adams degrees, the degree-zero part $B_0$ is the derived
	quotient of $kG$ by $e$. As $e$ is a central idempotent in
	$kG$, the algebra $kG$ splits as $kGe\times kG(1-e)$, whence
	$B_0\simeq kG(1-e)\simeq kG/(kG)e(kG)$ is proper. By
	\Cref{Thm:finite-Adams-window}, the silting object is the image
	of $\bigoplus_{\ell=0}^{n-1}B(-\ell)$, which corresponds to
	$\bigoplus_{\ell=0}^{n-1}\underline{S(-\ell)}$.
\end{proof}

\bigskip
Reversing the Adams grading gives the standard polynomial grading
$\deg x_i=1$ and turns $M(t)$ into $M(-t)$. In this grading,
the stable Serre functor is $(-n)[n-1]$ and the silting object is
$\bigoplus_{\ell=0}^{n-1}\underline{S(\ell)}$; this recovers the
results of Iyama--Takahashi
\cite[Theorems~1.5, 1.6 and~4.2]{IyamaTakahashiQuotient}. For an ungraded
Higgs realization via a higher McKay quiver with potential, see
\cite[Theorem~4.10(b)]{LiuHigherMcKay}.
\begin{Ex}
	Let $G=\mu_3$ act on $S=k[x,y,z]$ by scalar multiplication.
	Writing $S_j^{\mathrm{pol}}$ for the space of homogeneous polynomials
	of ordinary degree $j$, we have $S=R\oplus M_1\oplus M_2$, where
	$R=M_0$ is the third Veronese subring and
	\[ M_i=\bigoplus_{q\geq 0}S_{3q+i}^{\mathrm{pol}}\qquad (0\leq i\leq 2), \]
	with $S_{3q+i}^{\mathrm{pol}}$ placed in Adams degree $-(3q+i)$.
	Since the frozen summands are the shifts of $R$, the canonical
	$2$-cluster-tilting subcategory of
	$\underline{\operatorname{CM}}^{\mathbb Z}(R)$ is
	$\add\{M_1(t),M_2(t)\mid t\in\mathbb Z\}$. In the McKay quiver of
	$G$, the variables $x,y,z$ give three arrows from vertex $i$ to
	$i+1\pmod{3}$ for each $i$; deleting the trivial vertex leaves two
	vertices joined by three parallel arrows and no paths of length two, so
	$H^0(B)=A/AeA$ is the $3$-arrow Kronecker algebra. In our Adams
	grading, the stable Serre functor is $(3)[2]$.
\end{Ex}

\subsection{Relative singularity categories and Koszul duality}
We return to the general setup of \Cref*{subsec:relative-t-structures} and assume \textup{(A0)} and \textup{(A1)} of \Cref{Ass}, together with $\dim_k H^0(A)_0<\infty$. The remaining hypotheses of \Cref{Thm:CTC-in-gradedHiggs} are supposed only when that theorem is used. 

Put
\[S:=H^0(A)_0/\operatorname{rad}H^0(A)_0,\]
viewed as a dbg right $A$-module concentrated in bidegree $(0,0)$ via $A\to H^0(A)_0$. Write $S=S_{\mathrm u}\oplus S_{\mathrm f}$, where ${S_{\mathrm u}}$ is the sum of the simple direct summands annihilated by $e$, and ${S_{\mathrm f}}$ is the sum of the remaining simple direct summands. Choose cofibrant resolutions $P_{\mathrm u}\to S_{\mathrm u}$ and $P_{\mathrm f}\to S_{\mathrm f}$, put $P=P_{\mathrm u}\oplus P_{\mathrm f}$, and set
\begin{equation}\label{eq:dual-algebra}E:=\End_A^{\mathbb Z}(P)\simeq\REnd_A^{\mathbb Z}(S),\qquad u:=\mathrm{id}_{P_{\mathrm u}}\oplus0\in E^0_0.\end{equation}
Thus $u$ is a closed idempotent.

The truncation and composition-series argument in the proof of \Cref{Thm:restriction-to-perfect} gives
\begin{equation}\label{eq:finite-generators}\pvd^{\mathbb Z}(A)=\thick^{\mathbb Z}(S),\qquad \pvd_e^{\mathbb Z}(A)=\thick^{\mathbb Z}(S_{\mathrm u}).\end{equation}

By \textup{(A0)}, $S$ is perfect. The identities $\RHom_A^{\mathbb Z}(A(i),S)\simeq S(-i)$ show, by thick closure, that $\RHom_A^{\mathbb Z}(X,S)$ has finite total cohomology for every perfect $X$. Taking $X=S$ proves that $E$ is proper. Hence $\per^{\mathbb Z}(E)\subseteq\pvd^{\mathbb Z}(E)$.

\begin{Def}\rm\label{def:relative-singularity}
For a proper dbg algebra $E$ and a closed idempotent $u\in E^0_0$, define
\begin{equation}\label{eq:relative-singularity}\operatorname{sg}^{\mathbb Z}(E;u):=\Kar\!\left(\pvd^{\mathbb Z}(E)/\thick^{\mathbb Z}(uE)\right).\end{equation}
We call this the \emph{graded singularity category} of $E$ relative to the projective summand $uE$, and write $\operatorname{sg}^{\mathbb Z}(E):=\operatorname{sg}^{\mathbb Z}(E;1)$.
\end{Def}

\begin{Lem}\label{lem:koszul}
The functor
\begin{equation}\label{eq:covariant-koszul}D\RHom_A^{\mathbb Z}(-,S)\colon\per^{\mathbb Z}(A)\xrightarrow{\sim}\pvd^{\mathbb Z}(E)\end{equation}
is a triangle equivalence commuting with both $[1]$ and $(1)$. Moreover,
\begin{equation}\label{eq:morita-finite}\RHom_A^{\mathbb Z}(S,-)\colon\pvd^{\mathbb Z}(A)\xrightarrow{\sim}\per^{\mathbb Z}(E)\end{equation}
is a triangle equivalence satisfying $\RHom_A^{\mathbb Z}(S,S_{\mathrm u})\simeq uE$.
\end{Lem}

\begin{proof}
For the first assertion, retain the Adams grading and regrade by
\[\widetilde A_p^n=A_p^{n-2p},\qquad \widetilde M_p^n=M_p^{n-2p}.\]
This invertible regrading preserves Koszul signs, since $2p$ is even. It preserves perfect and perfectly valued objects and satisfies
\[\widetilde{M(i)[j]}\cong\widetilde M(i)[j+2i].\]
The third-quadrant condition gives $H^{>0}(\widetilde A)=0$ and $H^0(\widetilde A)=H^0(A)_0$, and
\[H^n(\widetilde A)\cong\bigoplus_{\lceil n/2\rceil\le p\le0}H^{n-2p}(A)_p\qquad(n\le0).\]
The sum is finite, and each summand is finite-dimensional by \Cref{lem:automatic-bigraded-local-finiteness}. Hence $\dim_k H^n(\widetilde A)<\infty$ for every $n$. Applying \cite[Theorem~4.4(1)]{Fushimi} to $\widetilde A$ with its Adams grading and the same semisimple module $S$ gives
\[\RHom_{\widetilde A}^{\mathbb Z}(-,S)\colon\per^{\mathbb Z}(\widetilde A)^{\mathrm{op}}\xrightarrow{\sim}\pvd^{\mathbb Z}(\widetilde E^{\mathrm{op}}),\]
where $\widetilde E$ is the corresponding regrading of $E$. Regrading back and composing with graded $k$-duality gives \Cref{eq:covariant-koszul}. Both $\RHom_A^{\mathbb Z}(-,S)$ and $D$ reverse cohomological and Adams shifts, so their composite commutes with $[1]$ and $(1)$.

For \Cref{eq:morita-finite}, derived Morita theory identifies $\thick^{\mathbb Z}(S)$ with $\per^{\mathbb Z}(E)$ via $\RHom_A^{\mathbb Z}(S,-)$. Equation~\Cref{eq:finite-generators} identifies the source with $\pvd^{\mathbb Z}(A)$. Under this equivalence, the projection onto $S_{\mathrm u}$ corresponds to left multiplication by $u$ on $E$, whose image is $uE$. Thus $\RHom_A^{\mathbb Z}(S,S_{\mathrm u})\simeq uE$.
\end{proof}

\begin{Lem}\label{lem:relative-nakayama}
For every $X\in\pvd_e^{\mathbb Z}(A)$, there is an isomorphism in $\cd^{\mathbb Z}(E)$
\begin{equation}\label{eq:relative-nakayama}D\RHom_A^{\mathbb Z}(X,S)\simeq\RHom_A^{\mathbb Z}(S,\mathbb S X).\end{equation}
Consequently, the equivalence in \Cref{eq:covariant-koszul} restricts to
\begin{equation}\label{eq:denominator}D\RHom_A^{\mathbb Z}(-,S)\colon\pvd_e^{\mathbb Z}(A)\xrightarrow{\sim}\thick^{\mathbb Z}(uE).\end{equation}
In particular, as derived right $E$-modules,
\begin{equation}\label{eq:relative-projective-injective}D(Eu)\simeq uE(a)[d+1].\end{equation}
\end{Lem}

\begin{proof}
For $X\in\per^{\mathbb Z}(A)$, put $\nu_A X:=D\RHom_A^{\mathbb Z}(X,A)$. Since $X$ is perfect, derived evaluation identifies $\RHom_A^{\mathbb Z}(X,Y)$ with $Y\lten_A\RHom_A^{\mathbb Z}(X,A)$. Tensor--Hom adjunction therefore gives, naturally in $Y\in\cd^{\mathbb Z}(A)$,
\begin{equation}\label{eq:derived-nakayama}D\RHom_A^{\mathbb Z}(X,Y)\simeq\RHom_A^{\mathbb Z}(Y,\nu_A X).\end{equation}
These identifications are natural on morphism complexes, so for $Y=S$ they respect the right $E$-actions defined using the chosen resolution $P$.

Let $X\in\pvd_e^{\mathbb Z}(A)$. \Cref{Ass}\textup{(A1)} and \Cref{eq:derived-nakayama} give, naturally in $Y\in\per^{\mathbb Z}(A)$,
\[\Hom_{\cd^{\mathbb Z}(A)}(Y,\mathbb S X)\cong D\Hom_{\cd^{\mathbb Z}(A)}(X,Y)\cong\Hom_{\cd^{\mathbb Z}(A)}(Y,\nu_A X).\]
Since $\mathbb S X$ is perfect, the Yoneda lemma yields a morphism $\theta_X\colon\mathbb S X\to\nu_A X$ such that $\Hom(Y,\theta_X)$ is an isomorphism for every $Y\in\per^{\mathbb Z}(A)$. In particular, taking $Y=A(i)[j]$ for all $i,j\in\mathbb Z$ shows that $\theta_X$ induces isomorphisms on all bigraded cohomology groups. Hence $\theta_X$ is an isomorphism in $\cd^{\mathbb Z}(A)$. Applying $\RHom_A^{\mathbb Z}(S,-)$ and using \Cref{eq:derived-nakayama} with $Y=S$ proves \Cref{eq:relative-nakayama} in $\cd^{\mathbb Z}(E)$.

By \Cref{Ass}\textup{(A1)}, $\mathbb S$ restricts to an autoequivalence of $\pvd_e^{\mathbb Z}(A)$. Equations~\Cref{eq:finite-generators}, \Cref{eq:morita-finite} and \Cref{eq:relative-nakayama} therefore give \Cref{eq:denominator}.

Finally, $\alpha^{-1}$ induces an isomorphism $S_\alpha\simeq S$ which restricts to $(S_{\mathrm u})_\alpha\simeq S_{\mathrm u}$ because $\alpha(e)=e$. Hence $\mathbb S S_{\mathrm u}\simeq S_{\mathrm u}(a)[d+1]$. Taking $X=S_{\mathrm u}$ in \Cref{eq:relative-nakayama} and using $\RHom_A^{\mathbb Z}(S_{\mathrm u},S)\simeq Eu$ and $\RHom_A^{\mathbb Z}(S,S_{\mathrm u})\simeq uE$ gives \Cref{eq:relative-projective-injective}.
\end{proof}

Combining \Cref{lem:koszul,lem:relative-nakayama} and passing to the idempotent completions of the Verdier quotients, we obtain the following theorem.
\begin{Thm}[Relative Koszul duality]\label{thm:relative-koszul}
The equivalence $D\RHom_A^{\mathbb Z}(-,S)$ induces a triangle equivalence
\begin{equation}\label{eq:main-equivalence}\cc^{\mathbb Z}(A,e)\xrightarrow{\sim}\operatorname{sg}^{\mathbb Z}(E;u).\end{equation}
These equivalences fit into a commutative diagram
\begin{equation}\label{eq:main-diagram} \begin{tikzcd}[column sep=large,row sep=large] \pvd_e^{\mathbb Z}(A) \arrow[r,hook] \arrow[d,"\sim"'] & \per^{\mathbb Z}(A) \arrow[r] \arrow[d,"{D\RHom_A^{\mathbb Z}(-,S)}"',"\sim"] & \cc^{\mathbb Z}(A,e) \arrow[d,"\sim"]\\
\thick^{\mathbb Z}(uE) \arrow[r,hook] & \pvd^{\mathbb Z}(E) \arrow[r] & \operatorname{sg}^{\mathbb Z}(E;u). \end{tikzcd} \end{equation}
The rows are Verdier quotient sequences up to idempotent completion of the last terms. All functors commute with the Adams and cohomological shifts.
\end{Thm}

For $e=0$ (hence $u=1$), \Cref{eq:main-equivalence} recovers $\cc^{\mathbb Z}(A)\simeq\operatorname{sg}^{\mathbb Z}(E)$, while for $e=1$ (hence $u=0$), it reduces to \Cref{eq:covariant-koszul}. In general, since $D\RHom_A^{\mathbb Z}(S,S)\simeq DE$, the equivalence in \Cref{eq:covariant-koszul} sends $\pvd^{\mathbb Z}(A)$ onto $\thick^{\mathbb Z}(DE)$ by \Cref{eq:finite-generators}.

\begin{Prop}[Compatibility with the frozen quotient]\label{prop:frozen-quotient}
Recall that $B=\overline A$ and put $E_{\mathrm u}=uEu$. Then
\begin{equation}\label{eq:corner-dual}E_{\mathrm u}\simeq\REnd_B^{\mathbb Z}(S_{\mathrm u}),\end{equation}
and $D\RHom_B^{\mathbb Z}(-,S_{\mathrm u})$ induces a triangle equivalence fitting into the commutative diagram
\begin{equation}\label{eq:reduction-diagram} \begin{tikzcd}[column sep=large,row sep=large] \cc^{\mathbb Z}(A,e) \arrow[r,"\sim"] \arrow[d,"\pi"'] & \operatorname{sg}^{\mathbb Z}(E;u) \arrow[d,"(-)u"]\\
\cc^{\mathbb Z}(B) \arrow[r,"\sim"] & \operatorname{sg}^{\mathbb Z}(E_{\mathrm u}). \end{tikzcd} \end{equation}
Let $Q_e$ be the image of $D\RHom_A^{\mathbb Z}(eA,S)$ in $\operatorname{sg}^{\mathbb Z}(E;u)$. Then
\begin{equation}\label{eq:second-quotient}\Kar\!\left(\operatorname{sg}^{\mathbb Z}(E;u)/\thick^{\mathbb Z}(Q_e)\right)\simeq\operatorname{sg}^{\mathbb Z}(E_{\mathrm u}).\end{equation}
\end{Prop}

\begin{proof}
By \Cref{prop:rec}, restriction along $A\to B$ identifies $\pvd^{\mathbb Z}(B)$ with $\pvd_e^{\mathbb Z}(A)$. Under this identification, the semisimple top of $H^0(B)_0=H^0(A)_0/H^0(A)_0eH^0(A)_0$ corresponds to $S_{\mathrm u}$. Thus $P_{\mathrm u}\otimes_A B$ is a cofibrant resolution of $S_{\mathrm u}$ over $B$, 
which induces a quasi-isomorphism
\[E_{\mathrm u}=\End_A^{\mathbb Z}(P_{\mathrm u})\xrightarrow{\sim}\End_B^{\mathbb Z}(P_{\mathrm u}\otimes_A B).\]
This proves \Cref{eq:corner-dual}.

The horizontal equivalences in \Cref{eq:reduction-diagram} are given by \Cref{thm:relative-koszul} for $(A,e)$ and $(B,0)$, using \Cref{eq:corner-dual} in the latter case. The functor $(-)u$ preserves perfectly valued objects and sends $uE$ to $E_{\mathrm u}$, hence induces the right vertical functor in \Cref{eq:reduction-diagram}. Moreover for $X\in\per^{\mathbb Z}(A)$, we have
\begin{equation}\label{eq:adjunction-square}\bigl(D\RHom_A^{\mathbb Z}(X,S)\bigr)u\simeq D\RHom_A^{\mathbb Z}(X,S_{\mathrm u})\simeq D\RHom_B^{\mathbb Z}(\pi X,S_{\mathrm u}).\end{equation}
Thus the commutativity of \Cref{eq:reduction-diagram} holds by passing to the quotients and their idempotent completions proves.

By \Cref{prop:rec} and compact localization, $\pi$ induces an equivalence
\[\Kar\!\left(\per^{\mathbb Z}(A)/\thick^{\mathbb Z}(eA)\right)\xrightarrow{\sim}\per^{\mathbb Z}(B).\]
Since $\pi$ identifies $\pvd_e^{\mathbb Z}(A)$ with $\pvd^{\mathbb Z}(B)$, the passage to Verdier localization yields an equivalence
\[\Kar\!\left(\cc^{\mathbb Z}(A,e)/\thick^{\mathbb Z}(eA)\right)\xrightarrow{\sim}\cc^{\mathbb Z}(B),\]
where $eA$ denotes its image in $\cc^{\mathbb Z}(A,e)$. Applying the relative Koszul duality equivalences in \Cref{eq:reduction-diagram} gives \Cref{eq:second-quotient}.
\end{proof}

Put $\cq_f=\add\{Q_e(-i)[f(i)]\mid i\in\mathbb Z\}$ and
\begin{equation}\label{eq:dual-higgs}\ch_{E,f}^{\mathbb Z}:=\{Y\in\operatorname{sg}^{\mathbb Z}(E;u)\mid\Hom(Y,\cq_f[>0])=0=\Hom(\cq_f,Y[>0])\}.\end{equation}
The equivalence in \Cref{eq:main-equivalence} sends $\cp_f$ onto $\cq_f$ and hence restricts to an extriangulated equivalence between $\ch_f^{\mathbb Z}(A,e)$ and $\ch_{E,f}^{\mathbb Z}$. By \Cref{Thm:CTC-in-gradedHiggs,Thm:finite-Adams-window} and the lower horizontal equivalence in \Cref{eq:reduction-diagram}, we obtain the following corollary.

\begin{Cor}\label{cor:higgs}\label{cor:window}
\begin{itemize}
\item[(1)] Under the hypotheses of \Cref{Thm:CTC-in-gradedHiggs}\textup{(1)}, there is an equivalence of Frobenius extriangulated categories
\[\ch_f^{\mathbb Z}(A,e)\xrightarrow{\sim}\ch_{E,f}^{\mathbb Z},\]
with projective--injective subcategory $\cq_f$ on the right, and triangle equivalences
\begin{equation}\label{eq:stable-higgs-singularity}\frac{\ch_f^{\mathbb Z}(A,e)}{[\cp_f]}\simeq\frac{\ch_{E,f}^{\mathbb Z}}{[\cq_f]}\simeq\operatorname{sg}^{\mathbb Z}(E_{\mathrm u}).\end{equation}
If $f(i)-f(i-a)=r$ for all $i\in\mathbb Z$ and $d+r\geq1$,  the image of $\cm_f$ under \Cref{eq:main-equivalence} is a $(d+r)$-cluster-tilting subcategory of $\ch_{E,f}^{\mathbb Z}$.

\item[(2)] Under the hypotheses of \Cref{Thm:finite-Adams-window}, let $T$ and $\Lambda$ be as above. The image of $D\RHom_B^{\mathbb Z}(T,S_{\mathrm u})$ in $\operatorname{sg}^{\mathbb Z}(E_{\mathrm u})$ is a silting generator, and there is a triangle equivalence
\[\operatorname{sg}^{\mathbb Z}(E_{\mathrm u})\simeq\per(\Lambda).\]
Under the hypotheses of both \textup{(1)} and \textup{(2)}, we have
\[\ch_f^{\mathbb Z}(A,e)/[\cp_f]\simeq\per(\Lambda).\]
\end{itemize}
\end{Cor}

\begin{Ex}\label{ex:a1}
Let $Q=A_1$, $A=\Gamma_Q^{\mathrm{fr}}$ and $e=e_{1'}$, as in \Cref*{subsec:stable-graded-higgs}. Write $\eta\colon1'\to1$ for the framing arrow, $\eta^*$ for its reverse and $t=t_1$. Thus
\[\operatorname{bideg}(\eta)=(0,0),\quad\operatorname{bideg}(\eta^*)=(-1,0),\quad\operatorname{bideg}(t)=(-1,-1),\qquad d(t)=\eta\eta^*.\]
The canonical morphism
\[A\longrightarrow R:=k\bigl(1'\underset{\eta^*}{\overset{\eta}{\rightleftarrows}}1\bigr)/(\eta\eta^*)\]
is a graded quasi-isomorphism. 

Put $P_i=e_iR$ for $i\in\{1',1\}$. The simple right modules have minimal graded projective resolutions
\begin{equation}\label{eq:a1-resolutions} \begin{aligned} 0&\longrightarrow P_1(1)\xrightarrow{\eta^*}P_{1'}\longrightarrow S_{1'}\longrightarrow0,\\
0&\longrightarrow P_1(1)\xrightarrow{\eta^*}P_{1'}\xrightarrow{\eta}P_1\longrightarrow S_1\longrightarrow0. \end{aligned} \end{equation}
Here the maps are given by left multiplication. Denote the corresponding projective complexes by $L_{1'}$ and $L_1$, concentrated in degrees $[-1,0]$ and $[-2,0]$, respectively.

Consider the closed maps of cohomological degree one
\[\xi\colon L_{1'}\to L_1(1),\qquad\zeta\colon L_1\to L_{1'},\]
where the only nonzero component of $\xi$ is in degree $-1$ with $\mathrm{id}_{P_1(1)}$ and 
those of $\zeta$ in degrees $-2$ and $-1$ are $-\mathrm{id}_{P_1(1)}$ and $\mathrm{id}_{P_{1'}}$, respectively. We have $\zeta\xi=0$, whereas $\xi\zeta$ represents a nonzero class in $\Ext_R^2(S_1,S_1(1))$.

Let $u$ and $v=1-u$ be the idempotents of $\End_R^{\mathbb Z}(L_{1'}\oplus L_1)$ corresponding to the summands $L_1$ and $L_{1'}$, respectively. 
Moreover, we have
\begin{equation}\label{eq:a1-dual}E\simeq k\bigl(v\underset{\zeta}{\overset{\xi}{\rightleftarrows}}u\bigr)/(\zeta\xi),\qquad\operatorname{bideg}(\xi)=(1,1),\quad\operatorname{bideg}(\zeta)=(0,1);\end{equation}
and
\[E_{\mathrm u}=uEu=k[\epsilon]/(\epsilon^2),\qquad\epsilon=\xi\zeta,\quad\operatorname{bideg}(\epsilon)=(1,2).\]
By \Cref{prop:frozen-quotient}, this is the Koszul dual of $B\simeq\Pi_2(k)=k[t]$, with $\operatorname{bideg}(t)=(-1,-1)$ and $d(t)=0$. We obtain
\[\cc^{\mathbb Z}(A,e)\simeq\operatorname{sg}^{\mathbb Z}(E;u),\qquad\operatorname{sg}^{\mathbb Z}(E_{\mathrm u})\simeq\cc^{\mathbb Z}(B)\simeq\cd^b(\operatorname{mod}k),\]
where the last equivalence is \Cref{eq:acyclic-derived-Morita}. In particular, $\operatorname{sg}^{\mathbb Z}(E_{\mathrm u})$ is triangle equivalent to the stable categories of both Higgs categories in \Cref{Ex:acyclic-diagonal-grading,Thm:Dynkin-constant-grading}.
\end{Ex}

\section{Weighted totalizations and orbit categories}\label{sec:weighted-totalization-singularity}

Fix integers $N\geq0$ and $c$, and put $f_N(p)=Np+c$. We retain the setup of \Cref{sec:graded-higgs}, including $\dim_k H^0(A)_0<\infty$, and \Cref{Ass} with $f=f_N$. Write $B=\overline A$. Condition~\Cref{eq:reduced-dualizing-condition} is automatic: for $N>0$, the function $f_N$ tends to $+\infty$ on the right; for $N=0$, \Cref{Ass}\textup{(A2)} requires $H^0(B)$ to be finite dimensional.

\subsection{Weighted totalization}
Let $C=\bigoplus_{p,q\in\mathbb Z}C_p^q$ be a dbg algebra.

\begin{Def}\rm\label{def:weighted-totalization}
The \emph{$N$-weighted totalization} of $C$ is the dg algebra $C_N^{\mathrm{tot}}=\operatorname{Tot}_N(C)$ with the original differential and
\begin{equation}\label{eq:weighted-total-space}(C_N^{\mathrm{tot}})^m=\bigoplus_{Np+q=m}C_p^q,\end{equation}
whose product is
\begin{equation}\label{eq:weighted-total-product}x\star y=(-1)^{Nps}xy\qquad(x\in C_p^q,\ y\in C_r^s).\end{equation}
For a bigraded right $C$-module $M$, its totalization has the action
\begin{equation}\label{eq:weighted-total-module-action}m\star x=(-1)^{Nps}mx\qquad(m\in M_p^q,\ x\in C_r^s).\end{equation}
We also use the degree-zero dg automorphism
\begin{equation}\label{eq:Adams-sign-automorphism}\varepsilon_C(x)=(-1)^p x\qquad(x\in C_p^q).\end{equation}
\end{Def}
Thus $\operatorname{Tot}_0$ forgets the Adams grading, and $\operatorname{Tot}_1$ is diagonal totalization.

\begin{Lem}\label{lem:weighted-totalization-basic-properties}
Totalization induces an exact functor $\operatorname{Tot}_N\colon\cd^{\mathbb Z}(C)\to\cd(C_N^{\mathrm{tot}})$ preserving direct sums and perfect objects. If $C$ is graded homologically smooth, then $C_N^{\mathrm{tot}}$ is homologically smooth. There are natural isomorphisms
\begin{equation}\label{eq:weighted-total-shifts}\operatorname{Tot}_N(M(r)[s])\simeq\operatorname{Tot}_N(M)[Nr+s].\end{equation}
In particular, for
\begin{equation}\label{eq:weighted-orbit-functor}F_N:=(-1)[N],\qquad\operatorname{Tot}_N F_N\simeq\operatorname{Tot}_N,\end{equation}
the induced morphism
\begin{equation}\label{eq:weighted-covering-Hom}\bigoplus_{\ell\in\mathbb Z}\RHom_C(X,F_N^\ell Y)\longrightarrow\RHom_{C_N^{\mathrm{tot}}}(\operatorname{Tot}_N X,\operatorname{Tot}_N Y)\end{equation}
is a quasi-isomorphism for $X\in\per^{\mathbb Z}(C)$ and $Y\in\cd^{\mathbb Z}(C)$.
\end{Lem}

\begin{proof}
The signs in \Cref{eq:weighted-total-product} and \Cref{eq:weighted-total-module-action} make $C_N^{\mathrm{tot}}$ a dg algebra and $\operatorname{Tot}_N M$ a dg module over it. Since
\begin{equation}\label{eq:weighted-total-cohomology}H^m(\operatorname{Tot}_N M)\cong\bigoplus_{Np+q=m}H^q(M)_p,\end{equation}
totalization preserves quasi-isomorphisms.

To prove \Cref{eq:weighted-total-shifts}, define the map by $m\mapsto(-1)^{Nrq}m$, where $m\in M_p^q$ is taken before shifting. This map is invertible and respects the differential and the right action. Totalization also commutes with cones and direct sums, so it induces an exact functor on derived categories. Since it sends $C(r)[s]$ to $C_N^{\mathrm{tot}}[Nr+s]$, it preserves perfect objects.

For the smoothness assertion, extend totalization to bimodules using the actions
\[a\star m=(-1)^{N\deg(a)|m|}am,\qquad m\star b=(-1)^{N\deg(m)|b|}mb.\]
The map
\[C_N^{\mathrm{tot}}\otimes_k C_N^{\mathrm{tot}}\xrightarrow{\sim}\operatorname{Tot}_N(C\otimes_k C),\qquad x\otimes y\longmapsto(-1)^{N\deg(x)|y|}x\otimes y,\]
is an isomorphism of outer free bimodules. Shifted free bimodules are also sent to perfect bimodules. Exactness and preservation of direct summands therefore show that totalization preserves perfect bimodules. Applied to the diagonal, this proves that $C_N^{\mathrm{tot}}$ is homologically smooth whenever $C$ is graded homologically smooth.

Finally, iterate the isomorphism $\operatorname{Tot}_N F_N\simeq\operatorname{Tot}_N$ and its inverse to define the map in \Cref{eq:weighted-covering-Hom}. For $X=C(r)[s]$, this map is a quasi-isomorphism by the decomposition of $Y$ into its Adams-degree components. The class of $X$ for which it is a quasi-isomorphism is thick, so the assertion follows for every perfect $X$.
\end{proof}

We use the same notation for the ungraded constructions, omitting the superscript $\mathbb Z$.

\begin{Prop}\label{prop:totalization-derived-quotient}
There is a canonical quasi-isomorphism
\begin{equation}\label{eq:totalization-derived-quotient}B_N^{\mathrm{tot}}\simeq A_N^{\mathrm{tot}}/^{\mathbf L}(A_N^{\mathrm{tot}}eA_N^{\mathrm{tot}}).\end{equation}
Totalization sends $\pvd_e^{\mathbb Z}(A)$ to $\pvd_e(A_N^{\mathrm{tot}})$ and $\pvd^{\mathbb Z}(B)$ to $\pvd(B_N^{\mathrm{tot}})$. It therefore induces functors on the corresponding cluster categories:
\[\operatorname{Tot}_N\colon\cc^{\mathbb Z}(A,e)\longrightarrow\cc(A_N^{\mathrm{tot}},e),\qquad\operatorname{Tot}_N\colon\cc^{\mathbb Z}(B)\longrightarrow\cc(B_N^{\mathrm{tot}}).\]
\end{Prop}

\begin{proof}
Use  $B=A\langle h\rangle/(he=eh=h)$ with $\operatorname{bideg}(h)=(0,-1)$ and $d(h)=e$. Under totalization, $h$ still has degree $-1$, and the relations $he=eh=h$ are unchanged. Thus the totalized algebra has the presentation
\[B_N^{\mathrm{tot}}\cong A_N^{\mathrm{tot}}\langle h\rangle/(he=eh=h),\qquad d(h)=e,\]
where the isomorphism fixes $A$ and $h$. The right-hand side is the ordinary Drinfeld model of the derived quotient in \Cref{eq:totalization-derived-quotient}.

Formula~\Cref{eq:weighted-total-cohomology} shows that totalization preserves perfectly valued objects, while $(\operatorname{Tot}_N M)e=\operatorname{Tot}_N(Me)$ preserves the condition $Me\simeq0$. The induced functors extend to the idempotent completions.
\end{proof}

For the Calabi--Yau calculation, put $\Theta_C=\RHom_{C^e}^{\mathbb Z}(C,C^e)$, with the bimodule structure given by the inner action. In the ungraded case we omit the superscript $\mathbb Z$.

\begin{Thm}\label{thm:weighted-totalization-CY}
Let $C$ be graded homologically smooth and suppose that
\begin{equation}\label{eq:weighted-CY-input}\Theta_C(a)[n]\simeq{}_{\alpha}C_1\qquad\text{in }\cd^{\mathbb Z}(C^e),\end{equation}
where $\alpha$ is a bidegree-preserving dg automorphism. Then
\begin{equation}\label{eq:weighted-CY-output}\Theta_{C_N^{\mathrm{tot}}}[n+Na]\simeq{}_{\beta_N}(C_N^{\mathrm{tot}})_1,\qquad\beta_N=\alpha^{\mathrm{tot}}\varepsilon_C^{N(a+1)}.\end{equation}
\end{Thm}

\begin{proof}
Write $R=C_N^{\mathrm{tot}}$, $\varepsilon=\varepsilon_C$, and $\mathsf T_N$ for bimodule totalization. We first construct a quasi-isomorphism
\begin{equation}\label{eq:total-dual-bimodule}\Phi_P\colon{}_{\varepsilon^N}\mathsf T_N\bigl(\RHom_{C^e}^{\mathbb Z}(P,C^e)\bigr)\longrightarrow\RHom_{R^e}(\mathsf T_N(P),R^e)\end{equation}
for every perfect bigraded bimodule $P$.

First take $P$ to be obtained from shifted free bimodules by finite direct sums and mapping cones. Both derived Hom complexes are then computed by ordinary Hom complexes. Write the outer free bimodules as $C\otimes C$ and $R\otimes R$, using the same symbols for elements before and after totalization. For a homogeneous map $f$ of bidegree $(r,s)$, write $f(m)=\sum_i x_i\otimes y_i$ and define
\begin{equation}\label{eq:explicit-total-dual-comparison}\Phi_P(f)(m)=\sum_i(-1)^{N\{r(|m|+\deg(x_i)+1)+\deg(x_i)|y_i|\}}x_i\otimes y_i.\end{equation}
This gives a homogeneous $R$-bimodule map of degree $Nr+s$. For homogeneous $u,v\in C$, the inner bimodule action is
\[(ufv)(m)=\sum_i(-1)^{|v|s+|u||v|+(|u|+|v|)|x_i|}x_i v\otimes uy_i.\]
Via $\deg(x_i)+\deg(y_i)=\deg(m)+r$ and $|x_i|+|y_i|=|m|+s$, a direct calculation gives
\[d\Phi_P(f)=\Phi_P(df),\qquad\Phi_P(\varepsilon^N(u)\star f\star v)=u\Phi_P(f)v.\]
Thus $\Phi_P$ is a map of dg $R$-bimodules. It is an isomorphism for shifted free $P$, and passage to cones and direct summands proves \Cref{eq:total-dual-bimodule} for every perfect $P$.

For a bigraded bimodule $M$, the map $m\mapsto(-1)^{Naq}m$, for $m\in M_p^q$, respects the differential and the actions after twisting the left action by $\varepsilon^{Na}$. It therefore gives
\begin{equation}\label{eq:total-bimodule-shift}\mathsf T_N(M(a)[n])\simeq{}_{\varepsilon^{Na}}\mathsf T_N(M)[Na+n].\end{equation}

Finally, taking $P=C$ in \Cref{eq:total-dual-bimodule} gives ${}_{\varepsilon^N}\mathsf T_N(\Theta_C)\simeq\Theta_R$. Using \Cref{eq:total-bimodule-shift} with $M=\Theta_C$ and $\varepsilon^2=\mathrm{id}$, we obtain
\[\mathsf T_N(\Theta_C(a)[n])\simeq{}_{\varepsilon^{Na}}\mathsf T_N(\Theta_C)[Na+n]\simeq{}_{\varepsilon^{N(a+1)}}\Theta_R[n+Na].\]

By \Cref{eq:weighted-CY-input}, the left-hand side is isomorphic to ${}_{\alpha^{\mathrm{tot}}}R_1$. Since $\varepsilon^2=\mathrm{id}$ and $\alpha^{\mathrm{tot}}$ commutes with $\varepsilon$, twisting by $\varepsilon^{N(a+1)}$ gives \Cref{eq:weighted-CY-output}.
\end{proof}

The totalizations of $A$ and $B$ are non-positive and homologically smooth. Moreover,
\begin{equation}\label{eq:totalized-H0-finite}H^0(B_N^{\mathrm{tot}})\cong\begin{cases}H^0(B)_0,&N>0,\\H^0(B),&N=0,\end{cases}\end{equation}
which is finite dimensional in both cases.

\begin{Cor}\label{cor:weighted-totalization-relative-Serre}
Assume that the Calabi--Yau structure induced on $B$ satisfied
\begin{equation}\label{eq:weighted-bimodule-hypothesis}\Theta_B(a)[d+1]\simeq{}_{\bar\alpha}B_1,\end{equation}
then
\begin{equation}\label{eq:weighted-CY-B}\Theta_{B_N^{\mathrm{tot}}}[d+1+Na]\simeq{}_{\beta_N}(B_N^{\mathrm{tot}})_1,\qquad\beta_N=\bar\alpha^{\mathrm{tot}}\varepsilon_B^{N(a+1)}.\end{equation}
The category $\cc(B_N^{\mathrm{tot}})$ is Hom-finite, with Serre functor $(-)_{\beta_N}[d+Na]$.
\end{Cor}

\begin{proof}
\Cref{thm:weighted-totalization-CY} gives \Cref{eq:weighted-CY-B} and hence the relative Serre duality on $(\per(B_N^{\mathrm{tot}}),\pvd(B_N^{\mathrm{tot}}))$. Since $B_N^{\mathrm{tot}}$ is non-positive and homologically smooth with finite-dimensional $H^0(B_N^{\mathrm{tot}})$, the argument of \Cref{Prop:relative Serre five-tuple} shows that
\[(\per(B_N^{\mathrm{tot}}),\pvd(B_N^{\mathrm{tot}}),(-)_{\beta_N}[d+1+Na],\add(B_N^{\mathrm{tot}}),0)\]
is a relative Serre five-tuple. \Cref{Thm:H is Frobenius and H/P=U/Ufd} therefore gives the Hom-finiteness of $\cc(B_N^{\mathrm{tot}})$ and its Serre functor $(-)_{\beta_N}[d+Na]$.
\end{proof}

If $\bar\alpha=\id$, the twist is trivial for even $N$. For odd $N$, it is trivial when $a$ is odd and equals $\varepsilon_B$ when $a$ is even.

\subsection{Orbit equivalences for cluster categories}

Write $\per_{\mathrm{dg}}^{\mathbb Z}(B)$ and $\per_{\mathrm{dg}}(B_N^{\mathrm{tot}})$ for the canonical dg enhancements of the corresponding perfect derived categories. We use the dg enhancements of the cluster categories given by the idempotent completions of the dg quotients, with compatible dg lifts of shifts and totalization.

\begin{Prop}\label{prop:weighted-covering-cluster-quotient}\label{cor:relative-weighted-covering-cluster}
Totalization induces natural quasi-isomorphisms
\begin{equation}\label{eq:weighted-covering-relative-cluster-dg}\bigoplus_{\ell\in\mathbb Z}\RHom_{\cc^{\mathbb Z}(A,e)}(X,F_N^\ell Y)\xrightarrow{\sim}\RHom_{\cc(A_N^{\mathrm{tot}},e)}(\operatorname{Tot}_N X,\operatorname{Tot}_N Y)\end{equation}
for $X,Y\in\cc^{\mathbb Z}(A,e)$, and
\begin{equation}\label{eq:weighted-covering-cluster-quotient-dg}\bigoplus_{\ell\in\mathbb Z}\RHom_{\cc^{\mathbb Z}(B)}(X,F_N^\ell Y)\xrightarrow{\sim}\RHom_{\cc(B_N^{\mathrm{tot}})}(\operatorname{Tot}_N X,\operatorname{Tot}_N Y)\end{equation}
for $X,Y\in\cc^{\mathbb Z}(B)$.
\end{Prop}

\begin{proof}
By \Cref{eq:weighted-covering-Hom}, totalization induces a quasi-fully faithful dg functor from the dg orbit of $\per_{\mathrm{dg}}^{\mathbb Z}(B)$ to $\per_{\mathrm{dg}}(B_N^{\mathrm{tot}})$. Since its image contains the generator $B_N^{\mathrm{tot}}$, it induces a quasi-equivalence
\[\bigl(\operatorname{pretr}(\per_{\mathrm{dg}}^{\mathbb Z}(B)/F_N^{\mathbb Z})\bigr)^\sharp\xrightarrow{\sim}\per_{\mathrm{dg}}(B_N^{\mathrm{tot}}).\]
The same argument applies to $A$.

We claim that
\begin{equation}\label{eq:proper-generated-by-totalization}\pvd(B_N^{\mathrm{tot}})=\thick\{\operatorname{Tot}_N L\mid L\in\pvd^{\mathbb Z}(B)\}.\end{equation}
The inclusion $\supseteq$ follows from \Cref{prop:totalization-derived-quotient}. As in the proof of \Cref{Thm:restriction-to-perfect}, the left-hand side is generated by the simple $H^0(B_N^{\mathrm{tot}})$-modules. By \Cref{eq:totalized-H0-finite} and the same nilpotence argument for $N=0$, these simples factor through $H^0(B)_0$. They are therefore totalizations of objects of $\pvd^{\mathbb Z}(B)$ concentrated in bidegree $(0,0)$, proving the $\subseteq$ inclusion.

By \Cref{prop:rec,prop:totalization-derived-quotient}, restriction identifies $\pvd^{\mathbb Z}(B)$ with $\pvd_e^{\mathbb Z}(A)$ and $\pvd(B_N^{\mathrm{tot}})$ with $\pvd_e(A_N^{\mathrm{tot}})$, where these restriction functors commute with totalization. Hence
\begin{equation}\label{eq:relative-proper-generated-by-totalization}\pvd_e(A_N^{\mathrm{tot}})=\thick\{\operatorname{Tot}_N L\mid L\in\pvd_e^{\mathbb Z}(A)\}.\end{equation}

Write $\pvd_{\mathrm{dg}}^{\mathbb Z}(B)$ and $\pvd_{\mathrm{dg}}(B_N^{\mathrm{tot}})$ for the full dg subcategories of perfectly valued objects. As $\pvd^{\mathbb Z}(B)$ is stable under $F_N$, the preceding quasi-equivalence and \Cref{eq:proper-generated-by-totalization}, together with \cite[Corollary~3.18]{FKQ2024}, give
\[\Bigl(\operatorname{pretr}\bigl((\per_{\mathrm{dg}}^{\mathbb Z}(B)/\pvd_{\mathrm{dg}}^{\mathbb Z}(B))/F_N^{\mathbb Z}\bigr)\Bigr)^\sharp\simeq\bigl(\per_{\mathrm{dg}}(B_N^{\mathrm{tot}})/\pvd_{\mathrm{dg}}(B_N^{\mathrm{tot}})\bigr)^\sharp,\]
where passing to the thick closure does not affect the completed quotient. This gives \Cref{eq:weighted-covering-cluster-quotient-dg} for the images of perfect objects and, by naturality, for their direct summands. Using \Cref{eq:relative-proper-generated-by-totalization}, the same argument gives \Cref{eq:weighted-covering-relative-cluster-dg}.
\end{proof}

\begin{Thm}\label{thm:relative-weighted-orbit}\label{thm:weighted-orbit}
We have the following commutative diagram of triangulated categories:
\[\begin{tikzcd}[column sep=large,row sep=large] (\cc^{\mathbb Z}(A,e)/F_N)_\triangle^\sharp \arrow[r,"\operatorname{Tot}_N","\sim"'] \arrow[d,"\pi"'] & \cc(A_N^{\mathrm{tot}},e) \arrow[d,"\pi_N"]\\
(\cc^{\mathbb Z}(B)/F_N)_\triangle^\sharp \arrow[r,"\operatorname{Tot}_N","\sim"'] & \cc(B_N^{\mathrm{tot}}). \end{tikzcd}\]
Here $\pi$ and $\pi_N$ are induced by the derived quotient maps.
\end{Thm}

\begin{proof}
The horizontal equivalences are obtained by taking $H^0$ of the quasi-equivalences constructed in the proof of \Cref{prop:weighted-covering-cluster-quotient}.

The map $x\otimes b\mapsto(-1)^{N\deg(x)|b|}x\otimes b$ induces a natural isomorphism
\[\operatorname{Tot}_N(X)\lten_{A_N^{\mathrm{tot}}}B_N^{\mathrm{tot}}\xrightarrow{\sim}\operatorname{Tot}_N(X\lten_A B).\]
This natural isomorphism descends to the dg orbit categories by compatibility with \Cref{eq:weighted-orbit-functor}, and hence the diagram commutes.
\end{proof}

The equivalence of \Cref{Thm:finite-Adams-window} identifies $F_N=(-1)[N]$ with $-\lten_\Lambda U[N]$. Thus \Cref{thm:weighted-orbit} yields the following corollary.

\begin{Cor}\label{cor:finite-window-weighted-orbit}
Assume $a>0$ and $B_0$ proper. We have the following commutative diagram:
\[ \begin{tikzcd}[column sep=large] \cc^{\mathbb Z}(B)\arrow[r,"F_N"]\arrow[d,"\sim"'] &\cc^{\mathbb Z}(B)\arrow[r,"\operatorname{Tot}_N"]\arrow[d,"\sim"'] &\cc(B_N^{\mathrm{tot}})\arrow[d,"\sim"]\\
\per\Lambda\arrow[r,"{-\lten_\Lambda U[N]}"'] &\per\Lambda\arrow[r] &\left(\per\Lambda/(-\lten_\Lambda U[N])\right)_\triangle^\sharp. \end{tikzcd} \]
In particular, there is a triangle equivalence
\begin{equation}\label{eq:finite-window-weighted-orbit-equivalence} \cc(B_N^{\mathrm{tot}}) \simeq \left(\per\Lambda/(-\lten_\Lambda U[N])\right)_\triangle^\sharp. \end{equation}
\end{Cor}

\subsection{Orbit equivalences for Higgs categories}

In this subsection we assume the bimodule identity \Cref{eq:weighted-bimodule-hypothesis} and
\begin{equation}\label{eq:totalized-Krull-Schmidt}\per(A_N^{\mathrm{tot}})\text{ is Krull--Schmidt}.\end{equation}
For $N>0$, this follows from smoothness and $H^0(A_N^{\mathrm{tot}})=H^0(A)_0$ being finite dimensional; for $N=0$, it is an additional hypothesis. Write $\ch(A_N^{\mathrm{tot}},e)$ for the ungraded Higgs category, defined using the frozen subcategory $\add(eA_N^{\mathrm{tot}})$.

\begin{Lem}\label{lem:totalized-relative-Serre}
Put $n_N=d+Na$ and $\widetilde\beta_N=\alpha^{\mathrm{tot}}\varepsilon_A^{N(a+1)}$. Then
\[(\per(A_N^{\mathrm{tot}}),\pvd_e(A_N^{\mathrm{tot}}),(-)_{\widetilde\beta_N}[n_N+1],\add A_N^{\mathrm{tot}},\add eA_N^{\mathrm{tot}})\]
is a relative Serre five-tuple. Hence $\ch(A_N^{\mathrm{tot}},e)$ is Krull--Schmidt and Frobenius extriangulated, with projective--injectives $\add(eA_N^{\mathrm{tot}})$ and stable category $\cc(B_N^{\mathrm{tot}})$, whose Serre functor is $(-)_{\beta_N}[n_N]$. If $n_N\geq1$, then $\add A_N^{\mathrm{tot}}$ is an $n_N$-cluster-tilting subcategory of $\ch(A_N^{\mathrm{tot}},e)$.
\end{Lem}

\begin{proof}
We verify the conditions of \Cref{Def:Serre five-tuple}.

For \textup{(RS1)}, $\per(A_N^{\mathrm{tot}})$ is Krull--Schmidt by \Cref{eq:totalized-Krull-Schmidt}. Smoothness and \Cref{eq:totalized-H0-finite} imply that $\per(B_N^{\mathrm{tot}})$ is Hom-finite, and hence so is $\per(A_N^{\mathrm{tot}})/\thick(eA_N^{\mathrm{tot}})$ by compact localization.
For $N>0$, $\add eA_N^{\mathrm{tot}}$ is functorially finite in $\add A_N^{\mathrm{tot}}$ since $H^0(A_N^{\mathrm{tot}})=H^0(A)_0$ is finite dimensional. For $N=0$, the same holds by \Cref{lem:graded-noetherian-functorial-finiteness}, since finite homogeneous generating sets of $H^0(A)e$ and $eH^0(A)$ over $eH^0(A)e$ are still generating sets after forgetting the grading.

For \textup{(RS2)}, every $L\in\pvd_e^{\mathbb Z}(A)$ is left and right orthogonal to $\thick^{\mathbb Z}(eA)$. By \Cref{eq:weighted-covering-Hom}, $\operatorname{Tot}_N L$ is left and right orthogonal to $\thick(eA_N^{\mathrm{tot}})$. By \Cref{eq:relative-proper-generated-by-totalization}, the same holds for every object of $\pvd_e(A_N^{\mathrm{tot}})$.

For \textup{(RS3)}, by the two-sided orthogonality in \textup{(RS2)}, the relative Serre duality follows from \Cref{eq:weighted-CY-B} via localization:
\[D\Hom(X,Y)\cong\Hom(Y,X_{\widetilde\beta_N}[n_N+1])\qquad(X\in\pvd_e(A_N^{\mathrm{tot}}),\ Y\in\per(A_N^{\mathrm{tot}})).\]
Here $\widetilde\beta_N$ preserves $e$ and induces $\beta_N$ on the quotient, so it preserves the relative proper subcategory.

For \textup{(RS4)}, as in the proof of \Cref{Prop:relative Serre five-tuple}, the Verdier quotient is idempotent complete, and hence
\[\per(A_N^{\mathrm{tot}})/\thick(eA_N^{\mathrm{tot}})\simeq\per(B_N^{\mathrm{tot}}).\]
Under this equivalence, the image of $\add A_N^{\mathrm{tot}}$ is $\add B_N^{\mathrm{tot}}$. By \Cref{Thm:restriction-to-perfect} with the trivial Adams grading and $f=0$, the standard $t$-structure restricts to $\per(B_N^{\mathrm{tot}})$ and
\[\per(B_N^{\mathrm{tot}})^{\geq0}\subseteq\pvd(B_N^{\mathrm{tot}})=\pi_N(\pvd_e(A_N^{\mathrm{tot}})).\]

For \textup{(RS5)}, the same gluing and truncation arguments as in the proof of \Cref{Prop:relative-t-on per} give the required relative $t$-structure on $\per(A_N^{\mathrm{tot}})$. Moreover, $\add B_N^{\mathrm{tot}}\simeq\operatorname{proj}H^0(B_N^{\mathrm{tot}})$ is a dualizing $k$-variety by \Cref{eq:totalized-H0-finite}.

The Frobenius property and the Serre functor follow from \Cref{Thm:H is Frobenius and H/P=U/Ufd}. The same theorem gives the $n_N$-cluster-tilting property when $n_N\geq1$, since $(-)_{\widetilde\beta_N}$ preserves $\add A_N^{\mathrm{tot}}$.
\end{proof}

Denote the chosen enhancements by $\mathcal C_{\mathrm{dg}}^{\mathbb Z}(A,e)$ and $\mathcal C_{\mathrm{dg}}(A_N^{\mathrm{tot}},e)$, where
\[\mathcal C^{\mathbb Z}_{\mathrm{dg}}(A,e):=\left(\operatorname{pretr}\left(\operatorname{per}^{\mathbb Z}_{\mathrm{dg}}(A)/\operatorname{pvd}^{\mathbb Z}_{e,\mathrm{dg}}(A)\right)\right)^{\sharp},\qquad H^0\bigl(\mathcal C^{\mathbb Z}_{\mathrm{dg}}(A,e)\bigr)\simeq\cc^{\mathbb Z}(A,e),\]
ungraded enhancement $\mathcal C_{\mathrm{dg}}(A_N^{\mathrm{tot}},e)$ is defined in the same way. Let $\widetilde{\mathcal H}^{\mathbb Z}_{f_N,\mathrm{dg}}(A,e)$ be the full dg subcategory of $\mathcal C_{\mathrm{dg}}^{\mathbb Z}(A,e)$ with objects whose images in $\cc^{\mathbb Z}(A,e)$ belong to $\ch_{f_N}^{\mathbb Z}(A,e)$. Let
\[\mathcal E=\tau_{\leq0}\widetilde{\mathcal H}_{f_N,\mathrm{dg}}^{\mathbb Z}(A,e),\qquad\mathcal H_{\mathrm{dg}}(A_N^{\mathrm{tot}},e)=\tau_{\leq0}\widetilde{\mathcal H}_{\mathrm{dg}}(A_N^{\mathrm{tot}},e)\]
be the connective truncations of their full dg Higgs subcategories, endowed with the induced exact structures of \Cref{lem:induced-exact-enhancement}. Moreover, by \Cref{Thm:CTC-in-gradedHiggs}, $H^0(\mathcal E)=\ch_{f_N}^{\mathbb Z}(A,e)$ is Frobenius, with projective--injective subcategory
\[\cp_{f_N}=\add\{P_p\mid p\in\mathbb Z\},\qquad P_p=eA(-p)[Np+c].\]
Since $F_NP_p=P_{p+1}$, the functor $F_N$ preserves the Higgs orthogonalities and restricts to an exact dg autoequivalence of $\mathcal E$, preserving its projective--injectives. We use the exact dg orbit construction of \Cref*{appB:exact-dg-orbit} and put
\[\mathcal D=\cd^b_{\mathrm{dg}}(\mathcal E),\qquad\mathcal T_{F_N}=\bigl(\operatorname{pretr}(\mathcal D/\widehat F_N^{\mathbb Z})\bigr)^\sharp,\qquad\mathscr O_{F_N}=H^0\operatorname{Orb}^{\mathrm{ex}}_{\mathrm{dg}}(\mathcal E,F_N),\]
where $\widehat F_N$ is the induced autoequivalence of $\mathcal D$.

\begin{Prop}[Normalized totalization]\label{prop:normalized-totalization-ambient-comparison}
The functor $J_N=[-c]\circ\operatorname{Tot}_N$ restricts to an exact dg functor $\mathcal E\to\mathcal H_{\mathrm{dg}}(A_N^{\mathrm{tot}},e)$, with $J_NP_p\simeq eA_N^{\mathrm{tot}}$, and induces quasi-isomorphisms
\begin{equation}\label{eq:weighted-totalization-orbit-hom}\bigoplus_{\ell\in\mathbb Z}\mathcal E(X,F_N^\ell Y)\xrightarrow{\sim}\mathcal H_{\mathrm{dg}}(A_N^{\mathrm{tot}},e)(J_NX,J_NY).\end{equation}
It extends to a quasi-equivalence
\begin{equation}\label{eq:higgs-appB-ambient-comparison}\mathcal T_{F_N}\xrightarrow{\sim}\mathcal C_{\mathrm{dg}}(A_N^{\mathrm{tot}},e).\end{equation}
Under this equivalence, $\mathscr O_{F_N}$ is the extension-closed and summand-closed subcategory $\mathcal O$ generated by $J_N\ch_{f_N}^{\mathbb Z}(A,e)$ in $\ch(A_N^{\mathrm{tot}},e)$. It is Frobenius with projective--injectives $\add(eA_N^{\mathrm{tot}})$.
\end{Prop}

\begin{proof}
By \Cref{eq:weighted-total-shifts}, we have $J_NP_p\simeq eA_N^{\mathrm{tot}}$ for all $p\in\mathbb Z$. Since $F_N^\ell P_0=P_\ell$, \Cref{eq:weighted-covering-relative-cluster-dg} gives, for $X\in\ch_{f_N}^{\mathbb Z}(A,e)$ and $j>0$,
\[\begin{aligned} \Hom_{\cc(A_N^{\mathrm{tot}},e)}(J_NX,eA_N^{\mathrm{tot}}[j])&\cong\bigoplus_{\ell\in\mathbb Z}\Hom_{\cc^{\mathbb Z}(A,e)}(X,P_\ell[j])=0,\\
\Hom_{\cc(A_N^{\mathrm{tot}},e)}(eA_N^{\mathrm{tot}},J_NX[j])&\cong\bigoplus_{\ell\in\mathbb Z}\Hom_{\cc^{\mathbb Z}(A,e)}(P_{-\ell},X[j])=0. \end{aligned}\]
The vanishing follows from the definition of $\ch_{f_N}^{\mathbb Z}(A,e)$. Hence $J_NX\in\ch(A_N^{\mathrm{tot}},e)$. Its restriction is exact by \Cref{lem:induced-exact-enhancement}, and \Cref{eq:weighted-totalization-orbit-hom} follows from \Cref{eq:weighted-covering-relative-cluster-dg}, since connective truncation commutes with direct sums.

By the definition of $\ch_{f_N}^{\mathbb Z}(A,e)$, we have $\Hom(P,X[r])=0$ for $P\in\cp_{f_N}$, $X\in\ch_{f_N}^{\mathbb Z}(A,e)$ and $r>0$. Thus, by \Cref{lem:induced-exact-enhancement}, the inclusion of $\mathcal E$ induces a quasi-fully faithful dg functor
\[\rho\colon\mathcal D\longrightarrow\mathcal C_{\mathrm{dg}}^{\mathbb Z}(A,e)\]
with $\rho\widehat F_N\simeq F_N\rho$. Since $J_NF_N\simeq J_N$, it follows from \Cref{eq:weighted-covering-relative-cluster-dg} that $J_N\rho$ induces a quasi-fully faithful dg functor
\begin{equation}\label{eq:higgs-dg-orbit-functor} \mathcal D/\widehat F_N^{\mathbb Z}\longrightarrow\mathcal C_{\mathrm{dg}}(A_N^{\mathrm{tot}},e). \end{equation}
As $\cm_{f_N}$ is presilting and the quotient functor $\per^{\mathbb Z}(A)\to\cc^{\mathbb Z}(A,e)$ preserves Hom spaces involving frozen objects, we have $\cm_{f_N}\subseteq\ch_{f_N}^{\mathbb Z}(A,e)$. In particular, the image of \Cref{eq:higgs-dg-orbit-functor} contains $J_N(A[c])\simeq A_N^{\mathrm{tot}}$, which thickly generates $\cc(A_N^{\mathrm{tot}},e)$. Passing to the pretriangulated hull and idempotent completion, the induced dg functor is also quasi-essentially surjective, and hence a quasi-equivalence, giving \Cref{eq:higgs-appB-ambient-comparison}.

Finally, under \Cref{eq:higgs-appB-ambient-comparison}, \Cref{lem:orbit-extension-closure} identifies $\mathscr O_{F_N}$ with the extension-and-summand closure $\mathcal O$ of $J_N\ch_{f_N}^{\mathbb Z}(A,e)$ in $\ch(A_N^{\mathrm{tot}},e)$. This is an equivalence of extriangulated categories since both structures are induced by triangles in the corresponding triangulated categories. By \Cref{thm:frobenius-exact-orbit} and $J_NP_p\simeq eA_N^{\mathrm{tot}}$ for all $p\in\mathbb Z$, $\mathcal O$ is Frobenius with
\[\operatorname{proj}\mathcal O=\operatorname{inj}\mathcal O=\add(eA_N^{\mathrm{tot}}).\]
\end{proof}

\begin{Cor}\label{cor:weighted-totalization-higgs-orbit}
Under \Cref{Ass} for $f=f_N$ and conditions \Cref{eq:weighted-bimodule-hypothesis} and \Cref{eq:totalized-Krull-Schmidt}, write
\[(\ch_{f_N}^{\mathbb Z}(A,e)/F_N)^{\mathrm{ex},\sharp}:=\mathscr O_{F_N}.\]
The functor $J_N=[-c]\circ\operatorname{Tot}_N$ is fully faithful on the ordinary orbit and induces an exact equivalence
\begin{equation}\label{eq:higgs-orbit-equivalence} (\ch_{f_N}^{\mathbb Z}(A,e)/F_N)^{\mathrm{ex},\sharp}\xrightarrow{\sim}\ch(A_N^{\mathrm{tot}},e). \end{equation}
The stable category is triangle equivalent to $\cc(B_N^{\mathrm{tot}})$.
\end{Cor}

\begin{proof}
Taking $H^0$ in \Cref{eq:weighted-totalization-orbit-hom} gives full faithfulness on the ordinary orbit. By \Cref{prop:normalized-totalization-ambient-comparison}, it remains to show $\mathcal O=\ch(A_N^{\mathrm{tot}},e)$.

Since both categories are Frobenius with the same projective--injectives $\add(eA_N^{\mathrm{tot}})$, the inclusion $\mathcal O\hookrightarrow\ch(A_N^{\mathrm{tot}},e)$ induces a fully faithful triangle functor
\begin{equation}\label{eq:higgs-stable-inclusion} \underline{\mathcal O}\longrightarrow\underline{\ch(A_N^{\mathrm{tot}},e)}\simeq\cc(B_N^{\mathrm{tot}}). \end{equation}
where the equivalence follows from \Cref{lem:totalized-relative-Serre}.

Suppose that $X\in\ch(A_N^{\mathrm{tot}},e)$ is a direct summand of $Y\in\mathcal O$ in the stable category. Then there are morphisms $f\colon X\to Y$ and $g\colon Y\to X$ such that $1_X-gf=vu$ for some $u\colon X\to P$ and $v\colon P\to X$ with $P\in\add(eA_N^{\mathrm{tot}})$. Hence $X$ is a direct summand of $Y\oplus P\in\mathcal O$, so $X\in\mathcal O$, by the summand-closed property of $\mathcal O$.

Thus the essential image of $\underline{\mathcal O}$ is a thick subcategory of $\cc(B_N^{\mathrm{tot}})$. It contains the generator $B_N^{\mathrm{tot}}$, since $J_N(A[c])\simeq A_N^{\mathrm{tot}}$ belongs to $\mathcal O$. Hence the functor in \Cref{eq:higgs-stable-inclusion} is essentially surjective. By the preceding argument, $\mathcal O=\ch(A_N^{\mathrm{tot}},e)$, which proves \Cref{eq:higgs-orbit-equivalence} and hence the triangle equivalence of stable categories.
\end{proof}

\section{Multigraded Higgs categories}\label{sec:other-gradings}

We extend the constructions of \Cref{sec:graded-t-structure,sec:graded-higgs} to gradings by an abelian group. The main difference lies in the finiteness conditions on the associated $t$-structures and silting subcategories.

\subsection{Silting subcategories and finiteness conditions}

Let $L$ be an abelian group and $\Sigma\subseteq L$ a submonoid containing $0$. Throughout this section, $A$ is a differential $L$-graded algebra satisfying
\begin{equation}\label{eq:Lambda-third-quadrant}A_\lambda^q=0\qquad\text{if }\lambda\notin\Sigma\text{ or }q>0.\end{equation}
We use the grading shift $M(\rho)_\lambda^q=M_{\lambda+\rho}^q$ and the notation $\cd^L(A)$, $\per^L(A)$ and $\pvd^L(A)$ as in the $\mathbb Z$-graded case.

A function $f\colon L\to\mathbb Z$ is called \emph{$\Sigma$-monotone} if
\begin{equation}\label{eq:Sigma-monotone}f(\lambda+\sigma)\le f(\lambda)\qquad(\lambda\in L,\ \sigma\in\Sigma).\end{equation}
For $L=\mathbb Z$ and $\Sigma=\mathbb Z_{\le0}$, this is precisely the condition that $f$ is non-decreasing.

Put $M_\lambda=A(-\lambda)[f(\lambda)]$ and
\begin{equation}\label{eq:Mf-Lambda}\cm_f^L=\add\{M_\lambda\mid\lambda\in L\}.\end{equation}
The Hom calculation in \Cref{prop:graded-silting-subcategory} gives
\begin{equation}\label{eq:Lambda-Hom-M-M}\Hom(M_\lambda,M_\mu[s])\cong H^{f(\mu)-f(\lambda)+s}(A)_{\lambda-\mu}.\end{equation}
If $\lambda-\mu\in\Sigma$, then $f(\lambda)\le f(\mu)$. Hence $\cm_f^L$ is presilting, and since $\thick(\cm_f^L)=\per^L(A)$, it is silting.

The associated $t$-structure on $\cd^L(A)$ is given by
\[\cd_f^{\le n}=\{X\mid H^q(X)_\lambda=0\text{ if }q+f(\lambda)>n\},\qquad \cd_f^{\ge n}=\{X\mid H^q(X)_\lambda=0\text{ if }q+f(\lambda)<n\}.\]
Indeed, the construction in the proof of \Cref{lem:t-stru-onD(A)} applies with \Cref{eq:Sigma-monotone} in place of monotonicity on $\mathbb Z$. In particular,
\begin{equation}\label{eq:Lambda-f-truncation}H^q(\tau_f^{\ge n}X)_\lambda\cong\begin{cases}H^q(X)_\lambda&\text{if }q+f(\lambda)\ge n,\\0&\text{otherwise}.\end{cases}\end{equation}
The truncation formula above yields the following theorem by the argument for \Cref{Thm:restriction-to-perfect}\textup{(\Cref{item:perfect-truncation-all-degrees})}$\Leftrightarrow$\textup{(\Cref{item:perfect-truncation-finite-sum})}, with grading shifts indexed by $L$. 
\begin{Thm}\label{thm:Lambda-f-properness}
Let $f\colon L\to\mathbb Z$ be $\Sigma$-monotone. The following conditions are equivalent:
\begin{enumerate}
\item $\tau_f^{\ge n}(\per^L(A))\subseteq\pvd^L(A)$ for every $n\in\mathbb Z$;
\item\label{item:Lambda-f-properness-finite-sum} for every $\rho\in L$ and $n\in\mathbb Z$,
\begin{equation}\label{eq:Lambda-f-properness}\sum_{\substack{\lambda\in L,\ q\le0\\q+f(\lambda-\rho)\ge n}}\dim_k H^q(A)_\lambda<\infty.\end{equation}
\end{enumerate}
If, in addition, $\pvd^L(A)\subseteq\per^L(A)$, then the $t$-structure restricts to $\per^L(A)$, with $\per^L(A)\cap\cd_f^{\ge n}\subseteq\pvd^L(A)$ for every $n\in\mathbb Z$.
\end{Thm}

Suppose that $\dim_k H^q(A)_\lambda<\infty$ for all $\lambda,q$. A sufficient condition for \Cref{eq:Lambda-f-properness} is that, for every $\rho\in L$ and $n\in\mathbb Z$, the set
\begin{equation}\label{eq:f-superlevel-finite}\{\lambda\in\operatorname{Supp}_L H^*(A)\mid f(\lambda-\rho)\ge n\}\end{equation}
is finite. In fact, for each such $\lambda$, only the cohomological degrees between $n-f(\lambda-\rho)$ and $0$ contribute to the sum. This applies to $L=\mathbb Z^m$, $\Sigma=-\mathbb N^m$ and $f(\lambda)=\langle\mathbf w,\lambda\rangle+c$, with $w_i>0$ for every $i$. 

We next consider the dualizing property of the silting subcategory.

\begin{Lem}\label{lem:Lambda-dualizing-criterion}
Assume that $\dim_k H^0(A)_\lambda<\infty$ for every $\lambda\in L$. Then $\cm_f^L$ is a dualizing $k$-variety if, for every $\lambda\in L$, both sets
\begin{equation}\label{eq:Lambda-dualizing-criterion} \begin{aligned} &\{\sigma\in\operatorname{Supp}_L H^0(A)\mid f(\lambda+\sigma)=f(\lambda)\},\\
&\{\sigma\in\operatorname{Supp}_L H^0(A)\mid f(\lambda-\sigma)=f(\lambda)\} \end{aligned} \end{equation}
are finite. The converse holds if $(\lambda+\Sigma)\cap(\mu-\Sigma)$ is finite for all $\lambda,\mu\in L$. In particular, the converse holds for $L=\mathbb Z^m$ and $\Sigma=-\mathbb N^m$.
\end{Lem}

\begin{proof}
By \Cref{eq:Lambda-Hom-M-M} and \Cref{eq:Sigma-monotone}, we have
\[\Hom(M_\lambda,M_\mu)\cong\begin{cases}H^0(A)_{\lambda-\mu}&\text{if }f(\lambda)=f(\mu),\\0&\text{otherwise}.\end{cases}\]
Thus \Cref{eq:Lambda-dualizing-criterion} is equivalent to requiring that both covariant and contravariant representable functors have finite support on $(M_\lambda)_{\lambda\in L}$. Since their values are finite dimensional, the sufficiency follows by the same argument as in the proof of \Cref{lem:dualizing-criterion}.

Conversely, suppose that $\cm_f^L$ is dualizing. For each $\lambda\in L$, the functor $D\Hom(M_\lambda,-)$ is finitely generated, so there is a surjection
\[\bigoplus_{i=1}^s\Hom(-,M_{\mu_i})\twoheadrightarrow D\Hom(M_\lambda,-).\]
By the Hom formula above, the support of the target is contained in
\[(\lambda-\Sigma)\cap\bigcup_{i=1}^s(\mu_i+\Sigma),\]
which is finite by hypothesis. This proves the finiteness of the second set in \Cref{eq:Lambda-dualizing-criterion}. The same argument in the opposite category proves the finiteness of the first.
\end{proof}

If \Cref{eq:Lambda-f-properness} holds, the first set in \Cref{eq:Lambda-dualizing-criterion} is automatically finite, by taking $\rho=-\lambda$ and $n=f(\lambda)$. Thus only the second set needs to be considered.

For affine functions, the properness and dualizing conditions admit a common description.
\begin{Prop}\label{prop:Lambda-affine-finiteness}
Assume that $A$ is graded homologically smooth. Let $f=\ell+c$, where $\ell:L\to\mathbb Z$ is a group homomorphism with $\ell(\Sigma)\subseteq\mathbb Z_{\le0}$ and $c\in\mathbb Z$. Then \Cref{eq:Lambda-f-properness} holds if and only if
\begin{equation}\label{eq:Lambda-zero-weight-finiteness}\dim_k\bigoplus_{\ell(\lambda)=0}H^0(A)_\lambda<\infty.\end{equation}
Under this condition, $\cm_f^L$ is a dualizing $k$-variety.
\end{Prop}

\begin{proof}
Necessity follows by taking $\rho=0$ and $n=c$ in \Cref{eq:Lambda-f-properness}. Conversely, group together the grading along $\ell$:
\[\widetilde A_p^q=\bigoplus_{\ell(\lambda)=p}A_\lambda^q.\]
This preserves graded homological smoothness, since shifted free bimodules remain shifted free. The algebra $\widetilde A$ is concentrated in the third quadrant and has finite-dimensional $H^0(\widetilde A)_0$ by \Cref{eq:Lambda-zero-weight-finiteness}. The diagonal case of \Cref{Thm:restriction-to-perfect} therefore gives
\[\sum_{q+\ell(\lambda)\ge n}\dim_k H^q(A)_\lambda<\infty\qquad(n\in\mathbb Z).\]
Since $f(\lambda-\rho)=\ell(\lambda)+c-\ell(\rho)$, this proves \Cref{eq:Lambda-f-properness}.

Each $H^0(A)_\lambda$ is now finite dimensional, and both sets in \Cref{eq:Lambda-dualizing-criterion} equal $\{\sigma\in\operatorname{Supp}_L H^0(A)\mid\ell(\sigma)=0\}$, which is finite by \Cref{eq:Lambda-zero-weight-finiteness}. \Cref{lem:Lambda-dualizing-criterion} proves the last assertion.
\end{proof}

In particular, for $L=\mathbb Z^m$ and $\Sigma=-\mathbb N^m$, applying \Cref{prop:Lambda-affine-finiteness} to $\ell(\lambda)=\sum_i\lambda_i$ shows that graded homological smoothness and $\dim_k H^0(A)_0<\infty$ imply $\dim_k H^q(A)_\lambda<\infty$ for all $\lambda,q$.

Condition~\Cref{eq:f-superlevel-finite} is not necessary. For example, let $L=\mathbb Z^2$, $\Sigma=-\mathbb N^2$ and $A=k[x,y]$ with zero differential and
\[\operatorname{bideg}(x)=((-1,0),0),\qquad\operatorname{bideg}(y)=((0,-1),-2).\]
For $f(\lambda_1,\lambda_2)=\lambda_1$, we have $\bigoplus_{\lambda_2\in\mathbb Z}H^0(A)_{(0,\lambda_2)}=k$. Hence \Cref{eq:Lambda-f-properness} holds by \Cref{prop:Lambda-affine-finiteness}, whereas the set in \Cref{eq:f-superlevel-finite} contains $(0,-j)$ for every $j\ge0$ when $\rho=0$ and $n=0$.
\subsection{Multigraded Higgs categories}

Let $e\in A_0^0$ be a closed idempotent and put $B=\overline A$. We retain the notation of \Cref{subsec:relative-t-structures}, replacing $\mathbb Z$ by $L$. In particular,
\[\cp_f^L=\add\{eA(-\lambda)[f(\lambda)]\mid\lambda\in L\},\qquad \overline{\cm}_f^L=\add\{B(-\lambda)[f(\lambda)]\mid\lambda\in L\}.\]
The graded relative cluster category $\cc^L(A,e)$ and the graded Higgs category $\ch_f^L(A,e)$ are defined as in \Cref{Def:grcc-ghiggs}.

By the same argument as in \Cref{lem:graded-noetherian-functorial-finiteness}, $\cp_f^L$ is functorially finite in $\cm_f^L$ if $H^0(A)$ is graded left and right Noetherian. The same holds if $\dim_k H^0(A)_0<\infty$ and $f(\lambda+\sigma)<f(\lambda)$ for all $\lambda\in L$ and $0\ne\sigma\in\Sigma$.

\begin{Thm}\label{thm:Lambda-graded-Higgs}
Let $A$ be graded homologically smooth, with $\dim_k H^q(A)_\lambda<\infty$ for every $\lambda,q$. Let $f$ be $\Sigma$-monotone and assume that
\begin{enumerate}
\item\label{item:Lambda-Higgs-properness} $B$ satisfies \Cref{eq:Lambda-f-properness};
\item\label{item:Lambda-Higgs-dualizing} $\overline{\cm}_f^L$ is a dualizing $k$-variety;
\item\label{item:Lambda-Higgs-frozen-finiteness} $\cp_f^L$ is functorially finite in $\cm_f^L$.
\end{enumerate}
Suppose further that, for some $L$-graded dg automorphism $\alpha$ of $A$ with $\alpha(e)=e$, $\mathbf a\in L$ and $d\in\mathbb Z$, the functor
\begin{equation}\label{eq:Lambda-relative-Serre}\mathbb S=(-)_\alpha(\mathbf a)[d+1]\end{equation}
restricts to an autoequivalence of $\pvd_e^L(A)$ and gives bifunctorial isomorphisms
\[D\Hom_{\per^L(A)}(X,Y)\cong\Hom_{\per^L(A)}(Y,\mathbb S X)\qquad(X\in\pvd_e^L(A),\ Y\in\per^L(A)).\]
Then the following assertions hold.
\begin{enumerate}
\item\label{item:Lambda-Higgs-Frobenius} The category $\ch_f^L(A,e)$ is Frobenius extriangulated with projective--injective subcategory $\cp_f^L$. There is a triangle equivalence
\[\ch_f^L(A,e)/[\cp_f^L]\xrightarrow{\sim}\cc^L(B),\]
and the stable Serre functor is induced by $\mathbb S[-1]$.
\item\label{item:Lambda-Higgs-cluster-tilting} If there exists $r\in\mathbb Z$ such that
\begin{equation}\label{eq:Lambda-f-Serre-compatibility}f(\lambda)-f(\lambda-\mathbf a)=r\qquad(\lambda\in L)\end{equation}
and $d+r\ge1$, then $\cm_f^L$ is a $(d+r)$-cluster-tilting subcategory of $\ch_f^L(A,e)$.
\end{enumerate}
\end{Thm}

\begin{proof}
The proof of \Cref{Prop:relative Serre five-tuple} applies using \Cref{thm:Lambda-f-properness} for the truncation condition. In the factorization step, \Cref{eq:Sigma-monotone} gives $f(\nu)=f(\lambda)$ whenever $\lambda-\nu,\nu-\mu\in\Sigma$ and $f(\lambda)=f(\mu)$. This yields the required relative Serre five-tuple. Since $\alpha(e)=e$ implies $\mathbb S^{\pm1}(\cp_f^L)\subseteq\thick(\cp_f^L)$, \Cref{Thm:H is Frobenius and H/P=U/Ufd} gives assertion~\textup{(\Cref{item:Lambda-Higgs-Frobenius})}.

Under \Cref{eq:Lambda-f-Serre-compatibility}, we have
\[\mathbb S[-d-r-1](M_\lambda)\simeq A_\alpha(-(\lambda-\mathbf a))[f(\lambda-\mathbf a)]\simeq M_{\lambda-\mathbf a}.\]
Thus $\mathbb S[-d-r-1](\cm_f^L)=\cm_f^L$, and the same theorem gives assertion~\textup{(\Cref{item:Lambda-Higgs-cluster-tilting})}.
\end{proof}

By \Cref{prop:Lambda-affine-finiteness}, for an affine function $f=\ell+c$, conditions~\textup{(\Cref{item:Lambda-Higgs-properness})} and~\textup{(\Cref{item:Lambda-Higgs-dualizing})} in the hypotheses of \Cref{thm:Lambda-graded-Higgs} can be replaced by the single condition
\begin{equation}\label{eq:Lambda-reduced-zero-weight-finiteness}\dim_k\bigoplus_{\ell(\lambda)=0}H^0(B)_\lambda<\infty.\end{equation}
Moreover, $f(\lambda)-f(\lambda-\mathbf a)=\ell(\mathbf a)$, so we obtain the following consequence.

\begin{Cor}\label{cor:Lambda-linear-cluster-dimension}
Under the hypotheses of \Cref{thm:Lambda-graded-Higgs}, let $f=\ell+c$ be affine. If $d+\ell(\mathbf a)\ge1$, then $\cm_f^L$ is a $(d+\ell(\mathbf a))$-cluster-tilting subcategory of $\ch_f^L(A,e)$.

In particular, for $L=\mathbb Z^m$, $\Sigma=-\mathbb N^m$ and $\ell(\lambda)=\langle\mathbf w,\lambda\rangle$ with $w_i>0$ for every $i$, conditions~\textup{(\Cref{item:Lambda-Higgs-properness})--(\Cref{item:Lambda-Higgs-frozen-finiteness})} in the hypotheses of \Cref{thm:Lambda-graded-Higgs} follow from graded homological smoothness and $\dim_k H^0(A)_0<\infty$.
\end{Cor}

\appendix
\section{Relative Serre five-tuple}\label[appendix]{app:relative-Serre-five-tuple}
	
Let $\ct$ be a triangulated category. Let $ \mathcal{P} $ be a pre-silting subcategory of $ \mathcal{T} $. Let $ \mathcal{S} $ be the thick subcategory $\thick\mathcal{P} $ of $ \ct $ and $ \mathcal{U} $ the quotient category $ \mathcal{T}/\mathcal{S} $. We call $ \mathcal{U} $ the \emph{silting reduction} of $ \mathcal{T} $ with respect to $ \mathcal{P} $ (see \cite{IY2018}). Denote by $\pi\colon\ct\ra\cu=\ct/\cs$ the canonical quotient functor.

For an integer $ l $, there is a bounded co-t-structure $ (\mathcal{S}_{\geqslant l},\mathcal{S}_{\leqslant l}) $ on $ \mathcal{S} $, where
	$$ \mathcal{S}_{\geqslant l}=\mathcal{S}_{>l-1}\coloneqq\bigcup_{i\geqslant0}\mathcal{P}[-l-i]*\cdots*\mathcal{P}[-l-1]*\mathcal{P}[-l],$$
	$$ \mathcal{S}_{\leqslant l}=\mathcal{S}_{<l+1}\coloneqq\bigcup_{i\geqslant0}\mathcal{P}[-l]*\mathcal{P}[-l+1]\cdots*\mathcal{P}[-l+i].$$
	
	Let $ \cz $ be the following subcategory of $ \ct $
	$$ \mathcal{Z}={}^{\perp_{\mathcal{T}}}(\mathcal{S}_{<0})\cap(\mathcal{S}_{>0})^{\perp_{\mathcal{T}}}={}^{\perp_{\mathcal{T}}}(\mathcal{P}[>0])\cap(\mathcal{P}[<0])^{\perp_{\mathcal{T}}} .$$

	\begin{Prop}\cite[Lemma 3.4]{IY2018}\label{Prop:siltingreduction}
		The composition $\cz\subseteq\ct\ra\cu$ induces a fully faithful embedding
		$$\pi\colon\cz/[\cp]\hookrightarrow\cu.$$
	\end{Prop}
	
	\begin{Rem}
	In the proof of \cite[Lemma 3.4]{IY2018}, we don't need conditions (P1) and (P2) in \cite[Section 3.1]{IY2018}.
	\end{Rem}

	Let $\cm$ be a silting subcategory of $\ct$ and assume that $\cp$ is a subcategory of $\cm$. By the proof of \cite[Theorem 3.7]{IY2018}, the image of $\cm$ in $\cu$ is a silting subcategory of $\cu$. By \cite[Proposition 2.8]{IY2018}, the pair $(\ct_{\geqslant0},\ct_{\leqslant0})$ is a bounded co-t-structure on $\ct$, where, for an integer $ l $,
	$$ \mathcal{T}_{\geqslant l}=\mathcal{T}_{>l-1}\coloneqq\bigcup_{i\geqslant0}\mathcal{M}[-l-i]*\cdots*\mathcal{M}[-l-1]*\mathcal{M}[-l],$$
	$$ \mathcal{T}_{\leqslant l}=\mathcal{T}_{<l+1}\coloneqq\bigcup_{i\geqslant0}\mathcal{M}[-l]*\mathcal{M}[-l+1]\cdots*\mathcal{M}[-l+i].$$
	Similarly, the pair $(\cu_{\geqslant0},\cu_{\leqslant0})$ is a bounded co-t-structure on $\cu$ defined by the image of $\cm$ in $\cu$.

\begin{Def}\rm\label{Def:Serre five-tuple}

Let $\ct^{\fd}$ be a thick subcategory of $\ct$. 

		We say that $(\ct,\ct^{\fd},\mathbb{S},\cm,\cp)$ is a \emph{relative Serre five-tuple} if the following conditions are satisfied
		\begin{itemize}
			\item[(RS1)] The category $\ct$ is Krull--Schmidt, $\cp$ is functorially finite in $\cm$ and $\cu=\ct/\thick(\cp)$ is Hom-finite.
			\item[(RS2)] $\thick(\cp)$ is left orthogonal to $\ct^{\mathrm{fd}}$, i.e. $\ct^{\fd}\subseteq\thick(\cp)^{\perp_{\mathcal{T}}}$.
			\item[(RS3)] $\mathbb{S}\colon \ct^{\fd}\ra\ct^{\fd} $ is a triangle equivalence such that there exists a bifunctorial isomorphism
			\[ D \operatorname{Hom}_{\mathcal{T}}(X, Y) \simeq \operatorname{Hom}_{\mathcal{T}}(Y, \mathbb{S}(X)) \]
			for any $X\in\ct^{\fd}$ and $Y\in\ct$.
			\item[(RS4)] Put $\overline{\cm}:=\add\pi(\cm)\subseteq\cu$. $\cu$ admits a $t$-structure
\[ (\cu^{\leqslant0},\cu^{\geqslant0}) \coloneqq \bigl(\overline{\cm}[<\!0]^{\perp_{\mathcal U}},\overline{\cm}[>\!0]^{\perp_{\mathcal U}}\bigr) \]
with $\cu^{\geqslant0}\subseteq\pi(\ct^{\fd})$. 
			\item[(RS5)] $\ct$ admits a $t$-structure $$(\ct^{\leqslant0},\ct^{\geqslant0})$$ 
			with $\ct^{\leqslant0}=\pi^{-1}(\cu^{\leqslant0})$, $\ct^{\geqslant0}=\cm[>\!0]^{\perp_{\mathcal{T}}}\cap\thick(\cp)^{\perp_{\mathcal{T}}}\cap\,^{\perp_{\mathcal{T}}}\thick(\cp)$ and $\ct^{\geqslant0}\subseteq\ct^{\fd}$. Moreover, $\overline{\cm}$ is a dualizing $k$-variety.
		\end{itemize}
	\end{Def}
By \textup{(RS2)} and \textup{(RS3)}, $\ct^{\mathrm{fd}}$ is left and right orthogonal to $\thick(\cp)$.
\begin{Ex}
Let $(\ct,\ct^{\fd},\mathbb{S},\cm)$ be a relative Serre quadruple in the sense of Iyama-Yang \cite[Definition 2.10]{IY2020}. Then $(\ct,\ct^{\fd},\mathbb{S},\cm,\cp=0)$ is a relative Serre five-tuple.

\end{Ex}

\bigskip

Let $(\ct,\ct^{\fd},\mathbb{S},\cm,\cp)$ be a relative Serre five-tuple. 
Let $\cu^{\fd}=\pi(\ct^{\fd})$. Then $(\cu,\cu^{\fd},\mathbb{S},\overline{\cm})$ is a relative Serre quadruple in the sense of \cite[Definition 2.10]{IY2020} by (RS1)--(RS5).

Indeed, the two-sided orthogonality to $\thick(\cp)$ identifies Hom spaces whenever one argument belongs to $\ct^{\fd}$. Thus $\pi$ identifies $\ct^{\fd}$ with the thick subcategory $\cu^{\fd}$ and transports the relative Serre duality. The image $\overline{\cm}$ is silting as noted above. It is idempotent complete by the dualizing hypothesis, so \cite[Theorem~2.9]{IY2018} makes $\cu$ idempotent complete; together with (RS1), this makes $\cu$ Hom-finite and Krull--Schmidt. The remaining conditions are (RS4)--(RS5).

\begin{Lem}
The assumptions of \Cref{Def:Serre five-tuple} imply
\begin{equation}\label{eq:frozen-approximation-halves} \begin{aligned} &\cp\text{ is contravariantly finite in }\cp[<0]^{\perp_{\ct}},\\
&\cp\text{ is covariantly finite in }{}^{\perp_{\ct}}\cp[>0]. \end{aligned} \end{equation}
In particular, every object of $\cz$ has both frozen approximations.
\end{Lem}
\begin{proof}
Write
\[ \begin{aligned} \ct_+&=\cm[<0]^{\perp_{\ct}} =\bigcup_{m\geq0}\cm*\cm[1]*\cdots*\cm[m],\\
\ct_-&={}^{\perp_{\ct}}\cm[>0] =\bigcup_{m\geq0}\cm[-m]*\cdots*\cm[-1]*\cm. \end{aligned} \]
These are the two halves of the bounded co-$t$-structure associated with $\cm$.

Let $X\in\cp[<0]^{\perp_{\ct}}$. The relative truncation gives a triangle
\[ Y\xrightarrow{u}X\longrightarrow L\longrightarrow Y[1],\qquad Y\in\ct^{\leq0},\quad L\in\ct^{\geq1}\subseteq\ct^{\fd}. \]
Since $L$ is orthogonal to every shift of $\cp$, the map $u$ induces isomorphisms on all functors $\Hom(\cp,-[j])$. In particular $Y\in\cp[<0]^{\perp_{\ct}}$. For $M\in\cm$ and $j>0$, the positive-degree comparison in \cite[Lemma~3.4]{IY2018} gives
\[ \Hom_{\ct}(M,Y[j])\cong\Hom_{\cu}(\pi M,\pi Y[j])=0, \]
because $\pi Y\in\cu^{\leq0}$. The Hom comparison applies since $M\in{}^{\perp}\cp[>0]$ and $Y\in\cp[<0]^{\perp}$. Hence $Y\in\ct_+$. Its finite co-$t$-decomposition provides a triangle
$Y_1\to M_0\xrightarrow{v}Y\to Y_1[1]$ with $M_0\in\cm$ and $Y_1\in\ct_+$. Thus $\Hom(\cp,v)$ is surjective. Choose a right $\cp$-approximation $P_0\to M_0$ using (RS1). The composite
\[ P_0\longrightarrow M_0\xrightarrow{v}Y\xrightarrow{u}X \]
is a right $\cp$-approximation of $X$.

For the other direction let $X\in{}^{\perp_{\ct}}\cp[>0]$. The reduced relative Serre quadruple has the left adjacent $t$-structure supplied by \cite[Theorem~4.10]{IY2018}. The cut with coaisle $\cu_-={}^{\perp_{\cu}}\overline{\cm}[>0]$ gives a triangle
\[ L\longrightarrow\pi X\longrightarrow V\longrightarrow L[1],\qquad L\in{}^{\perp_{\cu}}\overline{\cm}[\leq0]\subseteq\cu^{\fd},\quad V\in\cu_-. \]
Choose $L_A\in\ct^{\fd}$ with $\pi L_A=L$. The Hom comparison for a proper argument lifts the first map to $L_A\to X$; let $Y$ be its cone. Then $\pi Y\simeq V$, and $X\to Y$ induces isomorphisms on all $\Hom(-,\cp[j])$. In particular $Y\in{}^{\perp}\cp[>0]$. The same positive-degree comparison, now with source $Y$ and target $M\in\cm$, gives
\[ \Hom_{\ct}(Y,M[j])\cong\Hom_{\cu}(V,\pi M[j])=0\qquad(j>0). \]
Thus $Y\in\ct_-$. A finite co-$t$-decomposition gives a triangle
$K\to Y\xrightarrow{w}M^0\to K[1]$ with $M^0\in\cm$ and $K$ a finite extension of negative shifts of objects of $\cm$. Since $\Hom(K,\cp)=0$, every map $Y\to P$, $P\in\cp$, factors through $w$. Choose a left $\cp$-approximation $M^0\to P^0$ using (RS1). The composite
\[ X\longrightarrow Y\xrightarrow{w}M^0\longrightarrow P^0 \]
is a left $\cp$-approximation of $X$. The dualizing condition provides the left adjacent structure on the reduction used in the second construction.
\end{proof}

The two conclusions in \Cref{eq:frozen-approximation-halves} are precisely condition (P1) of \cite[Section~3.1]{IY2018}. Condition (P2) follows from $\ct=\thick(\cm)$ and $\cp\subseteq\cm$: for each $X\in\ct$, both $\Hom(X,\cp[j])$ and $\Hom(\cp,X[j])$ vanish for all sufficiently large $j$. Hence silting reduction gives
\begin{equation}\label{eq:full-silting-reduction} \cz/[\cp]\xrightarrow{\sim}\cu=\ct/\thick(\cp). \end{equation}
The equivalence is triangulated by \cite[Theorem~3.6]{IY2018}. The quotient $\cm/[\cp]$ is Krull--Schmidt: the nonzero images of its indecomposable objects retain local endomorphism rings. The degree-zero Hom comparison therefore identifies $\cm/[\cp]$ with $\overline{\cm}$; in particular, each object of $\overline{\cm}$ has a lift in $\cm$.

Let $\cx_{\cu}^{0}=\overline{\cm}[<0]^{\perp_{\cu}}\cap\cu^{\fd}$. We can view $\cx_{\cu}^{0}$ as a full subcategory of $\ct^{\fd}$. Define $\ct^{'\leqslant0}=\cx_{\cu}^{0}$ and $\ct^{'\geqslant0}=\pi^{-1}((\cx_{\cu}^{0})^{\perp})[1]$, where $\pi\colon\ct\ra\cu$ is the quotient functor.

Put $\cx_{\cu}^{1}=\cu^{\fd}\cap\cu^{\leq1}=\cx_{\cu}^{0}[-1]$. We use the same notation for its full preimage in $\ct^{\fd}$.

\begin{samepage}
\begin{Lem}\label{Lem:left t-structure}
Assume that one of the following conditions holds.
\begin{itemize}
	\item[(a)] $\mathbb{S}\colon\ct^{\fd}\ra\ct^{\fd}$ extends to a triangle equivalence $\mathbb{S}\colon\ct\ra\ct$ such that $\mathbb{S}^{\pm}(\cp)\subseteq\thick(\cp)$;
	\item[(b)] $\cm$ has an additive generator $M$.	
\end{itemize}	
	The pair $(\ct^{'\leqslant0},\ct^{'\geqslant0})$ is a $t$-structure on $\ct$ with $\ct^{'\leqslant0}\subseteq\ct^{\fd}$.
\end{Lem}
\end{samepage}

\begin{proof}
By our assumption, the Serre quadruple $(\cu,\cu^{\fd},\mathbb{S},\overline{\cm})$ satisfies condition (i) or (ii) of \cite[Theorem 2.12]{IY2020}. By the proof of \cite[Theorem 2.12]{IY2020}, the category $\cu$ admits another $t$-structure $(\cx_{\cu}^{0},(\cx_{\cu}^{0})^{\perp}[1])$ where $\cx_{\cu}^{0}=\overline{\cm}[<0]^{\perp_{\cu}}\cap\cu^{\fd}\subseteq\cu^{\fd}$. 

It is clear that $\Hom_{\ct}(\ct^{'\leqslant0},\ct^{'\geqslant1})$ vanishes. Then it is enough to show 
$$\ct=\ct^{'\leqslant0}*\pi^{-1}((\cx_{\cu}^{0})^{\perp}).$$

Since $\ct=\bigcup_{l\geqslant0}\ct_{\geqslant-l}$, it is sufficient to show $\ct_{\geqslant -l}\subseteq\cx_{\cu}^{0}*\pi^{-1}((\cx_{\cu}^{0})^{\perp})$ for any $l\geqslant0$.  

Let $X$ be an object of $\ct$. Then we have a triangle in $\cu$
$$X_{1}\xrightarrow{\alpha} \pi(X)\ra X_{2}\ra X_{1}[1]$$
with $X_{1}\in \cx_{\cu}^{0}$ and $X_{2}\in(\cx_{\cu}^{0})^{\perp}$. Since $\Hom_{\cu}(X_{1},\pi(X))\cong\Hom_{\ct}(X_{1},X)$, there exists a map $\tilde{\alpha}\colon X_{1}\ra X$ in $\ct$ such that $\pi(\tilde{\alpha})=\alpha$. Therefore we have the following triangle in $\ct$
$$X_{1}\xrightarrow{\tilde{\alpha}}X\ra X'_{2}\ra X_{1}[1]$$ with $X_{1}\in\ct^{'\leqslant0}$ and $\pi(X'_{2})\cong X_{2}\in(\cx_{\cu}^{0})^{\perp}$. This means that $X'_{2}$ lies in $\pi^{-1}((\cx_{\cu}^{0})^{\perp})$. Hence $\ct=\ct^{'\leqslant0}*\pi^{-1}((\cx_{\cu}^{0})^{\perp})$.

The required shift closure follows from the aisle property of $\cx_{\cu}^{0}$. Shifting this torsion pair by $[-1]$ gives the torsion pair with left part $\cx_{\cu}^{1}$ used below.

\end{proof}

\begin{Def}\rm\label{Def:AGK-Higgs}
	We define the associated \emph{relative AGK category} as the idempotent completion of the following Verdier quotient
	$$\cc=\operatorname{Kar}\bigl(\ct/\ct^{\fd}\bigr).$$
    We define the associated \emph{relative fundamental domain} as the following full subcategory of $\ct$
	$$\cf=\cx^{\perp}\cap^{\perp}\!\cy\cap\cz,$$ 
	where $$\cx=\pi^{-1}\bigl(\overline{\cm}[<\!-1]^{\perp_{\cu}}\bigr)\cap\ct^{\fd}=\pi^{-1}(\cu^{\leqslant1})\cap\ct^{\fd}=\ct^{\leqslant1}\cap\ct^{\fd},$$ $$\cy=\cm[\geqslant0]^{\perp_{\ct}}\cap\thick(\cp)^{\perp_{\mathcal{T}}}\cap\,^{\perp_{\mathcal{T}}}\thick(\cp)=\ct^{\geqslant1}$$ and $$\cz=^{\perp_{\ct}}\!\cp[>0]\cap\cp[<0]^{\perp_{\ct}}.$$ We define the associated \emph{Higgs category} as the following full subcategory of $\cc$
	$$\ch=^{\perp_{\cc}}\!\cp[>0]\cap\cp[<0]^{\perp_{\cc}}.$$
\end{Def}

\begin{Rem}
	
Let $(\ct,\ct^{\fd},\cm,\cp)$ be a $(n+1)$-Calabi--Yau quadruple in the sense of \cite[Definition 3.1]{Wu2025}, with $n\geq1$. Then \Cref{Def:AGK-Higgs} agrees with the relative Higgs construction of \cite{Wu2025}, with the ambient quotient taken up to idempotent completion.
	If moreover $\cp=0$, we have 
	$\cz=\ct$, $\cf=\cx^{\perp}\cap^{\perp}\!\cy$ and $\ch=\cc$. And $(\ct,\ct^{\fd},\cm)$ is an $(n+1)$-Calabi--Yau triple in the sense of Iyama-Yang \cite[Section 5]{IY2018}. In this case, we have 
	$$\cf=(\cm[\leqslant-1])^{\perp}\cap^{\perp}\!(\cm[\geqslant n])=\cm*\cm[1]*\cdots\cm[n-1]\quad\text{and}\quad\ch=\cc.$$	
\end{Rem}

\bigskip

It is clear that $\cp$ is also a pre-silting subcategory of $\cc$. Denote by $\tau\colon\ct\ra\cc$ the quotient functor. It induces a functor $\cf\ra\ch$ which we still denote it by $\tau$. We have the following commutative diagram
\[ \begin{tikzcd} &\thick(\cp)\arrow[r,equal]\arrow[d,hook]&\thick(\cp)\arrow[d,hook]\\
\ct^{\fd}\arrow[r,hook]\arrow[d,equal]&\ct\arrow[r,two heads,"\tau"]\arrow[d,two heads,"\pi"]&\cc\arrow[d,two heads,"\pi"]\\
\cu^{\fd}\arrow[r,hook]&\cu\arrow[r,two heads,"\tau"]&\cu/\cu^{\fd}. \end{tikzcd} \]

\begin{Rem}\label{Rem:H/Pfullyfaithful}
	Equation~\Cref{eq:frozen-approximation-halves} gives condition (P1) of \cite[Section~3.1]{IY2018} for $\cp\subseteq\ct$. 
	The silting-reduction comparison gives a fully faithful embedding 
	\[ \pi\colon\mathcal{H}/[\mathcal{P}] \hookrightarrow \operatorname{Kar}\bigl(\mathcal{U}/\mathcal{U}^{\mathrm{fd}}\bigr). \]
	Under the hypotheses of \Cref{Thm:H is Frobenius and H/P=U/Ufd}, we will show that this embedding is an equivalence of triangulated categories.	
\end{Rem}

\begin{Lem}\label{Lem:t-structure on T^fd}
The pair $(\ct^{\fd}\cap\ct^{\leqslant0},\ct^{\geqslant0})$ is a bounded $t$-structure on $\ct^{\fd}$.	
\end{Lem}

\begin{proof}
	Let $X\in\ct^{\fd}$, there is a triangle
	$$\sigma^{\leqslant0}X\ra X\ra\sigma^{\geqslant1}X\ra\sigma^{\leqslant0}X[1].$$
	Since $\sigma^{\geqslant1}X\in\ct^{\geqslant0}\subseteq\ct^{\fd}$, then both $X$ and $\sigma^{\geqslant1}X$ belong to $\ct^{\fd}$. It follows that $\sigma^{\leqslant0}X$ lies in
	$\ct^{\fd}\cap\ct^{\leqslant0}$. This shows that $(\ct^{\fd}\cap\ct^{\leqslant0},\ct^{\geqslant0})$ is a $t$-structure on $\ct^{\fd}$.
	
	Let $X$ be any object of $\ct^{\fd}$. Then we have $\Hom_{\mathcal{T}}(\cm,X[i])=0$ and
\[ \Hom_{\mathcal{U}}(\overline{\cm},\pi(X)[-i])\simeq\Hom_{\ct}(\cm,X[-i])\simeq D\Hom_{\ct}(\mathbb{S}^{-1}(X),\cm[i])=0 \]
for $i\gg0$. Namely, $X$ belongs to $\ct^{\fd}\cap\ct^{\leqslant i}\cap\ct^{\geqslant -i}$ for some positive integers $i$.

The bounds are uniform on $\cm$, since both $X$ and $\mathbb S^{-1}X$ are built from finitely many shifts of objects of $\cm$.

\end{proof}

\begin{Lem}\label{Lem:morohism in C}\cite[Lemma 3.20]{Wu2025}
	Let $X$ be an object of $\ct$ such that $\pi(X)\in\cu_{\leqslant0}$ and $Y\in\ct$. Let $g$ be an element of $\Hom_{\cc}(X,Y)$. Then it has a representative of the form $X\xleftarrow{s} Z\xrightarrow{f} Y$ such that the cone of $s$ belongs to $\pi^{-1}(\cu_{\leqslant0})\cap\ct^{\fd}$.
\end{Lem}

\begin{Lem}\label{Lem:faithful functor of F}
	The projection functor $\tau\colon\ct\ra\cc$ induces a bijection
	$$\Hom_{\ct}(U,V)\ra\Hom_{\cc}(U,V)$$ for any $V\in\cx^{\perp}$ and $U\in\ct$ which satisfies $\pi(U)\in\cu_{\leqslant0}=\bigcup_{i\geqslant0}\overline{\mathcal M}*\overline{\mathcal M}[1]\cdots*\overline{\mathcal M}[i]$.
	
	Consequently, it restricts to a fully faithful functor $\tau\colon\cf\hookrightarrow\cc$.
\end{Lem}

\begin{proof}
The silting co-$t$-structure and (RS4) give
$\cu_{\leqslant0}=\overline{\cm}[<0]^{\perp_{\cu}}=\cu^{\leqslant0}$.
	Let $U$ be an object of $\ct$ such that $\pi(U)\in\cu_{\leqslant0}$ and $V$ an object of $\cx^{\perp}$. We use $\Hom_{\ct}(\cx,V)=0$ and $\cx=\ct^{\fd}\cap\ct^{\leqslant1}$.
	
	We first show that $\Hom_{\ct}(U,V)\ra \Hom_{\cc}(U,V)$ is injective. Assume that $f\in\Hom_{\ct}(U,V)$ becomes 
	zero in $\cc$. By \cite[Lemma 2.1.26]{Neeman2001}, $f$ factors through some
	$W\in\ct^{\fd}$. By {\let\CrefLemmaSpace\space\Cref{Lem:t-structure on T^fd}}, there exists
	a triangle
	$$\sigma^{\leqslant0}W\ra W\ra\sigma^{\geqslant1}W\ra\sigma^{\leqslant0}W[1]$$
	such that $\pi(\sigma^{\leqslant0}W)\in\cu^{\leqslant0}$ and
	$\sigma^{\geqslant1}W\in\cm[>\!0]^{\perp_{\mathcal{T}}}\cap\thick(\cp)^{\perp_{\mathcal{T}}}\cap\,^{\perp_{\mathcal{T}}}\thick(\cp)$. 
	
	Since $\Hom_{\ct}(U,\sigma^{\geqslant1}W)\cong\Hom_{\cu}(\pi(U),\pi(\sigma^{\geqslant1}W))$ vanishes, the map $f$ also factors through
	$\sigma^{\leqslant0}W$ along $\sigma^{\leqslant0}W\ra W$, i.e. we have the following
	commutative diagram in $\ct$
	\[ \begin{tikzcd}[ampersand replacement=\&] U\arrow[rr,"f"]\arrow[d]\arrow[rd]\&\&V\\
\sigma^{\leqslant0}W\arrow[r]\&W\arrow[r]\arrow[ur]\&\sigma^{\geqslant1}W\arrow[r]\&\sigma^{\leqslant0}W[1]. \end{tikzcd} \]
	
	Notice that we have
	$$\Hom_{\ct}(\sigma^{\leqslant0}W,V)=0$$
	as $\sigma^{\leqslant0}W\in\ct^{\fd}\cap\ct^{\leqslant0}\subseteq\cx$ and $V\in\cx^{\perp}$.
	Thus $f$ is zero.
	
	Next we show that $\Hom_{\ct}(U,V)\ra \Hom_{\cc}(U,V)$ is surjective.
	Let $g\colon U\ra V$ be a morphism in $\cc$. By {\let\CrefLemmaSpace\space\Cref{Lem:morohism in C}}, it has a representative of the form $U\xleftarrow{s}Z\xrightarrow{f}V$ such that the cone $W$ of $s$ belongs to $\pi^{-1}(\cu_{\leqslant0})\cap\ct^{\fd}$. Then we have an exact sequence
	$$\Hom_{\ct}(U,V)\xrightarrow{s^{*}}\Hom_{\ct}(Z,V)\ra\Hom_{\ct}(W[-1],V).$$

	Since $W\in\pi^{-1}(\cu_{\leqslant0})\cap\ct^{\fd}$, (RS5) gives
	$W[-1]\in\ct^{\fd}\cap\ct^{\leqslant1}=\cx$. As $V\in\cx^{\perp}$, we have $\Hom_{\ct}(W[-1],V)=0$.
	Therefore there exits $h\in\Hom_{\ct}(U,V)$ such that $f=hs$. Hence $U\xleftarrow{s}Z\xrightarrow{f}V$ is equivalent to $h\colon U\ra V$. This implies that $\Hom_{\ct}(U,V)\ra \Hom_{\cc}(U,V)$ is surjective. Hence $\tau\colon\ct\ra\cc$ induces a bijection
	$$\Hom_{\ct}(U,V)\ra\Hom_{\cc}(U,V)$$ for any $V\in\cx^{\perp}$ and $U\in\ct$ which satisfies $\pi(U)\in\cu_{\leqslant0}=\bigcup_{i\geqslant0}\overline{\mathcal M}*\overline{\mathcal M}[1]\cdots*\overline{\mathcal M}[i]$.

For $N\in\cf$, the condition $N\in{}^{\perp}\cy$ gives $\pi N\in{}^{\perp}\cu^{\geq1}=\cu^{\leq0}$: every object of $\cu^{\geq1}$ lifts to $\cy\subseteq\ct^{\fd}$ by (RS4)--(RS5), and the Hom comparison applies. Hence the preceding bijection gives a fully faithful functor $\tau\colon\cf\hookrightarrow\cc$. The inclusion of the Verdier quotient into its idempotent completion does not change these Hom spaces.

\end{proof}

\bigskip

Let $\cf_{\cu}$ be the shifted fundamental domain \cite[Definition 2.11]{IY2020} of the  relative Serre quadruple $(\cu,\cu^{\fd},\mathbb S,\overline{\cm})$ defined as
$$\cf_{\cu}=(\cx_{\cu}^{1})^{\perp}\cap^{\perp}\!\cy_{\cu},$$
where $\cx_{\cu}^{1}=\overline{\cm}[<\!-1]^{\perp_{\cu}}\cap\cu^{\fd}$ and $\cy_{\cu}=\overline{\cm}[\geqslant0]^{\perp_{\cu}}$. 
\begin{Lem}\label{Lem:F/P=F_U}
	We have the following equivalence of additive categories
	$$\pi\colon\frac{\cf}{[\cp]}\iso\cf_{\cu}.$$
\end{Lem}

\begin{proof}

For $Z\in\cz$, the Hom comparison and (RS4)--(RS5) give
\[ \Hom_{\ct}(\cx,Z)\cong\Hom_{\cu}(\cx_{\cu}^{1},\pi Z),\qquad \Hom_{\ct}(Z,\cy)\cong\Hom_{\cu}(\pi Z,\cy_{\cu}). \]
Thus $Z\in\cf$ if and only if $\pi Z\in\cf_{\cu}$.

It is easy to see that $\pi\colon\ct\ra\cu$ induces a functor $\pi\colon\frac{\cf}{[\cp]}\ra\cf_{\cu}$. By \Cref{Prop:siltingreduction}, we see that $\pi\colon\frac{\cf}{[\cp]}\ra\cf_{\cu}$ is fully faithful. It remains to show that it is dense. Let $X$ be an object of $\cf_{\cu}$. Since $(\overline{\cm}[<0]^{\perp_{\cu}},\cy_{\cu})$ is a torsion pair on $\cu$, then $$^{\perp}\cy_{\cu}=\overline{\cm}[<0]^{\perp_{\cu}}=\cu_{\leqslant0}=\bigcup_{k\geqslant0}\overline{\mathcal M}*\overline{\mathcal M}[1]\cdots*\overline{\mathcal M}[k].$$ Hence there exists $k\geqslant0$ such that $X$ belongs to $\overline{\cm}*\overline{\mathcal M}[1]\cdots*\overline{\mathcal M}[k]\subseteq\cu$, i.e. there exist $ k $ triangles in $\cu$
$$P_{1}\xrightarrow{f_{1}} Q_{0}\to X\to P_{1}[1],$$
$$ P_{2}\xrightarrow{f_{2}} Q_{1}\to P_{1}\to P_{2}[1],$$
$$ \cdots $$
$$P_{k}\xrightarrow{f_{k}} Q_{k-1}\to P_{k-1}\to P_{k}[1] ,$$
where $ Q_{0} $, $ Q_{1} $, $ \cdots $, $ Q_{k-1} $ and $ P_{k} $ are in $ \pi(\cm) $.	

If $k=0$, lift $X\in\overline{\cm}$ to an object of $\cm$; it lies in $\cf$ by the preceding Hom comparison. We may therefore assume $k\geq1$.

Since we have an equivalence of additive categories
$$\pi\colon\cm/[\cp]\ra\pi(\cm),$$ there exits a morphism $g'_{k}\colon M_{k}\ra N'_{k-1}$ such that $\pi(M_{k})= P_{k}$, $\pi(N'_{k-1})= Q_{k-1}$ and $\pi(g'_{k})=f_{k}$. Let $\gamma_{k}\colon M_{k}\ra R_{k-1}$ be a left $\cp$-approximation of $M_{k}$. Define
$$(M_{k}\xrightarrow{g_{k}}N_{k-1})\coloneq(M_{k}\xrightarrow{[g'_{k}\,\,\gamma_{k}]^{t}}N'_{k-1}\oplus R_{k-1}).$$
Then we still have $\pi(g_{k})=f_{k}$. Let $M_{k-1}=\cone(g_{k})$. It is easy to see that $M_{k-1} $ lies in $\cz=^{\perp_{\ct}}\!\cp[>0]\cap\cp[<0]^{\perp_{\ct}}$ and $\pi(M_{k-1})\cong P_{k-1}$. Indeed, the approximation gives $\Hom(M_{k-1},\cp[1])=0$, and the other defining vanishings of $\cz$ follow from the long exact Hom sequences. Hence we obtain a triangle in $\ct$ $$M_{k}\xrightarrow{g_{k}} N_{k-1}\ra M_{k-1}\ra M_{k}[1]$$
which lifts $P_{k}\xrightarrow{f_{k}} Q_{k-1}\to P_{k-1}\to P_{k}[1]$.

By \Cref{Prop:siltingreduction}, there exists a morphism $g'_{k-1}\colon M_{k-1}\ra N'_{k-2}$ such that $\pi(N'_{k-2})=Q_{k-2}$ and $\pi(g'_{k-1})=f_{k-1}$. Since $M_{k-1}\in\cz$, its left $\cp$-approximation exists by \Cref{eq:frozen-approximation-halves}. By the same augment above, we find a morphism $g_{k-1}\colon M_{k-1}\ra N_{k-2}$ such that $\pi(g_{k-1})=f_{k-1}$, $\pi(N_{k-2})=Q_{k-2}$ and $M_{k-2}\coloneqq\cone(g_{k-1})$ belongs to $\cz$. 

By induction, we have the following triangles
in $\ct$
$$M_{1}\xrightarrow{g_{1}} N_{0}\to Y\to M_{1}[1],$$
$$ M_{2}\xrightarrow{g_{2}} N_{1}\to M_{1}\to M_{2}[1],$$
$$ \cdots $$
$$M_{k}\xrightarrow{g_k} N_{k-1}\to M_{k-1}\to M_{k}[1] ,$$
where $ N_{0} $, $ N_{1} $, $ \cdots $, $ N_{k-1} $ and $ M_{k} $ are in $ \cm $ and $Y, M_{1},\cdots,M_{k-1}$ are in $\cz$. Moreover we have $\pi(Y)\cong X$ and $Y\in\cf$ by the Hom comparison at the start of the proof. This shows that $\pi\colon\frac{\cf}{[\cp]}\ra\cf_{\cu}$ is dense. Then the functor $\pi\colon\frac{\cf}{[\cp]}\ra\cf_{\cu}$ is an equivalence. 

At every stage the source lies in $\cz$, so \Cref{eq:frozen-approximation-halves} supplies the next approximation. Since all added frozen terms lie in $\cm$, the lift uses only $\cm,\cm[1],\ldots,\cm[k]$, as required for the bound on the cohomological shifts in \Cref{Prop:rf^Z=add} .

\end{proof}

\begin{Rem}\label{Rem:equivalence of ex-quotient}
	It is clear that $\cf$ is an extension closed subcategory of
	$\ct$. Hence it carries a canonical extriangulated structure $(\cf,\mathbb{E}_{\cf},\mathfrak{s})$ in the sense of Nakaoka-Palu \cite{NP2019} and it is described as follows$\colon$
	\begin{itemize}
		\item[(1)]  For any two objects $X,Z\in\cf$, the $\mathbb{E}$-extension space $\mathbb{E}_{\cf}(Z,X)$ is given by $\Hom_{\ct}(Z,X[1])$.
		\item[(2)] For any $\delta\in\mathbb{E}_{\cf}(Z,X)=\Hom_{\ct}(Z,X[1])$, take a distinguished triangle
		$$X\xrightarrow{f}Y\xrightarrow{g}Z\xrightarrow{\delta}X[1]$$ and define $\mathfrak{s}(\delta)=[X\xrightarrow{f}Y\xrightarrow{g}Z]$. This $\mathfrak{s}(\delta)$ does not depend on the choice of the distinguished
		triangle above. 
	\end{itemize}
Similarly we have a canonical extriangulated structure $(\cf_{\cu},\mathbb{E}_{\cf_{\cu}},\mathfrak{s})$ on $\cf_{\cu}$. It is not hard to see that $\pi\colon\frac{\cf}{[\cp]}\iso\cf_{\cu}$ is an equivalence of extriangulated categories.	

Indeed, $\cf\subseteq\cz$, and \cite[Lemma~3.4]{IY2018} identifies
\[ \mathbb E_{\cf}(Z,X)\cong\Hom_{\cu}(\pi Z,\pi X[1])=\mathbb E_{\cf_{\cu}}(\pi Z,\pi X). \]
Localization sends realizing triangles to realizing triangles, so this isomorphism is compatible with the two extriangulated structures.

\end{Rem}

\begin{Thm}\label{Thm:F=H}
	Assume that one of the following conditions holds.
	\begin{itemize}
		\item[(a)] $\mathbb{S}\colon\ct^{\fd}\ra\ct^{\fd}$ extends to a triangle equivalence $\mathbb{S}\colon\ct\ra\ct$ such that $\mathbb{S}^{\pm}(\cp)\subseteq\thick(\cp)$;
		\item[(b)] $\cm$ has an additive generator $M$.
	\end{itemize} Then the functor $\tau\colon\ct\ra\cc$ induces an equivalence of additive categories.
	$$\tau\colon\cf\iso\ch.$$
\end{Thm}	

\begin{proof}

Put $\cc_0=\ct/\ct^{\fd}$ and let $\ch_0\subseteq\cc_0$ be defined by the same frozen orthogonalities as $\ch$. We first prove the equivalence with $\ch_0$.

	By \Cref{Lem:faithful functor of F} and the full faithfulness of the canonical embedding $\cc_0\hookrightarrow\cc$, the natural quotient functor $\ct\to\cc_0$ restricts to a fully faithful functor $\cf\hookrightarrow\cc_0$. Since $\cp$ is left and right orthogonal to $\ct^{\fd}$, the image $\tau(\cf)$ lies in $\ch_{0}$. Hence we obtain a fully faithful functor $\cf\hookrightarrow\ch_0$. It is enough to show that it is dense. 

	Let $Z\in\ct$ represent an object of $\ch_0$. Since $\thick(\cp)$ and $\ct^{\fd}$ are orthogonal on both sides, the quotient functor $\ct\to\cc_0$ preserves Hom spaces whenever one argument belongs to $\thick(\cp)$. Thus the defining frozen orthogonalities of $\ch_0$ imply $Z\in\cz={}^\perp_{\ct}(\cp[>0])\cap(\cp[<0])^{\perp_{\ct}}$.

	Since $(\ct^{\leqslant0},\ct^{\geqslant0})$ is also a $t$-structure on $\ct$, there exits a triangle in $\ct$
	$$\sigma^{\leqslant0}Z\ra Z\ra\sigma^{\geqslant1}Z\ra\sigma^{\leqslant0}Z[1]$$
	with $\sigma^{\leqslant0}Z\in\ct^{\leqslant0}=\pi^{-1}(\cu^{\leqslant0})$ and $\sigma^{\geqslant1}Z\in\ct^{\geqslant1}=\cm[\geqslant\!0]^{\perp_{\mathcal{T}}}\cap\thick(\cp)^{\perp_{\mathcal{T}}}\cap\,^{\perp_{\mathcal{T}}}\thick(\cp)\subseteq\ct^{\fd}$. Since $\sigma^{\geqslant1}Z\in\ct^{\fd}$, we have $Z\simeq\sigma^{\leqslant0}Z$ in $\cc_0$. Moreover, applying Hom to the truncation triangle and using the two-sided orthogonality between $\ct^{\fd}$ and $\cp$, we obtain $\sigma^{\leqslant0}Z\in\cz$. 

	By the $[-1]$-shift of the torsion pair in \Cref{Lem:left t-structure}, $\ct$ has a torsion pair with left part $\cx_{\cu}^{1}\subseteq\ct^{\fd}$ and right part $\pi^{-1}((\cx_{\cu}^{1})^{\perp})$. Then we have a triangle in $\ct$
	$$X\ra\sigma^{\leqslant0}Z\ra Z'\ra X[1] $$
	with $X\in\cx_{\cu}^{1}\subseteq\ct^{\fd}\subseteq\cz$ and $Z'\in\pi^{-1}((\cx_{\cu}^{1})^{\perp})$. Since $X\in\ct^{\fd}$, we have $Z\simeq\sigma^{\leqslant0}Z\simeq Z'$ in $\cc_0$. Moreover, the two-sided orthogonality between $\ct^{\fd}$ and $\cp$ gives $Z'\in\cz$. By the Hom comparison, the condition $\pi Z'\in(\cx_{\cu}^{1})^{\perp}$ implies $Z'\in\cx^{\perp}$. Finally, $X[1]$ and $\sigma^{\leqslant0}Z$ lie in $^{\perp}\cy=\ct^{\leqslant0}$, which is extension-closed. Hence $Z'\in{}^{\perp}\cy$, and therefore $Z'\in\cf$. Thus $\cf\simeq\ch_0$.

The category $\ch_0$ is idempotent complete and Krull--Schmidt because $\cf$ is summand-closed in $\ct$. It remains to treat the completion $\cc=\operatorname{Kar}(\cc_0)$. The two Verdier localizations commute after completion, giving
\[ \operatorname{Kar}(\cc/\thick(\cp))\simeq \operatorname{Kar}(\cu/\cu^{\fd}). \]
The right-hand category is already idempotent complete by \cite[Theorem~2.12]{IY2020}. Since the quotient functor preserves Hom spaces involving $\cp$, the image of $\cp$ in $\cc$ remains presilting. As $\ch$ satisfies the corresponding frozen orthogonalities, the fully faithful part of \cite[Lemma~3.4]{IY2018} gives a fully faithful functor $\ch/[\cp]\to\cc/\thick(\cp)$. Composing with the canonical embedding into its idempotent completion and the preceding equivalence of quotients gives a fully faithful functor $\ch/[\cp]\to\cu/\cu^{\fd}$. Its restriction to $\ch_0/[\cp]$ is essentially surjective: combine $\ch_0\simeq\cf$, \Cref{Lem:F/P=F_U}, and the shifted fundamental-domain equivalence of \cite[Theorem~2.12]{IY2020}.

Given $X\in\ch$, choose $Y\in\ch_0$ with the same image in this quotient. Full faithfulness provides maps $f\colon X\to Y$ and $g\colon Y\to X$ for which $gf-1_X=vu$ factors through some $P\in\cp$, with $u\colon X\to P$ and $v\colon P\to X$. Then
\[ X\xrightarrow{(f,u)}Y\oplus P\xrightarrow{(g,-v)}X \]
has composite $1_X$. The category $\ch_0$ is full, idempotent complete, and contains $Y\oplus P$, so it also contains an object isomorphic to the retract $X$. Hence $\ch=\ch_0$, proving the result.

\end{proof}

The main result of this Appendix is the following theorem, which extends \cite[Theorem 2.12]{IY2020} to the relative setting.

\begin{Thm}\label{Thm:H is Frobenius and H/P=U/Ufd}
Let $(\ct,\ct^{\fd},\mathbb{S},\cm,\cp)$ be a relative Serre five-tuple.
Assume that one of the following conditions holds.
\begin{itemize}
	\item[(a)] $\mathbb{S}\colon\ct^{\fd}\ra\ct^{\fd}$ extends to a triangle equivalence $\mathbb{S}\colon\ct\ra\ct$ such that $\mathbb{S}^{\pm}(\cp)\subseteq\thick(\cp)$;
	\item[(b)] $\cm$ has an additive generator $M$.
\end{itemize}
Then we have
	\begin{itemize}
		\item[(1)] The Higgs category $\ch$ is a Krull--Schmidt Frobenius extriangulated category whose projective-injective objects are $\cp$. Furthermore, its stable category $\frac{\ch}{[\cp]}$ is triangulated equivalent to $\frac{\cu}{\cu^{\fd}}$ and under (a) it admits a Serre functor given by $\mathbb{S}\circ[-1]$.
		\item[(2)] Under (a), let $n\geqslant1$, and assume that $\cm$ is stable under $\mathbb{S}_{n+1}=\mathbb{S}\circ[-n-1]$. Then $\cm$ is an $n$-cluster tilting subcategory of $\ch$, and its image is an $n$-cluster tilting subcategory of $\cu/\cu^{\fd}$.
	\end{itemize}

The stable category is Hom-finite. If $\ct$ is Hom-finite, then $\ch$ is Hom-finite as well.

\end{Thm}
\begin{proof}
\emph{Step 1: $\ch$ is an extriangulated category.} 

The defining orthogonalities of $\mathcal H$ are preserved under extensions. Hence $\mathcal H$ is an extension-closed subcategory of $\cc$. Then $\mathcal{H}$ is an extriangulated category and its extriangulated structure $(\mathcal{H}, \mathbb{E}_{\ch}, \mathfrak{s})$ is described as follows:
\begin{enumerate}
	\item[(1)] For any two objects $X, Z \in \mathcal{H} \subseteq \mathcal{C}$, the $\mathbb{E}_{\ch}$-extension space $\mathbb{E}_{\ch}(Z, X)$ is given by
	\[ \text{Hom}_{\mathcal{C}}(Z,X[1]). \]
	\item[(2)] For any $\delta \in \mathbb{E}_{\ch}(Z, X) = \text{Hom}_{\mathcal{C}}(Z, X[1])$, take a distinguished triangle
	\[ X \xrightarrow{f} Y \xrightarrow{g} Z \xrightarrow{\delta} X[1] \]
	and define $\mathfrak{s}(\delta) = [X \xrightarrow{f} Y \xrightarrow{g} Z]$. This $\mathfrak{s}(\delta)$ does not depend on the choice of the distinguished triangle above.
\end{enumerate}

\emph{Step 2: $\ch$ has  enough injectives and the full subcategory of injective objects in $\ch$ is $\cp$. And $\ch$ has  enough projectives and the full subcategory of projective objects in $\ch$ is $\cp$.}

The proof follows the lines of that of \cite[Theorem 3.25]{Wu2025}.

For $X\in\ch$, \Cref{Thm:F=H} gives $X\simeq\tau Z$ with $Z\in\cf\subseteq\cz$. By \Cref{eq:frozen-approximation-halves} and the Hom comparison, $X$ admits both left and right $\cp$-approximations. Completing them to triangles gives
\[ K\longrightarrow P_0\longrightarrow X\longrightarrow K[1], \qquad X\longrightarrow P^0\longrightarrow L\longrightarrow X[1], \]
with $P_0,P^0\in\cp$. As in the cited proof, $K,L\in\ch$; the boundary-degree vanishings follow from the surjections
\[ \Hom_{\cc}(P,P_0)\twoheadrightarrow\Hom_{\cc}(P,X), \qquad \Hom_{\cc}(P^0,P)\twoheadrightarrow\Hom_{\cc}(X,P) \quad(P\in\cp). \]
Since
\[ \mathbb E_{\ch}(\cp,\ch)=\mathbb E_{\ch}(\ch,\cp)=0, \]
these conflations give enough projectives and injectives. Splitting them for a projective or injective object yields
\[ \operatorname{proj}\ch=\operatorname{inj}\ch=\cp. \]

\emph{Step 3: the stable category $\frac{\ch}{[\cp]}$ is triangulated equivalent to $\frac{\cu}{\cu^{\fd}}$ and under (a) it admits a Serre functor given by $\mathbb{S}\circ[-1]$.} 

By \Cref{Thm:F=H} and {\let\CrefLemmaSpace\space\Cref{Lem:F/P=F_U}}, we have equivalences of additive categories $\cf_{\cu}\simeq\frac{\cf}{[\cp]}\simeq\frac{\ch}{[\cp]}$. Since  the  relative Serre quadruple $(\cu,\cu^{\fd},\mathbb{S},\overline{\cm})$ satisfies the Assumption in \cite[Theorem 2.12]{IY2020}, the composition $\cf_{\cu}\subseteq\cu\ra\frac{\cu}{\cu^{\fd}}$ is an equivalence of additive categories. Hence we have $\frac{\ch}{[\cp]}\simeq\frac{\cu}{\cu^{\fd}}$. This equivalence is triangulated: an injective conflation $X\to P^0\to Y$ is sent to a distinguished triangle $X\to0\to Y\to X[1]$, so $Y\simeq X[1]$ in $\cu/\cu^{\fd}$. Thus the stable suspension agrees with the triangulated shift, and conflations are sent to distinguished triangles. By \cite[Theorem 2.12]{IY2020}, under (a), $\frac{\cu}{\cu^{\fd}}$ admits a Serre functor given by $\mathbb{S}\circ[-1]$.

The same theorem makes $\cu/\cu^{\fd}$ Hom-finite and Krull--Schmidt. Also $\ch\simeq\cf$ is Krull--Schmidt, and is Hom-finite when $\ct$ is.

\emph{Step 4: $\cm$ is an $n$-cluster tilting subcategory of $\ch$.} Since $\cm$ is stable under $\mathbb{S}_{n+1}$, by \cite[Theorem 2.6]{Hanihara2022}, we see that the image of $\overline{\cm}$ is an
$n$-cluster tilting subcategory of $\frac{\cu}{\cu^{\fd}}$. Then it is easy to see that $\cm$ is an
$n$-cluster tilting subcategory of $\ch$.

To see this directly, first let $M\in\cm$ and $X\in\cx$. Silting gives $M\in\cz$ and $\Hom(M,\cy)=0$. Relative duality and $\mathbb S^{-1}(\cm)=\cm[-n-1]$ give
\[ D\Hom_{\ct}(X,M)\cong\Hom_{\ct}(\mathbb S^{-1}M,X)=0, \]
since $\pi X\in\cu^{\leq1}$ and $n+1\geq2$. Hence $\cm\subseteq\cf$, and its image lies in $\ch$.

The $\cm$-approximations in the stable category lift to approximations modulo $[\cp]$. Adjoining the frozen approximations from \Cref{eq:frozen-approximation-halves} gives actual left and right $\cm$-approximations in $\ch$. Moreover, dimension shifting gives
\[ \mathbb E_{\ch}^{i}(X,Y) \cong\Hom_{\ch/[\cp]}(X,Y[i]),\qquad i\ge1. \]
Hence the $n$-cluster-tilting orthogonality conditions lift from the stable category. Finally, if $X$ is stably isomorphic to some $M\in\cm$, the retract argument in \Cref{Thm:F=H} makes $X$ a direct summand of $M\oplus P\in\cm$, so $X\in\cm$.

\end{proof}

\section{Exact dg orbit categories}\label[appendix]{appB:exact-dg-orbit}

Throughout this appendix, let $k$ be a field. We use the term
\emph{exact dg category} in the sense of Chen
\cite{Chen20261,Chen20262}. Let $(\mathcal E,\mathscr S)$ be a small
connective exact dg category over $k$. For a dg category $\mathcal A$, we denote its homotopy category by
$\mathrm H^0(\mathcal A)$. We put
$\mathscr E:=\mathrm H^0(\mathcal E)$
and write
\[ \bigl(\mathscr E,\mathbb E_{\mathcal E}, \mathfrak s_{\mathcal E}\bigr) \]
for the extriangulated category associated with
$(\mathcal E,\mathscr S)$. The dependence of
$\mathbb E_{\mathcal E}$ and $\mathfrak s_{\mathcal E}$ on
$\mathscr S$ is suppressed from the notation.

We use $\operatorname{pretr}(-)$ for the pretriangulated dg hull and
$(-)^{\sharp}$ for a fixed dg idempotent completion. 

\begin{Lem}
\label{lem:induced-exact-enhancement}
Let $\mathcal V$ be a pretriangulated dg category and
$\mathscr X\subseteq H^0(\mathcal V)$ a full extension-closed additive
subcategory. Let $\widetilde{\mathcal X}\subseteq\mathcal V$ be the full
dg subcategory on its objects, and put
$\mathcal X=\tau_{\leq0}\widetilde{\mathcal X}$.
\begin{enumerate}
\item[(1)] The dg category $\mathcal X$ has a natural exact structure whose
associated extriangulated category is $\mathscr X$, with
\[ \mathbb E_{\mathcal X}(Z,X)=\Hom_{H^0(\mathcal V)}(Z,X[1]). \]
Its conflations are characterized by the cone comparison
\Cref{eq:induced-exact-cone-comparison} below.

\item[(2)] If $G:\mathcal V\to\mathcal V'$ is a dg functor with
$H^0(G)(\mathscr X)\subseteq\mathscr X'$, it induces an exact dg functor
$\mathcal X\to\mathcal X'$, defined from $\mathscr X'\subseteq H^0(\mathcal V')$
in the same way.

\item[(3)] Suppose, in addition, that $\mathscr X$ has enough projectives and
\begin{equation}\label{eq:induced-realization-vanishing} \Hom_{H^0(\mathcal V)}(P,X[r])=0 \qquad(P\in\operatorname{proj}\mathscr X,\ X\in\mathscr X,\ r\geq2). \end{equation}
Then the inclusion induces a quasi-fully faithful realization
\[ \rho_{\mathscr X}:\cd^b_{\mathrm{dg}}(\mathcal X)\longrightarrow\mathcal V \]
with homotopy image $\operatorname{tria}_{H^0(\mathcal V)}(\mathscr X)$.
If $H^0(\mathcal V)$ is idempotent complete and
$\thick(\mathscr X)=H^0(\mathcal V)$, the realization becomes a
quasi-equivalence after dg idempotent completion of its source.
\end{enumerate}
\end{Lem}

\begin{proof}
\noindent\emph{(1).}
The induced exact structure follows from the canonical exact structure
on $\mathcal V$, extension-closedness, and connective truncation
\cite{Chen20261}. Explicitly, a conflation is a three-term homotopy
complex given by closed degree-zero maps $f:X\to Y$, $g:Y\to Z$ and a
specified degree-$-1$ map $h:X\to Z$ with $d(h)=-gf$.
Writing $X^\wedge=\mathcal V(-,X)$, its cone comparison is
\begin{equation}\label{eq:induced-exact-cone-comparison} \begin{gathered} \operatorname{Cone}(f^\wedge)=Y^\wedge\oplus X^\wedge[1],\qquad d_{\operatorname{Cone}(f^\wedge)}= \begin{pmatrix}d_{Y^\wedge}&f^\wedge\\
0&-d_{X^\wedge}\end{pmatrix},\\
\chi_{g,h}=(g^\wedge,-h^\wedge): \operatorname{Cone}(f^\wedge)\longrightarrow Z^\wedge. \end{gathered} \end{equation}
It is homotopy short exact precisely when $\chi_{g,h}$ is a
quasi-isomorphism of dg $\mathcal V$-modules; in particular, the
specified $h$ matters.

\smallskip\noindent\emph{(2).}
Representable dg modules and their finite cones are homotopically
projective, so this quasi-isomorphism is a homotopy equivalence.
A dg functor preserves the specified homotopy and the corresponding
homotopy equivalence of cones. The restriction to connective covers is
therefore exact.

\smallskip\noindent\emph{(3).}
For $f:X\to Y[r]$ with $r\geq2$, choose a projective conflation with
ambient triangle
\begin{equation}\label{eq:frobenius-syzygy-triangle} \Omega X\longrightarrow P_X\xrightarrow{p_X}X \xrightarrow{\delta_X}\Omega X[1]. \end{equation}
The vanishing \Cref{eq:induced-realization-vanishing} gives
$fp_X=0$, so $f$ factors through $\delta_X$.

The conclusion follows directly from
\cite[Corollary~3.4]{Plogmann26}. Finally,
$\operatorname{Kar}(\operatorname{tria}(\mathscr X)) =\thick(\mathscr X)$, so completion of the fully faithful realization
is essentially surjective under the last hypothesis.
\end{proof}

The realization is induced by the dg inclusion and is compatible with
automorphisms preserving $\mathscr X$, using functorial twisted-complex
and dg-quotient models. The same induced exact structure is used for
the Higgs enhancement in \Cref*{sec:weighted-totalization-singularity} and the orbit hull below.

\subsection{The ordinary and derived dg orbit categories}

Let
\[ F\colon(\mathcal E,\mathscr S) \longrightarrow(\mathcal E,\mathscr S) \]
be an exact dg automorphism.

We first fix a strictly connective model. Here \emph{connective}
means $H^j\mathcal E(X,Y)=0$ for $j>0$, whereas \emph{strictly
connective} means $\mathcal E(X,Y)^j=0$ for $j>0$.
Let $\mathscr S_{\leq0}$ be the exact structure on
$\tau_{\leq0}\mathcal E$ induced by $\mathscr S$. The canonical
inclusion
\[ (\tau_{\leq0}\mathcal E,\mathscr S_{\leq0}) \longrightarrow(\mathcal E,\mathscr S) \]
is an exact quasi-equivalence, and it is equivariant for
$\tau_{\leq0}F$ and $F$. The functor $\tau_{\leq0}F$ is an exact
dg automorphism, with inverse $\tau_{\leq0}(F^{-1})$.
We replace $(\mathcal E,\mathscr S,F)$ by this strictly connective
model and retain the notation. Exact structures transport along this quasi-equivalence \cite{Chen20261},
and the bounded dg derived categories are quasi-equivalent by their
universal property \cite{Chen20262}. The equivariant inclusion also
induces quasi-isomorphisms on the direct-sum Hom complexes of the
ordinary dg orbit categories. 

Hence the choice of model affects neither
these constructions up to quasi-equivalence nor the associated
extriangulated categories up to equivalence.

Let
$\iota\colon\mathcal E\longrightarrow \mathcal \cd:=\cd^b_{\mathrm{dg}}(\mathcal E,\mathscr S)$
be Chen's universal exact morphism
\cite[Theorem~3.1]{Chen20262}. Recall that Chen's bounded dg derived
category is constructed as
\[ \cd^b_{\mathrm{dg}}(\mathcal E,\mathscr S) = \operatorname{pretr}(\mathcal E)/\mathcal N_{\mathrm{dg}}, \]
where $\mathcal N_{\mathrm{dg}}$ is the full pretriangulated dg
subcategory generated by the totalizations of the conflations in
$\mathscr S$.

Since $F$ is exact, it sends the totalization of a conflation to the
totalization of a conflation. Consequently, it preserves
$\mathcal N_{\mathrm{dg}}$ and induces a dg automorphism
\[ \widehat F\colon\mathcal \cd\longrightarrow\mathcal \cd \]
together with a dg natural isomorphism
\[ \widehat F\iota\overset{\sim}{\Longrightarrow}\iota F. \]

We use the twisted-complex model for $\operatorname{pretr}(\mathcal E)$
and the Drinfeld dg-quotient model for
$\cd^b_{\mathrm{dg}}(\mathcal E,\mathscr S)$. Since $F$ preserves
conflations, it acts on both constructions and induces a strict dg
automorphism $\widehat F$, with inverse induced by $F^{-1}$.
Comparisons with other dg models are understood up to quasi-equivalence.

\begin{Def}\cite{KellerOrbit}
	\label{def:ordinary-dg-orbit}
	The \emph{ordinary dg orbit category} $\mathcal E/F^{\mathbb Z}$
	has the same objects as $\mathcal E$, and its morphism complexes are
	\[ (\mathcal E/F^{\mathbb Z})(X,Y) := \bigoplus_{n\in\mathbb Z}\mathcal E(X,F^nY). \]
	Thus every element has finite support in the orbit direction. For
	homogeneous elements
	\[ f=(f_n)_n\in(\mathcal E/F^{\mathbb Z})(X,Y), \qquad g=(g_m)_m\in(\mathcal E/F^{\mathbb Z})(Y,Z), \]
	composition is defined by
	\[ (g\circ f)_r := \sum_{n+m=r}F^n(g_m)\circ f_n, \]
	and the differential is defined componentwise by
	\[ d(f)_n:=d_{\mathcal E}(f_n). \]
	The identity of $X$ has component $1_X$ in orbit degree zero and
	zero components in all other orbit degrees.
\end{Def}

For a strict dg automorphism, the ordinary dg orbit category above is
quasi-equivalent to the derived dg orbit category of
\cite[Definition~3.11]{FKQ2024}. Indeed, the canonical maps of the
left-lax quotient become strictly invertible in the ordinary orbit,
so dg localization induces a comparison quasi-functor to it.
The Hom-complex calculation in \cite[Remark~3.12]{FKQ2024} shows that
this comparison is quasi-fully faithful, and it is essentially
surjective since both orbit categories have the same objects.
We use this comparison also for $\widehat F$ and the induced
automorphisms below.

Let

\begin{equation}\label{eq:derived-orbit-hull} \mathcal T_F:=\left(\operatorname{pretr}(\mathcal D/\widehat F^{\mathbb Z})\right)^{\sharp}. \end{equation}
be the idempotent-completed
pretriangulated hull of $\mathcal D/\widehat F^{\mathbb Z}$. We put
$\mathscr T_F:=\mathrm H^0(\mathcal T_F).$
Let
$Q\colon\mathcal D\longrightarrow\mathcal T_F$ and
$J=Q\iota\colon\mathcal E\longrightarrow\mathcal T_F$
be the canonical dg functors. The canonical orbit isomorphism gives a
dg natural isomorphism
\[ JF\overset{\sim}{\Longrightarrow}J. \]

\subsection{The exact dg orbit hull}

For full subcategories $\mathscr X,\mathscr Y$ of a triangulated
category $\mathscr T$, let $\mathscr X*\mathscr Y$ be the strictly
full subcategory formed by the objects $M$ which occur in a
distinguished triangle
\[ X\longrightarrow M\longrightarrow Y\longrightarrow X[1] \]
with $X\in\mathscr X$ and $Y\in\mathscr Y$. Inside $\mathscr T_F$,
define recursively
\[ \begin{aligned} \mathscr X_1 &:={\operatorname{add}} J(\mathscr E),\\
\mathscr X_{r+1} &:={\operatorname{add}} (\mathscr X_r*\mathscr X_1), \qquad r\geq1, \end{aligned} \]
and put
\[ \mathscr O_F:=\bigcup_{r\geq1}\mathscr X_r. \]
Here ${\operatorname{add}}$ denotes closure under finite
direct sums, isomorphisms, and direct summands in the idempotent-complete
triangulated category $\mathscr T_F$.

\begin{Lem}
	\label{lem:orbit-extension-closure}
	The sequence $(\mathscr X_r)_{r\geq1}$ is increasing. Moreover, for
	all $r,s\geq1$,
	\[ {\operatorname{add}} (\mathscr X_r*\mathscr X_s) \subseteq \mathscr X_{r+s}. \]
	Consequently, $\mathscr O_F$ is the smallest strictly full,
	extension-closed, and idempotent-complete additive subcategory of
	$\mathscr T_F$ containing $J(\mathscr E)$.
\end{Lem}

\begin{proof}
	Since $0\in\mathscr X_1$, the split triangle
	\[ X\xrightarrow{1_X}X\longrightarrow0\longrightarrow X[1] \]
	shows that $\mathscr X_r\subseteq\mathscr X_{r+1}$. We prove the second assertion by induction on $s$. The case $s=1$
	is the definition of $\mathscr X_{r+1}$. Suppose the assertion is
	known for $s$. The octahedral axiom gives the associativity of the
	extension operation,
	\[ (\mathscr X_r*\mathscr X_s)*\mathscr X_1 = \mathscr X_r*(\mathscr X_s*\mathscr X_1), \]
	up to passage to strictly full subcategories. We can enlarge the end terms by finite direct sums using split
	triangles, and handle their direct summands by adding the complementary
	split triangles. Hence passing to $\operatorname{add}$ is compatible
	with the inclusion given by the octahedral axiom. We
	obtain
	\[ \begin{aligned} {\operatorname{add}} (\mathscr X_r*\mathscr X_{s+1}) &\subseteq {\operatorname{add}} \bigl((\mathscr X_r*\mathscr X_s)*\mathscr X_1\bigr)\\
&\subseteq {\operatorname{add}} (\mathscr X_{r+s}*\mathscr X_1)\\
&=\mathscr X_{r+s+1}. \end{aligned} \]
	This proves the induction step.
	
	It follows that the union $\mathscr O_F$ is closed under extensions,
	finite direct sums, and direct summands. Conversely, let
	$\mathscr U\subseteq\mathscr T_F$ be a strictly full,
	extension-closed additive subcategory which contains $J(\mathscr E)$
	and is closed under direct summands. Then $\mathscr U$ contains
	$\mathscr X_1$, and induction shows that it contains every
	$\mathscr X_r$. This proves the minimality assertion.
\end{proof}

Let $\widetilde{\mathcal O}_F$ be the full dg subcategory of
$\mathcal T_F$ whose objects are the objects of $\mathscr O_F$.

\begin{Def}\rm
	\label{def:exact-dg-orbit-hull}
	The \emph{exact dg orbit hull} of $(\mathcal E,\mathscr S)$ by
	$F$ is the connective dg category
	\[ \operatorname{Orb}^{\mathrm{ex}}_{\mathrm{dg}}(\mathcal E,F) := \tau_{\leq0}\widetilde{\mathcal O}_F. \]
	Its conflations are the three-term homotopy complexes whose images
	under the canonical dg functor
	\[ \jmath_F\colon \operatorname{Orb}^{\mathrm{ex}}_{\mathrm{dg}}(\mathcal E,F) \longrightarrow\mathcal T_F \]
	are homotopy short exact.

{\color{black}
Explicitly, such a complex is given by closed degree-zero morphisms
$f\colon X\to Y$ and $g\colon Y\to Z$, together with a specified
$h\colon X\to Z$ of degree $-1$ satisfying $d(h)=-gf$.
After applying $\jmath_F$ and the dg Yoneda embedding, the condition
is that
\[ (g^\wedge,-h^\wedge)\colon \operatorname{Cone}(f^\wedge)\longrightarrow Z^\wedge \]
be a quasi-isomorphism of dg $\mathcal T_F$-modules, with
$\operatorname{Cone}(f^\wedge)=Y^\wedge\oplus X^\wedge[1]$
and differential
$\left(\begin{smallmatrix}d_{Y^\wedge}&f^\wedge\\0&-d_{X^\wedge}\end{smallmatrix}\right)$.
Thus the induced triangle is determined using the specified homotopy,
not merely the two morphisms in $\mathscr T_F$.
}

\end{Def}

\begin{Prop}
	\label{prop:exact-orbit-is-exact}
	The dg category
	$\operatorname{Orb}^{\mathrm{ex}}_{\mathrm{dg}}(\mathcal E,F)$,
	endowed with the conflations of
	\Cref{def:exact-dg-orbit-hull}, is a connective exact dg
	category. Its associated extriangulated category is
	\[ \mathrm H^0\!\left( \operatorname{Orb}^{\mathrm{ex}}_{\mathrm{dg}}(\mathcal E,F) \right) =\mathscr O_F. \]
	The functor $J=Q\iota\colon\mathcal E\longrightarrow\mathcal T_F$ factors through an exact dg functor
	\[ J\colon(\mathcal E,\mathscr S) \longrightarrow \operatorname{Orb}^{\mathrm{ex}}_{\mathrm{dg}}(\mathcal E,F) \]
	equipped with a dg natural isomorphism $JF\cong J$.
\end{Prop}

\begin{proof}
{\color{black}
	The pretriangulated dg category $\mathcal T_F$ carries its canonical
	exact dg structure, whose conflations are the homotopy short exact
	sequences. \Cref{lem:orbit-extension-closure} shows that
	$\widetilde{\mathcal O}_F$ is a full extension-closed dg subcategory
	of $\mathcal T_F$. It therefore inherits the exact dg structure
	induced from $\mathcal T_F$. The connective truncation
	$\tau_{\leq0}\widetilde{\mathcal O}_F$ inherits the same class of
	three-term homotopy complexes and hence the induced exact dg structure.
	This is precisely the exact structure in
	\Cref{def:exact-dg-orbit-hull}. Since connective truncation
	does not change zeroth cohomology, its associated extriangulated
	category is $\mathscr O_F$.
}

The universal morphism $\iota$ is exact, so it sends each conflation
$(f,g,h)$ to a homotopy short exact complex in $\mathcal D$.
Its cone comparison is a homotopy equivalence after passing to the
pretriangulated hull. The dg functor $Q$, extended to that hull,
preserves this homotopy equivalence and the specified homotopy.
Consequently $(Jf,Jg,Jh)$ is homotopy short exact in $\mathcal T_F$.
All three terms lie in $\mathscr X_1$, so it is a conflation of the
induced exact structure.

Since $\mathcal E$ is strictly connective, we have
$\mathcal E=\tau_{\leq0}\mathcal E$. By functoriality of connective
truncation, $J$ factors through
$\tau_{\leq0}\widetilde{\mathcal O}_F$, and this factorization is exact
by the preceding argument. Finally, the canonical dg natural
isomorphism $Q\widehat F\cong Q$ and its inverse have closed
degree-zero components, so they survive the truncation and yield
$JF\cong J$.

\end{proof}

\begin{Prop}[Morphisms and first extensions]
	\label{prop:orbit-hom-ext}
	For $X,Y\in\mathscr E$ and $p\in\mathbb Z$, there are natural
	isomorphisms
	\[ \begin{aligned} \operatorname{Hom}_{\mathscr T_F}(JX,JY[p]) &\cong \bigoplus_{n\in\mathbb Z} \operatorname{Hom}_{\cd^b(\mathcal E,\mathscr S)} \bigl(\iota X,\iota F^nY[p]\bigr), \end{aligned} \]
	where
	$\cd^b(\mathcal E,\mathscr S) :=\mathrm H^0\bigl(\cd^b_{\mathrm{dg}}(\mathcal E,\mathscr S)\bigr).$
	In particular,
	\[ \operatorname{Hom}_{\mathscr O_F}(JX,JY) \cong \bigoplus_{n\in\mathbb Z} \operatorname{Hom}_{\mathscr E}(X,F^nY), \]
	and
	\[ \mathbb E_{\mathscr O_F}(JZ,JX) \cong \bigoplus_{n\in\mathbb Z} \mathbb E_{\mathcal E}(Z,F^nX). \]
\end{Prop}

\begin{proof}
	The canonical dg functors from a dg category to its pretriangulated
	hull and from a pretriangulated dg category to its dg idempotent
	completion are quasi-fully faithful on the objects. Hence
	there is a natural quasi-isomorphism
	\[ \mathcal T_F(JX,JY) \simeq \bigoplus_{n\in\mathbb Z} \mathcal D(\iota X,\widehat F^n\iota Y). \]
	Taking $p$-th cohomology and using
	$\widehat F^n\iota Y\simeq\iota F^nY$ gives the first formula.
	
	For $p=0$, Chen's embedding theorem
	\cite[Theorem~3.1]{Chen20262} gives
	$\operatorname{Hom}_{\cd^b(\mathcal E,\mathscr S)} (\iota X,\iota Y) \cong \operatorname{Hom}_{\mathscr E}(X,Y).$ For $p=1$, the same theorem gives
	\[ \operatorname{Hom}_{\cd^b(\mathcal E,\mathscr S)} (\iota Z,\iota X[1]) \cong \mathbb E_{\mathcal E}(Z,X). \]
	Since the extriangulated structure on $\mathscr O_F$ is induced
	from $\mathscr T_F$, we have
	\[ \mathbb E_{\mathscr O_F}(JZ,JX) \cong \operatorname{Hom}_{\mathscr T_F}(JZ,JX[1]). \]
	Combining these identifications proves the two special cases.
\end{proof}

\subsection{The Frobenius case}

Assume from now on that $\mathscr E$ is Frobenius. Put
\[ \mathscr P:=\operatorname{proj}\mathscr E =\operatorname{inj}\mathscr E, \]
let $\mathcal P\subseteq\mathcal E$ be the full dg subcategory on the
objects of $\mathscr P$, and assume that $F(\mathcal P)=\mathcal P.$ Inside $\mathscr O_F$, put
$ \mathscr P_F:={\operatorname{add}} J(\mathscr P).$

\begin{Thm}
	\label{thm:frobenius-exact-orbit}
	The extriangulated category $\mathscr O_F$ is Frobenius, and
	\[ \operatorname{proj}\mathscr O_F = \operatorname{inj}\mathscr O_F = \mathscr P_F. \]
\end{Thm}

\begin{proof}
	Let $P\in\mathscr P$ and $X\in\mathscr E$. By
	\Cref{prop:orbit-hom-ext},
	\[ \begin{aligned} \mathbb E_{\mathscr O_F}(JP,JX) &\cong \bigoplus_{n\in\mathbb Z} \mathbb E_{\mathcal E}(P,F^nX)=0,\\
\mathbb E_{\mathscr O_F}(JX,JP) &\cong \bigoplus_{n\in\mathbb Z} \mathbb E_{\mathcal E}(X,F^nP)=0. \end{aligned} \]
	The first vanishing uses projectivity of $P$, while the second uses
	$F^nP\in\mathscr P$ and hence injectivity of $F^nP$.
	
	Fix $P\in\mathscr P$. If
	$A\longrightarrow B\longrightarrow C\longrightarrow A[1]$
	is a distinguished triangle in $\mathscr T_F$ with
	$A,C\in\mathscr O_F$, then the long exact Hom sequences show that
	vanishing of
	$\operatorname{Hom}_{\mathscr T_F}(JP,A[1])$ and
	$\operatorname{Hom}_{\mathscr T_F}(JP,C[1])$ implies vanishing of
	$\operatorname{Hom}_{\mathscr T_F}(JP,B[1])$. The dual assertion
	holds for $\operatorname{Hom}_{\mathscr T_F}(-,JP[1])$.
	By induction on the extension length, and using additivity to pass to
	finite direct sums and direct summands, we conclude that
	\[
	\mathbb E_{\mathscr O_F}(JP,-)=0 \quad\text{and}\quad \mathbb E_{\mathscr O_F}(-,JP)=0
	\]
	on all of $\mathscr O_F$. The same holds for every object of
	$\mathscr P_F$. Thus the objects of $\mathscr P_F$ are
	projective-injective.
	
	We next prove that there are enough projectives. Let $\mathscr C$ be
	the full subcategory of $\mathscr O_F$ formed by the objects $M$
	for which there exists a deflation
	\[ P_M\twoheadrightarrow M, \qquad P_M\in\mathscr P_F. \]
	For $X\in\mathscr E$, choose a conflation
	\[ K_X\longrightarrow P_X\longrightarrow X \]
	in $\mathscr E$, with $P_X\in\mathscr P$. Applying $J$ gives
	a conflation in $\mathscr O_F$, so $JX\in\mathscr C$.
	Clearly, $\mathscr C$ is closed under finite direct sums.
	
	The subcategory $\mathscr C$ is also closed under direct summands.
	Indeed, suppose $M\in\mathscr C$ and that $X$ is a direct summand
	of $M$. Choose a deflation $p\colon P\twoheadrightarrow M$, with
	$P\in\mathscr P_F$, and a retraction $r\colon M\twoheadrightarrow X$. A retraction is a split deflation, and a composite of deflations
	is a deflation. Hence $rp\colon P\twoheadrightarrow X$ is a
	deflation, so $X\in\mathscr C$.
	
	Finally, $\mathscr C$ is extension-closed. Consider an extriangle
	\[ A\longrightarrow B\xrightarrow{b}C\dashrightarrow \]
	with $A,C\in\mathscr C$, and choose deflations
	\[ p_A\colon P_A\twoheadrightarrow A, \qquad p_C\colon P_C\twoheadrightarrow C, \]
	with $P_A,P_C\in\mathscr P_F$. Pulling the given extriangle back
	along $p_C$ yields a morphism of extriangles
	\[ \begin{array}{cccccc} A&\longrightarrow&B'&\longrightarrow&P_C&\dashrightarrow\\
\Vert&&\downarrow u&&\downarrow p_C&\\
A&\longrightarrow&B&\xrightarrow{b}&C&\dashrightarrow. \end{array} \]
	Since $P_C$ is projective, the upper extriangle splits, and hence
	$B'\simeq A\oplus P_C$. Moreover, $u$ is a deflation since it is the
	base change of the deflation $p_C$ along $b$. Therefore
	\[ P_A\oplus P_C \xrightarrow{\,p_A\oplus1_{P_C}\,} A\oplus P_C\simeq B' \xrightarrow{\,u\,}B \]
	is a composite of deflations. Thus $B\in\mathscr C$.
	
	The subcategory $\mathscr C$ is extension-closed and
	idempotent-complete and contains $J(\mathscr E)$. By
	\Cref{lem:orbit-extension-closure}, it equals $\mathscr O_F$.
	Thus $\mathscr O_F$ has enough projectives in $\mathscr P_F$.
	The dual argument shows that it has enough injectives in
	$\mathscr P_F$.
	
	It remains to identify all projective and injective objects. If $Q$
	is projective, choose a deflation $P\twoheadrightarrow Q$, with
	$P\in\mathscr P_F$. Since $Q$ is projective, this deflation
	splits, so $Q$ is a direct summand of $P$. Hence
	$Q\in\mathscr P_F$. The dual argument applies to injective objects.
	This proves the assertion.
\end{proof}

\begin{Lem}
	\label{lem:frobenius-orbit-higher-orthogonality}
	For every $P\in\mathscr P_F$, $X\in\mathscr O_F$, and $r>0$,
	one has
	\[ \operatorname{Hom}_{\mathscr T_F}(P,X[r])=0 \qquad\text{and}\qquad \operatorname{Hom}_{\mathscr T_F}(X,P[r])=0. \]
\end{Lem}

\begin{proof}
	We first take $P=JP_0$ and $X=JX_0$, with
	$P_0\in\mathscr P$ and $X_0\in\mathscr E$. By
	\Cref{prop:orbit-hom-ext} and Chen's comparison with higher
	extensions \cite[Proposition~3.18]{Chen20262},
	\[ \begin{aligned} \operatorname{Hom}_{\mathscr T_F}(JP_0,JX_0[r]) &\cong \bigoplus_{n\in\mathbb Z} \mathbb E_{\mathcal E}^{r}(P_0,F^nX_0),\\
\operatorname{Hom}_{\mathscr T_F}(JX_0,JP_0[r]) &\cong \bigoplus_{n\in\mathbb Z} \mathbb E_{\mathcal E}^{r}(X_0,F^nP_0). \end{aligned} \]
	Here $\mathbb E_{\mathcal E}^{r}$ denotes the $r$-th higher
	extension bifunctor of $\mathscr E$. It is constructed as the
	$r$-fold tensor product of $\mathbb E_{\mathcal E}$ over
	$\mathscr E$. Since $P_0$ is projective,
	\[ \mathbb E_{\mathcal E}^{r}(P_0,-)=0 \qquad(r>0). \]
	Since $F^nP_0$ is injective for every $n\in\mathbb Z$,
	\[ \mathbb E_{\mathcal E}^{r}(-,F^nP_0)=0 \qquad(r>0). \]
	This proves both vanishings for the objects $JP_0$ and $JX_0$.
	
	Fix $P_0\in\mathscr P$ and $r>0$. If
	\[ A\longrightarrow B\longrightarrow C\longrightarrow A[1] \]
	is a triangle with $A,C\in\mathscr O_F$, the long exact Hom
	sequences show that vanishing for $A$ and $C$ implies vanishing
	for $B$, both for
	$\operatorname{Hom}_{\mathscr T_F}(JP_0,-[r])$ and for
	$\operatorname{Hom}_{\mathscr T_F}(-,JP_0[r])$. Induction on the
	extension length therefore extends both vanishings from
	$J(\mathscr E)$ to all of $\mathscr O_F$. Finally, the vanishings
	are preserved under finite direct sums and direct summands. Since
	$\mathscr P_F={\operatorname{add}} J(\mathscr P)$, the
	result follows.
\end{proof}

\begin{Thm}\label{thm:frobenius-orbit-realization}
{\color{black}The canonical exact dg functor
$\jmath_F\colon \operatorname{Orb}^{\mathrm{ex}}_{\mathrm{dg}}(\mathcal E,F) \longrightarrow\mathcal T_F$
induces a quasi-fully faithful realization morphism}
\[ \operatorname{real}_F\colon \cd^b_{\mathrm{dg}}\bigl(\operatorname{Orb}^{\mathrm{ex}}_{\mathrm{dg}}(\mathcal E,F)\bigr)\longrightarrow\mathcal T_F, \]
whose essential image on homotopy categories is $\operatorname{tria}_{\mathscr T_F}(\mathscr O_F)$. After dg idempotent completion, it induces a quasi-equivalence
\[ \left(\cd^b_{\mathrm{dg}}\bigl(\operatorname{Orb}^{\mathrm{ex}}_{\mathrm{dg}}(\mathcal E,F)\bigr)\right)^{\sharp}\xrightarrow{\sim}\mathcal T_F. \]
Without idempotent completion of the source, only the stated essential-image description is asserted.
\end{Thm}
{\color{black}
\begin{proof}
Let $X,Y\in\mathscr O_F$ and $f\colon X\to Y[r]$, with $r\geq2$. By \Cref{thm:frobenius-exact-orbit} choose a projective conflation, whose ambient triangle is
\begin{equation}\label{eq:orbit-realization-syzygy-triangle} \Omega X\longrightarrow P_X\xrightarrow{p_X}X\xrightarrow{\delta_X}\Omega X[1]. \end{equation}
\Cref{lem:frobenius-orbit-higher-orthogonality} gives $\Hom(P_X,Y[r])=0$. Hence $fp_X=0$, and the long exact Hom sequence gives $f=f'\delta_X$ for $f'\colon\Omega X[1]\to Y[r]$. Induction on $r$, starting with degree one, factors $f'[-1]$ into $r-1$ degree-one maps through objects of $\mathscr O_F$. Shifting and composing with $\delta_X$ proves that $\mathscr O_F$ is $1$-generated in $\mathscr T_F$.
The full extension-closed dg subcategory $\widetilde{\mathcal O}_F\subseteq\mathcal T_F$ therefore satisfies \cite[Theorem~3.2]{Plogmann26}. The realization is quasi-fully faithful and has essential image exactly $\operatorname{tria}_{\mathscr T_F}(\mathscr O_F)$.
By construction $\cd^b_{\mathrm{dg}}(\mathcal E)=\operatorname{pretr}(\mathcal E)/\mathcal N_{\mathrm{dg}}$. Its image is triangulatedly generated by the objects of $\mathcal E$. The dg orbit has the same objects before taking its hull, so $J(\mathscr E)$ thickly generates the completed category $\mathscr T_F$. Since $J(\mathscr E)\subseteq\mathscr O_F$, we have $\thick_{\mathscr T_F}(\mathscr O_F)=\mathscr T_F$. Completing the source of the quasi-fully faithful realization consequently makes it essentially surjective. This proves the last assertion, without identifying $\operatorname{tria}$ with $\thick$.
\end{proof}
}

\subsection{Comparison with the stable orbit category}

Let $\mathcal P_F$ be the full dg subcategory of
$\operatorname{Orb}^{\mathrm{ex}}_{\mathrm{dg}}(\mathcal E,F)$ on
the objects of $\mathscr P_F$, and put
\[ \underline{\mathscr O_F} := \mathscr O_F/[\mathscr P_F]. \]
The category $\underline{\mathscr O_F}$ is triangulated because
$\mathscr O_F$ is Frobenius.

Put
\[\mathcal N:=\operatorname{pretr}(\iota\mathcal P)\subseteq\mathcal D,\qquad \mathcal D_{\mathrm{st}}:=\mathcal D/\mathcal N,\]
and let
$q_{\mathrm{st}}\colon\mathcal D\to\mathcal D_{\mathrm{st}}$
be the dg quotient functor.
Since $F(\mathcal P)=\mathcal P$, both $\widehat F$ and
$\widehat F^{-1}$ preserve $\mathcal N$.
Write
\[\widehat F_{\mathcal N}:=\widehat F|_{\mathcal N},\qquad \widehat F_{\mathrm{st}}\colon\mathcal D_{\mathrm{st}}\longrightarrow\mathcal D_{\mathrm{st}}\]
for the restriction and the induced dg automorphism, respectively.
Using a functorial dg quotient model, these satisfy
\[q_{\mathrm{st}}\widehat F=\widehat F_{\mathrm{st}}q_{\mathrm{st}}.\]
Thus $\widehat F_{\mathrm{st}}$ is the action induced by
$\widehat F$ after quotienting out $\mathcal N$.

\begin{Thm}\label{thm:stable-orbit-comparison}
There is a canonical triangle equivalence
\[ \operatorname{Kar}(\underline{\mathscr O_F}) \xrightarrow{\sim} H^0\left(\left(\operatorname{pretr}(\mathcal D_{\mathrm{st}}/\widehat F_{\mathrm{st}}^{\mathbb Z})\right)^{\sharp}\right). \]
If $\underline{\mathscr O_F}$ is already idempotent complete, its Karoubi envelope can be omitted. In particular this applies to the Krull--Schmidt Higgs categories used in \Cref*{sec:weighted-totalization-singularity}.
\end{Thm}
\begin{proof}
Put $\mathcal A=\operatorname{Orb}^{\mathrm{ex}}_{\mathrm{dg}}(\mathcal E,F)$.

Since $\mathcal A$ is connective Frobenius,
\cite[Theorems~3.1 and~3.25]{Chen20262}
give the following triangle equivalence:

\begin{equation}\label{eq:stable-orbit-derived-quotient} \underline{\mathscr O_F}\simeq H^0\bigl(\cd^b_{\mathrm{dg}}(\mathcal A)/\operatorname{pretr}(\mathcal P_F)\bigr). \end{equation}
No idempotent completion is taken in
\Cref{eq:stable-orbit-derived-quotient}.
Here and below, $\operatorname{pretr}(\mathcal P_F)$ denotes the
full pretriangulated dg subcategory generated by the canonical
images of the objects of $\mathcal P_F$ in
$\cd^b_{\mathrm{dg}}(\mathcal A)$.

We use the following quotient--completion comparison.
Let $\mathcal C$ be a pretriangulated dg category and
$\mathcal B\subseteq\mathcal C$ a full pretriangulated dg
subcategory. The canonical maps induce a morphism
\[u\colon(\mathcal C/\mathcal B)^{\sharp}\longrightarrow\left(\mathcal C^{\sharp}/\operatorname{thick}_{\mathrm{dg}}(\mathcal B)\right)^{\sharp}\]
in $\mathrm{Hqe}$.
Here $\operatorname{thick}_{\mathrm{dg}}(\mathcal B)$ is the full
dg subcategory of $\mathcal C^{\sharp}$ whose homotopy category
is the thick closure of the image of $H^0(\mathcal B)$.

The morphism $u$ is an isomorphism in $\mathrm{Hqe}$.
Indeed, the dg quotients enhance the corresponding Verdier
quotients. For every idempotent-complete triangulated category
$\mathscr V$, restriction identifies exact functors from either
$H^0$-category in the comparison with exact functors
$H^0(\mathcal C)\to\mathscr V$ annihilating $H^0(\mathcal B)$.
This follows from the universal properties of Verdier quotients
and of the triangulated Karoubi envelope
\cite[Theorem~1.5]{BalmerSchlichting2001}.
Consequently, $H^0(u)$ is a triangle equivalence.
Since $u$ is a dg comparison between pretriangulated dg
categories, full faithfulness on all shifted Hom spaces gives
quasi-full faithfulness. Together with essential surjectivity
of $H^0(u)$, this proves
\[(\mathcal C/\mathcal B)^{\sharp}\simeq\left(\mathcal C^{\sharp}/\operatorname{thick}_{\mathrm{dg}}(\mathcal B)\right)^{\sharp}.\]

Apply the quotient--completion comparison with
$\mathcal C=\cd^b_{\mathrm{dg}}(\mathcal A)$ and
$\mathcal B=\operatorname{pretr}(\mathcal P_F)$.
By \Cref{thm:frobenius-orbit-realization}, the completed
realization is a quasi-equivalence
\[\operatorname{real}_F^{\sharp}\colon\left(\cd^b_{\mathrm{dg}}(\mathcal A)\right)^{\sharp}\xrightarrow{\sim}\mathcal T_F.\]
Its restriction along the canonical morphism from $\mathcal A$
agrees with $\jmath_F$ in $\mathrm{Hqe}$.
On homotopy categories, the images of the objects of
$\mathcal P_F$ therefore form
$\mathscr P_F=\operatorname{add}J(\mathscr P)$.
Hence $\operatorname{real}_F^{\sharp}$ identifies the full thick
dg subcategory generated by the canonical image of
$\mathcal P_F$ with
$\operatorname{thick}_{\mathrm{dg}}(J\mathcal P)$.
Combining this with
\Cref{eq:stable-orbit-derived-quotient} gives
\begin{equation}\label{eq:realization-stable-quotient} \operatorname{Kar}(\underline{\mathscr O_F})\simeq H^0\left(\left(\mathcal T_F/\operatorname{thick}_{\mathrm{dg}}(J\mathcal P)\right)^{\sharp}\right). \end{equation}

We now apply the orbit--quotient comparison to the
$\widehat F$-invariant inclusion $\mathcal N\subseteq\mathcal D$.
Put
\[\mathcal U_F:=\operatorname{pretr}(\mathcal D/\widehat F^{\mathbb Z}),\qquad \mathcal V_F:=\operatorname{pretr}(\mathcal N/\widehat F_{\mathcal N}^{\mathbb Z}).\]
The direct-sum formula for orbit morphism complexes, together
with fullness and $\widehat F^{\pm1}$-invariance of
$\mathcal N$, shows that the induced inclusion of dg orbit
categories is fully faithful in the chosen strict models.
The same holds after taking pretriangulated hulls.
We therefore regard $\mathcal V_F$ as a full pretriangulated
dg subcategory of $\mathcal U_F$.

Since $\mathcal D_{\mathrm{st}}=\mathcal D/\mathcal N$ and
$\widehat F_{\mathrm{st}}$ is the induced action,
\cite[Corollary~3.18]{FKQ2024} gives a canonical isomorphism
\[\operatorname{pretr}(\mathcal D_{\mathrm{st}}/\widehat F_{\mathrm{st}}^{\mathbb Z})\simeq\mathcal U_F/\mathcal V_F\]
in $\mathrm{Hqe}$.
By definition, $\mathcal U_F^{\sharp}=\mathcal T_F$.

We identify the thick closure of $\mathcal V_F$ in
$\mathcal T_F$.
Every object of $\mathcal N=\operatorname{pretr}(\iota\mathcal P)$
is obtained from $\iota\mathcal P$ by finite sums, shifts and
cones. Since $H^0(Q)$ is exact, its image belongs to
$\operatorname{tria}_{\mathscr T_F}(J\mathscr P)$.
Passing to the dg orbit does not change the objects, and taking
its pretriangulated hull only adds finite triangulated
constructions. Thus the image of $H^0(\mathcal V_F)$ in
$\mathscr T_F$ is contained in
$\operatorname{tria}_{\mathscr T_F}(J\mathscr P)$.
It also contains $J\mathscr P$.
Consequently, the full thick dg subcategories generated by
$\mathcal V_F$ and by $J\mathcal P$ coincide:
\[\operatorname{thick}_{\mathrm{dg}}(\mathcal V_F)=\operatorname{thick}_{\mathrm{dg}}(J\mathcal P),\]
where both closures are taken in $\mathcal T_F$.

Apply the quotient--completion comparison to
$\mathcal C=\mathcal U_F$ and $\mathcal B=\mathcal V_F$.
Together with the preceding orbit--quotient comparison, this
gives
\begin{equation}\label{eq:stable-orbit-to-quotient} \left(\operatorname{pretr}(\mathcal D_{\mathrm{st}}/\widehat F_{\mathrm{st}}^{\mathbb Z})\right)^{\sharp}\simeq\left(\mathcal T_F/\operatorname{thick}_{\mathrm{dg}}(J\mathcal P)\right)^{\sharp}. \end{equation}
All dg comparisons here are taken in $\mathrm{Hqe}$, with the
orbit comparison using the equivariant structure of the
quotient functor $q_{\mathrm{st}}$.
Taking $H^0$ in \Cref{eq:stable-orbit-to-quotient} and combining
it with \Cref{eq:realization-stable-quotient} proves the
canonical triangle equivalence.

If $\underline{\mathscr O_F}$ is idempotent complete, the
canonical functor to its Karoubi envelope is an equivalence.
In particular, if $\mathscr O_F$ is Krull--Schmidt, its stable
category is again Krull--Schmidt.
Indeed, every object has a finite decomposition into
indecomposables, and the projective--injective summands become
zero in the stable category.
For a non-projective indecomposable object $X$, the ring
\[\operatorname{End}_{\underline{\mathscr O_F}}(X)=\operatorname{End}_{\mathscr O_F}(X)/[\mathscr P_F](X,X)\]
is a quotient of a local ring by a proper ideal, and is
therefore local. The ideal is proper because a factorization
of $1_X$ through a projective object would make $X$ projective.
Thus the stable category is Krull--Schmidt and hence
idempotent complete. This applies to the Higgs categories
specified in the statement.
\end{proof}

\bibliographystyle{plain}

\begin{thebibliography}{}
	
	\bibitem{AmiotCluster}
	Claire Amiot,
	\emph{Cluster categories for algebras of global dimension 2
		and quivers with potential},
	Ann. Inst. Fourier (Grenoble) 59 (2009), no.~6, 2525--2590.
	
	
	\bibitem{AmiotOppermann}
	Claire Amiot and Steffen Oppermann, \emph{Cluster equivalence and
		graded derived equivalence}, Doc. Math. 19 (2014), 1155--1206.
	
	\bibitem{BalmerSchlichting2001}
	Paul Balmer and Marco Schlichting,
	\emph{Idempotent completion of triangulated categories},
	J. Algebra 236 (2001), no.~2, 819--834.
	
	\bibitem{BMRRT}
	Aslak Bakke Buan, Bethany Rose Marsh, Markus Reineke,
	Idun Reiten, and Gordana Todorov,
	\emph{Tilting theory and cluster combinatorics},
	Adv. Math. 204 (2006), no.~2, 572--618.
	
	\bibitem{BuchweitzMCM}
	R.-O. Buchweitz, \emph{Maximal Cohen--Macaulay modules and Tate
		cohomology}, Mathematical Surveys and Monographs, vol.~262,
	American Mathematical Society, Providence, RI, 2021.
	
	\bibitem{Chen20261}
	Xiaofa Chen, \emph{Exact dg categories I: Foundations},
	Adv. Math. 489 (2026), Paper No. 110809, 85 pp.
	
	\bibitem{Chen20262}
	Xiaofa Chen, \emph{Exact dg categories II: The embedding theorem},
	J. Pure Appl. Algebra 230 (2026), no. 7, Paper No. 108298, 56 pp.
	
	\bibitem{ChenChen2021}
	Xiaofa Chen and Xiao-Wu Chen,
	\emph{An informal introduction to dg categories},
	\href{https://arxiv.org/abs/1908.04599v2}{arXiv:1908.04599v2} [math.RT].
	
	\bibitem{FKQ2024}
	Li Fan, Bernhard Keller, Yu Qiu, \emph{Dg enhanced orbit categories and applications}, arXiv:2405.00093 [math.RT].
	
	\bibitem{FuKellerCoefficients}
	Changjian Fu and Bernhard Keller, \emph{On cluster algebras with
		coefficients and $2$-Calabi--Yau categories},
	Trans. Amer. Math. Soc. 362 (2010), no.~2, 859--895.
	
	\bibitem{Fushimi}
	Riku Fushimi,
	\emph{Contravariant Koszul duality between non-positive and positive dg algebras},
	\href{https://arxiv.org/abs/2409.08842v3}{arXiv:2409.08842v2}.
	
	\bibitem{GuoHigherCluster}
	Lingyan Guo, \emph{Cluster tilting objects in generalized higher
		cluster categories}, J. Pure Appl. Algebra 215 (2011), no.~9,
	2055--2071.
	
	\bibitem{Hanihara2022}
	Norihiro Hanihara, \emph{Cluster categories of formal DG algebras and singularity categories}, Forum of Mathematics, Sigma (2022), Vol. 10:e35 1–50.
	
	\bibitem{Hanihara24}
	Norihiro Hanihara, \emph{Calabi-Yau completions for roots of dualizing dg bimodules}, arXiv:2412.18753 [math.RT].
	

	\bibitem{IyamaTakahashiQuotient}
	O. Iyama and R. Takahashi, \emph{Tilting and cluster tilting for
		quotient singularities}, Math. Ann. \textbf{356} (2013), 1065--1105.

	\bibitem{IY2018}
	Osamu Iyama, Dong Yang, \emph{Silting reduction and Calabi–Yau reduction of triangulated categories},
	Trans. Amer. Math. Soc. 370 (2018), no. 11, 7861–7898.
	
	\bibitem{IY2020}
	Osamu Iyama, Dong Yang, \emph{Quotients of triangulated categories and equivalences of Buchweitz, Orlov, and Amiot-Guo-Keller}, Amer. J. Math. 142 (2020), no. 5, 1641–1659.
	
	\bibitem{Kalck-Yang2016}
	Martin Kalck, Dong Yang,
	\emph{Relative singularity categories I: Auslander resolutions},
	Adv. Math. 301 (2016), 973–1021.
	
	\bibitem{Kalck-Yang2018}
	Martin Kalck, Dong Yang, \emph{Relative singularity categories II: DG models}, arXiv:1803.08192 [math.AG].
	
	\bibitem{KellerOrbit}
	Bernhard Keller,
	\emph{On triangulated orbit categories},
	Doc. Math. 10 (2005), 551--581.
	
	\bibitem{KellerCYCompletion}
	Bernhard Keller,
	\emph{Deformed Calabi--Yau completions},
	with an appendix by Michel Van den Bergh,
	J. Reine Angew. Math. 654 (2011), 125--180.
	
	\bibitem{KellerScherotzke2016}
	Bernhard Keller and Sarah Scherotzke, \emph{Graded quiver varieties and derived categories}, arXiv:1303.2318 [math.RT].
	
	\bibitem{KellerWang2023}
	Bernhard Keller and Yu Wang, \emph{An introduction to relative Calabi--Yau structures}, Representations of algebras and related structures, EMS Ser. Congr. Rep., EMS Press, Berlin, 2023, 279--304; arXiv:2111.10771.
	
	\bibitem{KellerWuHiggs}
	Bernhard Keller and Yilin Wu, \emph{Relative cluster categories and
		Higgs categories with infinite-dimensional morphism spaces},
	with an appendix by Chris Fraser and Bernhard Keller,
	\href{https://arxiv.org/abs/2307.12279}{arXiv:2307.12279}.
	
		\bibitem{LiuHigherMcKay}
	J. Liu, \emph{Singularity categories via higher McKay quivers with
		potential}, Selecta Math. (N.S.) \textbf{32} (2026), Paper No.~54.
	
	\bibitem{NP2019}
	Hiroyuki Nakaoka, Yann Palu, \emph{Extriangulated categories, Hovey twin cotorsion pairs and model
	structures}, Cah. Topol. G\'{e}om. Diff\'{e}r. Cat\'{e}g. 60 (2019), no. 2, 117–193.
	
	\bibitem{Neeman2001}
	Amnon Neeman, \emph{Triangulated categories}, Annals of Mathematics Studies, vol. 148, Princeton University Press, 2001.
	
	\bibitem{Plogmann26}
	Marvin Plogmann, \emph{Complicial simple-minded collections}, arXiv:2603.03122 [math.RT].
	
	\bibitem{Wu2023}
	Yilin Wu, \emph{Relative cluster categories and Higgs categories}, Adv. Math. 424 (2023), Paper No. 109040, 112 pp.
	
	\bibitem{Wu2025}
	Yilin Wu, \emph{Silting reduction, relative AGK's construction and Higgs construction}, arXiv:2510.00470 [math.RT].
	

\end{thebibliography}

\end{document}